\documentclass[11pt]{article}

\usepackage{amsmath,amssymb}
\usepackage{pstricks,pst-node,xy}
\xyoption{all}

\numberwithin{equation}{section}
\newtheorem{thm}{Theorem}[section] 
\newtheorem{prp}[thm]{Proposition}
\newtheorem{lmm}[thm]{Lemma}   

\newtheorem{dfn}[thm]{Definition}
\newtheorem{mythm}{Theorem}

\def\lan{\langle}
\def\ran{\rangle}
\def\lr#1{\lan#1\ran}
\def\blr#1{\big\lan#1\big\ran}
\def\ov#1{\overline{#1}}
\def\ti#1{\tilde{#1}}
\def\wt#1{\widetilde{#1}}
\def\e_ref#1{(\ref{#1})}
\def\smsize#1{\begin{small}#1\end{small}}
\def\lgsize#1{\begin{large}#1\end{large}}
\def\sf#1{\textsf{#1}}

\def\lra{\longrightarrow}
\def\Lra{\Longrightarrow}
\def\Llra{\Longleftrightarrow}

\def\al{\alpha}
\def\be{\beta}
\def\ga{\gamma}
\def\de{\delta}
\def\ka{\kappa}
\def\la{\lambda}
\def\si{\sigma}
\def\th{\theta}
\def\vp{\varpi}

\def\De{\Delta}
\def\Ga{\Gamma}
\def\Th{\Theta}

\def\eset{\emptyset}
\def\i{\infty}
\def\hb{\hbar}

\def\cA{\mathcal A}
\def\cB{\mathcal B}
\def\C{\mathbb C}
\def\cC{\mathcal C}
\def\d{\mathfrak d}

\def\cD{\mathcal D}
\def\E{\mathbf e}
\def\F{\mathcal F}

\def\H{\mathcal H}
\def\I{\mathfrak i}
\def\cI{\mathcal I}
\def\J{\mathcal J}
\def\K{\mathcal K}
\def\L{\mathbb L}
\def\cL{\mathcal L}
\def\cM{\mathcal M}
\def\M{\mathfrak M}
\def\N{\mathcal N}
\def\P{\mathbb P}
\def\Pf{\mathbb P^4}
\def\Pn{\mathbb P^n}

\def\fR{\mathfrak R}
\def\O{\mathcal O}
\def\Q{\mathbb Q}
\def\T{\mathbb T}
\def\U{\mathfrak U}
\def\V{\mathcal V}
\def\cY{\mathcal Y}
\def\Z{\mathbb Z}
\def\cZ{\mathcal Z}

\def\Aut{\textnormal{Aut}}

\def\Dom{\textnormal{Dom}}
\def\Edg{\textnormal{Edg}}
\def\eff{\textnormal{eff}}
\def\ev{\textnormal{ev}}
\def\GW{\textnormal{GW}}
\def\mod{\textnormal{mod~}}

\def\Sym{\textnormal{Sym}}
\def\Span{\textnormal{Span}}
\def\val{\textnormal{val}}
\def\Ver{\textnormal{Ver}}

\begin{document}

\thispagestyle{empty}

\title{The Reduced Genus One Gromov-Witten Invariants\\
of Calabi-Yau Hypersurfaces}
\author{Aleksey Zinger\thanks{Partially supported by a Sloan fellowship 
and DMS Grant 0604874}}
\date{\today}
\maketitle

\begin{abstract}
\noindent
We compute the reduced genus~$1$ Gromov-Witten invariants of Calabi-Yau hypersurfaces.
As a consequence, we confirm the 1993 Bershadsky-Cecotti-Ooguri-Vafa (BCOV) prediction 
for the standard genus~$1$ GW-invariants of a quintic threefold.
We combine constructions from a series of previous papers with 
the classical localization theorem to relate the reduced genus~$1$ invariants 
of a CY-hypersurface to previously computed integrals on moduli spaces of stable 
genus~$0$ maps into projective space.
The resulting, rather unwieldy, expressions for a genus~$1$ equivariant
generating function simplify drastically, using a regularity property of 
a genus~$0$ equivariant generating function in half of the cases. 
Finally, by disregarding terms that cannot effect the non-equivariant part
of the former, we relate the answer to an explicit hypergeometric 
series in a simple way. 
The approach described in this paper is systematic.
It is directly applicable to computing reduced genus~$1$ GW-invariants
of other complete intersections and should apply to 
higher-genus localization computations. 
\end{abstract}

\tableofcontents

\setcounter{section}{-1}

\section{Introduction}
\label{intro_sec}

\subsection{Mirror Symmetry Predictions for a Quintic Threefold}
\label{mirsym_subs}

\noindent
Gromov-Witten invariants of a smooth projective variety~$X$
are certain counts of curves in~$X$.
In many cases, these invariants are known or conjectured to possess
rather amazing structure which is often completely unexpected from
the classical point of view.
For example, a generating function for the genus~$0$ GW-invariants
solves a third-oder PDE in two variables.
In the case of the complex projective space~$\Pn$,
the resulting PDE condition on the generating function is equivalent to a recursion
for counts of rational curves in $\Pn$ and solves
the classical problem of enumeration of such curves in~$\Pn$;
see  \cite[Section~10]{RT} and \cite[Section~5]{KoM}.\\

\noindent
The above mentioned PDE property\footnote{It is equivalent to 
the associativity of the multiplication in quantum cohomology.}
is just one type of structure of GW-invariants motivated by their relation
to string theory.
The mirror symmetry principle of string theory predicts yet another type of structure
whenever $X$ is a Calabi-Yau threefold.
It relates GW-invariants of $X$ to an integral on
the moduli space of Kahler structures of the ``mirror'' of~$X$.
In the case $X$ is a quintic threefold (degree~$5$ hypersurface in~$\Pf$),
this integral was computed in~\cite{CaDGP}, leading to stunning predictions
concerning the Gromov-Witten theory of~$X$.
These predictions have been only partially verified.\\

\noindent
For each pair of nonnegative integers $(g,d)$, let $N_{g,d}$ denote the genus~$g$
degree~$d$ GW-invariant of a quintic threefold~$X_5$; 
see equation~\e_ref{CYGW_e} below.
For $q\!=\!0,1,\ldots$, we define a degree~$q$ polynomial $I_q(t)$ in $t$
with coefficients in the power series in~$e^t$ by 
\begin{equation}\label{Ifuncdfn5_e}
\sum_{q=0}^{\i} I_q(t)w^q \equiv e^{wt}\sum_{d=0}^{\i}e^{dt}
\frac{\prod_{r=1}^{r=5d}(5w\!+\!r)}{\prod_{r=1}^{r=d}((w\!+\!r)^5\!-\!w^5)}.
\end{equation}
For example,
\begin{equation}\label{Ifuncdfn5_e2}
I_0(t)=1+\sum_{d=1}^{\i}e^{dt}\frac{(5d)!}{(d!)^5}\,, \qquad
I_1(t)=tI_0(t)+\sum_{d=1}^{\i}e^{dt}
\bigg(\frac{(5d)!}{(d!)^5}\sum_{r=d+1}^{5d}\!\frac{5}{r}\bigg).
\end{equation}
Let
\begin{equation}\label{Jfuncdfn5_e}
J_q(t)=I_q(t)\big/I_0(t) ~~~\forall\,q=1,2,\ldots, \qquad T=J_1(t).
\end{equation}
The mirror symmetry prediction of~\cite{CaDGP} for the genus~$0$ GW-invariants 
of~$X_5$ can be stated~as
\begin{equation}\label{g0ms_e}
\frac{5}{6}T^3+\sum_{d=1}^{\i}N_{0,d}e^{dT}
=\frac{5}{2}\big(J_1(t)J_2(t)-J_3(t)\big);
\end{equation}
see Appendix~\ref{comparison_app} for a comparison of statements of mirror symmetry.
A prediction for the genus~$1$ GW-invariants of $X_5$ was made in~\cite{BCOV},
building up on~\cite{CaDGP}. 
Both of these predictions date back to the early days of the Gromov-Witten theory.
More recently, predictions for higher-genus GW-invariants of~$X_5$
have been made; the approach of~\cite{HKlQ} generates mirror formulas for
GW-invariants of~$X_5$ up to genus~$51$.\\

\noindent
While the ODE condition on GW-invariants mentioned above is proved directly,
the mathematical approach to the mirror principle has been to compute 
the relevant GW-invariants in each specific case.
However, this is rarely a simple task.
The prediction for genus~$0$ invariants was confirmed mathematically 
in the mid-1990s.
The prediction for genus~$1$ invariants is verified in this paper.

\begin{mythm}
\label{bcov_thm}
If $N_{1,d}$ denotes the degree $d$ genus $1$ Gromov-Witten invariant of
a quintic threefold,
\begin{equation}\label{g1ms_e}
2\sum_{d=1}^{\i}N_{1,d}e^{dT}
=\frac{25}{6}\big(J_1(t)-t\big)+
\ln\bigg(I_0(t)^{-62/3}\big(1\!-\!5^5e^t\big)^{-1/6}J_1'(t)^{-1}\bigg).
\end{equation}\\
\end{mythm}

\noindent
This theorem is deduced from Theorem~\ref{main_thm0} in Subsection~\ref{lowdeg_subs}.
An outline of this paper is contained in the next subsection.\\

\noindent
The author is very grateful to D.~Zagier,
for taking interest in the statements of and collaborating on the proof of
Propositions~\ref{hypergeom_prp} and~\ref{hypergeom_prp2}, 
to Jun Li, for bringing~\cite{BCOV} to the author's attention and 
collaborating on related work concerning genus~$1$ GW-invariants, 
to R.~Pandharipande, for generously sharing his expertise in the Gromov-Witten theory,
and the referee, for suggestions on improving the exposition.
The author would also like to thank S.~Katz, Y.-P.~Lee, D.~Maulik, 
A.~Mustata, A.~Postnikov, and R.~Vakil for many enlightening discussions.

\subsection{Computing GW-Invariants of Hypersurfaces}
\label{appr_subs}

\noindent
One approach to computing GW-invariants of a projective hypersurface
(and more generally, of a complete intersection) is to relate them to
GW-invariants of the ambient projective space as follows.
Whenever $g$, $d$, and $k$ are nonnegative integers and $X$ is a smooth subvariety of~$\Pn$, 
denote by $\ov\M_{g,k}(X,d)$ the moduli space of stable degree-$d$ maps into $X$
from genus~$g$ curves with $k$ marked points; see \cite[Chapter~24]{MirSym}.
Let $\U$ be the universal curve over $\ov\M_{g,k}(\Pn,d)$,
with structure map~$\pi$ and evaluation map~$\ev$:
$$\xymatrix{\U \ar[d]^{\pi} \ar[r]^{\ev} & \Pn \\
\ov\M_{g,k}(\Pn,d).}$$
In other words, the fiber of $\pi$ over $[\cC,f]\!\in\!\ov\M_{g,k}(\Pn,d)$,
where $\cC$ is a nodal curve with $k$ marked points and $f\!:\cC\!\lra\!\Pn$
is a stable morphism, is~$\cC/\Aut(\cC,f)$, while 
$$\ev\big([\cC,f;z]\big)=f(z) \qquad\hbox{if}\quad 
z\!\in\!\cC.$$\\

\noindent
A smooth degree-$a$ hypersurface $X$ in $\P^n$ is determined by a section $s$ of $\O_{\Pn}(a)$ 
which is transverse to the zero set:
$$X=s^{-1}(0) \qquad\hbox{for some}~~ 
s\in H^0\big(\Pn;\O_{\Pn}(a)\big).\footnote{In other words, 
$s$ is a holomorphic section of $\ga_n^{*\otimes a}$, 
where $\ga_n\!\lra\!\Pn$ is the tautological line bundle.}$$
The section $s$ induces a section $\ti{s}$ of the sheaf 
$\pi_*\ev^*\O_{\Pn}(a)\!\lra\!\ov\M_{g,k}(\Pn,d)$ by
$$\ti{s}\big([\cC,f]\big)=[s\circ f].$$
It is immediate that 
\begin{equation}\label{hyperpl_e1}
\ov\M_{g,k}(X,d) \equiv\big\{[\cC,f]\!\in\!\ov\M_{g,k}(\Pn,d)\!: f(\cC)\!\subset\!X\big\}
=\ti{s}^{-1}(0).
\end{equation}
On the other hand, GW-invariants of $X$ are defined by integration against
the virtual fundamental class (VFC) of $\ov\M_{g,k}(X,d)$ constructed in 
\cite{BeFa}, \cite{FuO}, and~\cite{LiT}:
\begin{equation}\label{GW_e}
\GW_{g,k}^X(d;\eta)=
\blr{\eta,\big[\ov\M_{g,k}(X,d)\big]^{vir}}
\qquad \forall~\eta\!\in\!H^*\big(\ov\M_{g,k}(X,d);\Q\big).
\end{equation}
If $X$ is a quintic threefold (or another Calabi-Yau threefold), the 
cycle $[\ov\M_{g,0}(X,d)]^{vir}$ is zero-dimensional and 
its degree is denoted by~$N_{g,d}$:
\begin{equation}\label{CYGW_e}
N_{g,d}= \GW_{g,0}^X(d;1)\equiv\blr{1,\big[\ov\M_{g,0}(X,d)\big]^{vir}}.
\end{equation}\\

\noindent
In light of Poincare Duality, equations~\e_ref{hyperpl_e1} and~\e_ref{GW_e}
suggest that $\GW_{g,k}^X(d;\eta)$ should be expressible as an integral
against $[\ov\M_{g,k}(\Pn,d)]^{vir}$ via some sort of euler class of 
the sheaf $\pi_*\ev^*\O_{\Pn}(a)$, whenever $\eta$ comes from $\ov\M_{g,k}(\Pn,d)$.
As can be easily seen from the definition of VFC,
this is indeed the case if $g\!=\!0$:
\begin{equation}\label{genus0_e}
\GW_{0,k}^X(d;\eta)
=\blr{\eta\cdot e\big(\pi_*\ev^* \O_{\Pn}(a)\big), \big[\ov\M_{0,k}(\Pn,d)\big]}
\qquad\forall~\eta\!\in\!H^*(\ov\M_{0,k}(\Pn,d);\Q).
\end{equation}
The moduli space $\ov\M_{0,k}(\Pn,d)$ is a smooth stack (orbifold)  and
$$\pi_*\ev^* \O_{\Pn}(a) \lra\ov\M_{0,k}(\Pn,d)$$
is a locally free sheaf (vector bundle).
Thus, the right-hand side of~\e_ref{genus0_e}
can be computed via the classical localization theorem of~\cite{ABo},
though the complexity of this computation increases rapidly with the degree~$d$.
Nevertheless, it has been completed in full generality in a number of different of ways,
verifying the genus~$0$ mirror symmetry prediction~\e_ref{g0ms_e}.
At the end of this subsection we briefly recall 
Givental's approach~\cite{Gi} to dealing with this complexity and describe 
its interplay with our approach in genus~$1$;
other proofs of~\e_ref{g0ms_e} can be found in \cite{Ber}, \cite{Ga}, 
\cite{Le}, and~\cite{LLY}.\\

\noindent
The most algebraically natural generalization of~\e_ref{genus0_e} to positive genera
fails.
In the genus~$1$ case, by \cite[Theorem~1.1]{LiZ}, the most 
topologically natural analogue of~\e_ref{genus0_e} does hold for the {\it reduced} 
GW-invariants $\GW_{1,k}^{0;X}$ defined in~\cite{g1comp2}:
\begin{equation}\label{genus1_e1}
\GW_{1,k}^{0;X}(d;\eta)
=\blr{\eta\cdot e\big(\pi_*\ev^* \O_{\Pn}(a)\big), \big[\ov\M_{1,k}^0(\Pn,d)\big]}
\qquad\forall~\eta\!\in\!H^*(\ov\M_{1,k}(\Pn,d);\Q),
\end{equation}
where $\ov\M_{1,k}^0(\Pn,d)$ is the \sf{main component} of $\ov\M_{1,k}(\Pn,d)$,
i.e.~the closure of the locus in $\ov\M_{1,k}(\Pn,d)$ consisting of maps from
smooth domains.
While $\ov\M_{1,k}^0(\Pn,d)$ is not a smooth stack, it is an equi-dimensional
orbi-variety and has a well-defined fundamental class.
While the sheaf
$$\pi_*\ev^* \O_{\Pn}(a) \lra \ov\M_{1,k}^0(\Pn,d)$$
is not locally free, it is shown in~\cite{g1cone} that its euler class 
is well-defined.
If 
\begin{equation}\label{desing_e}
\xymatrix{\V_1 \ar[d] \ar[r]^{p_*}& \pi_*\ev^*\O_{\Pn}(a) \ar[d]\\ 
\wt\M_{1,k}^0(\Pn,d)\ar[r]^p &  \ov\M_{1,k}^0(\Pn,d)}
\end{equation}
is a desingularization of $\ov\M_{1,k}^0(\Pn,d)$ and $\pi_*\ev^*\O_{\Pn}(a)$
(i.e.~$\wt\M_{1,k}^0(\Pn,d)$ is smooth, $p$ is a birational morphism, 
and $\V_1$ is locally free), then
\begin{equation}\label{genus1_e2}
\blr{\eta\cdot e\big(\pi_*\ev^* \O_{\Pn}(a)\big), \big[\ov\M_{1,k}^0(\Pn,d)\big]}
=\blr{\pi^*\eta\cdot e(\V_1), \big[\wt\M_{1,k}^0(\Pn,d)\big]}
\end{equation}
for all $\eta\!\in\!H^*(\ov\M_{1,k}^0(\Pn,d);\Q)$.
A natural desingularization~\e_ref{desing_e} that inherits every torus action from $\Pn$ 
is constructed in~\cite{VaZ}.
Thus, the classical localization theorem of~\cite{ABo} can be used to compute 
the right-hand side of~\e_ref{genus1_e1} via~\e_ref{genus1_e2}.
On the other hand, by \cite[Section~3]{g1comp2},
the reduced genus~$1$ GW-invariant differs from the standard one by a 
combination of genus~$0$ GW-invariants.
In particular, if $X$ is a quintic threefold, by \cite[Theorem~1.1]{g1comp2}
\begin{equation}\label{genus1_e3}
N_{1,d}\equiv \GW_{1,0}^X(d;1)
=N_{1,d}^0+\frac{1}{12}N_{0,d},
\end{equation}
where $N_{1,d}^0$ is the reduced genus~$1$ degree~$d$ 
invariant $\GW_{1,0}^{0;X}(d;1)$ of~$X$.\\

\noindent
In this paper we compute the numbers $\GW_{1,0}^{0;X}(d;1)$
for a smooth degree~$n$ hypersurface $X$ in~$\P^{n-1}$, with $n\!\ge\!3$.
By~\e_ref{genus1_e1}, \e_ref{genus1_e2}, and the divisor 
relation (see~\cite[Section 26.3]{MirSym}),
\begin{equation}\label{derivgen_e1}
d\,\GW_{1,0}^{0;X}(d;1)=\blr{e(\V_1)\,\ev_1^*H,\big[\wt\M_{1,1}^0(\P^{n-1},d)\big]},
\end{equation}
where $H\!\in\!H^2(\P^{n-1})$ is the hyperplane class and 
$$\V_1\lra \wt\M_{1,1}^0(\P^{n-1},d)$$
is the desingularization of the sheaf 
$$\pi_*\ev^* \O_{\P^{n-1}}(n) \lra \ov\M_{1,k}^0(\P^{n-1},d)$$
constructed in~\cite{VaZ}.
Analogously to \cite{Gi}, we package all reduced genus~$1$ GW-invariants~\e_ref{derivgen_e1} 
of a degree~$n$ hypersurface~$X_n$ in~$\P^{n-1}$ into a generating function~$\F$;
this is a power series with coefficients in the equivariant cohomology of~$\P^{n-1}$.
In~\cite{Gi}, an analogous power series provides a convenient way 
to describe a degree-recursive feature of the genus~$0$ GW-invariants of~$X_n$.
The resulting recursion, \cite[Lemma~30.1.1]{MirSym}, has 
a ``mystery correction term", which is intrinsically determined
by the recursion and another property of the genus~$0$ generating function
called polynomiality (see \cite[Section~30.2]{MirSym}).
By applying the Atiyah-Bott localization theorem, 
we relate~$\F$ to the genus~$0$ generating functions \e_ref{Z1ptdfn_e}-\e_ref{Z2ptdfn_e}.
In the process, we encounter seemingly unwieldy terms which turn out
to be  similar to the expressions encoded by
the genus~$0$ ``correction term" cleverly avoided in~\cite{Gi};
the latter are all of the form~\e_ref{regprp_e1}.
While none of these expressions by itself determines any of our unwieldy terms,
the entire series of expressions insures that the relevant genus~$0$ 
generating function has the remarkably rigid structure of Definition~\ref{regular_dfn}, 
which in turn determines all of our unwieldy terms via~\e_ref{regprp_e2}.
This leads to Propositions~\ref{contr_prp1} and~\ref{contr_prp2},
which describe the contributions to~$\F$ 
from the two different types of fixed loci in terms of known integrals on 
$\ov\M_{0,2}(\P^{n-1},d)$.
In Section~\ref{alg_sec}, we use Lemma~\ref{algebr_lmm1} to extract 
the non-equivariant part of the expressions in Propositions~\ref{contr_prp1} 
and~\ref{contr_prp2}, obtaining Theorem~\ref{main_thm} on page~\pageref{main_thm}.
Theorem~\ref{main_thm0} in the next subsection
follows immediately from Theorem~\ref{main_thm} and~\e_ref{derivgen_e2}.\\

\noindent
The approach of this paper to summing over all possible fixed loci involves
breaking the graph into trees at a special node.
As such trees contribute to certain genus~$0$ integrals, 
the desired sum is expressible in terms of these integrals.
The same approach directly carries over to computing reduced genus~$1$ GW-invariants
of any complete intersection and 
should be applicable to localization computations in higher 
genus.\footnote{This is not to say that higher-genus GW-invariants of projective 
hypersurfaces are now easily computable.
It is far from clear at this point what integrals should be localized in higher genus.
No higher-genus analogue of~\e_ref{genus1_e1} has been proved yet, though 
a conjectural version is stated in \cite[Subsection~1.1]{LiZ}.
Even with such a higher-genus hyperplane property, one would still need to either 
figure out how to apply the localization theorem in a singular setting or
construct a desingularization of the main component of $\ov\M_{g,k}(\Pn,d)$.}
In the latter case, there will be more ``special'' nodes, but their number will
be bounded above by the genus. 
Once the graphs are broken at the special nodes, 
there will be a number of distinguished trees and an arbitrary number of ``generic'' trees.
The number of the former will again be bounded by the genus.
On the other hand, it should be possible to sum over all possibilities
for the latter, using the regularity property of the relevant genus~$0$ integral
described in Subsection~\ref{Zreg_subs}.

\subsection{Mirror Symmetry Formulas for Projective CY-Hypersurfaces}
\label{lowdeg_subs}

\noindent
In this subsection we formulate a generalization of Theorem~\ref{bcov_thm}
to projective Calabi-Yau hypersurfaces of arbitrary dimensions;
see Theorem~\ref{main_thm0} below.
We then take a closer look at its low-dimensional cases, comparing some of them
with known results and others with the mirror symmetry predictions
of~\cite{BCOV} and~\cite{KlPa}.\\

\noindent
Let $n$ be a positive integer.
For each $q\!=\!0,1,\ldots$, define $I_{0,q}(t)$ by 
\begin{equation}\label{Ifuncdfn_e}
\sum_{q=0}^{\i} I_{0,q}(t)w^q \equiv e^{wt}\sum_{d=0}^{\i}e^{dt}
\frac{\prod_{r=1}^{r=nd}(nw\!+\!r)}{\prod_{r=1}^{r=d}((w\!+\!r)^n\!-\!w^n)}
\equiv R(w,t).
\end{equation}
Each $I_{0,q}(t)$ is a degree-$q$ polynomial in $t$ 
with coefficients that are power series in~$e^t$; 
see~\e_ref{Ifuncdfn5_e2} for explicit formulas for $I_0\!\equiv\!I_{0,0}$ 
and $I_1\!\equiv\!I_{0,1}$ in the $n\!=\!5$ case.
For $p,q\!\in\!\Z^+$ with $q\!\ge\!p$, let
\begin{equation}\label{Tden_e}
I_{p,q}(t)=\frac{d}{dt}\bigg(\frac{I_{p-1,q}(t)}{I_{p-1,p-1}(t)}\bigg).  
\end{equation}
By the first statement of Proposition~\ref{hypergeom_prp} below, 
each of the ``diagonal'' terms $I_{p,p}(t)$ is a power series in $e^t$ with constant 
term~$1$, whenever it is defined. 
Thus, the division in~\e_ref{Tden_e} is well-defined for all~$p$.
Let
\begin{equation}\label{mirmap_e}
T=\frac{I_{0,1}(t)}{I_{0,0}(t)}.
\end{equation}
By (i) of Proposition~\ref{hypergeom_prp}, the map $t\!\lra\!T$ is a change of variables;
it will be called the \sf{mirror map}.\\

\noindent
Let $\bar{R}(w,t)\!=\!R(w,t)/I_{0,0}(t)$. 
Then, $e^{-wt}\bar{R}(w,t)$ is a power series with $e^t$-constant term~$1$ and
$$\cD_w^p\ln\bar{R}(w,t)\equiv  \frac{1}{p!}\bigg\{\frac{d}{dw}\bigg\}^p
\Big(\ln\big(e^{-wt}\bar{R}(w,t)\big)\Big) \bigg|_{w=0} 
\in\Q[[e^t]]$$
for all $p\!\in\!\Z^+$.

\begin{mythm}\label{main_thm0}
For each $n\!\in\!\Z^+$,
the reduced genus~$1$ degree~$d$ Gromov-Witten invariants of a degree~$n$ 
hypersurface~$X$ in $\P^{n-1}$ are given~by
\begin{equation*}\begin{split}
\sum_{d=1}^{\i}e^{dT}\GW_{1,0}^{0;X}(d;1)&=
\bigg(\frac{(n\!-\!2)(n\!+\!1)}{48}+\frac{1-(1\!-\!n)^n}{24n^2}\bigg)(T\!-\!t)
+\frac{n^2\!-\!1+(1\!-\!n)^n}{24n}\ln I_{0,0}(t)\\
&\qquad -\begin{cases}
\frac{n-1}{48}\ln\big(1\!-\!n^ne^t\big)+
\sum_{p=0}^{(n-3)/2}\frac{(n-1-2p)^2}{8}\ln I_{p,p}(t),
&\hbox{if}~2\!\not|n;\\
\frac{n-4}{48}\ln\big(1\!-\!n^ne^t\big)+
\sum_{p=0}^{(n-4)/2}\frac{(n-2p)(n-2-2p)}{8}\ln I_{p,p}(t),
&\hbox{if}~2|n;
\end{cases}\\
&\qquad+
\frac{n}{24}\sum_{p=2}^{n-2}\bigg(\cD_w^{n-2-p}\frac{(1\!+\!w)^n}{(1\!+\!nw)}\bigg)
\big(\cD_w^p\ln\bar{R}(w,t)\big),
\end{split}\end{equation*}
where $t$ and $T$ are related by the mirror map~\e_ref{mirmap_e}.
\end{mythm}

\noindent
If $n\!=\!1$, both sides of the formula in Theorem~\ref{main_thm0} vanish.
If $n\!=\!2$, $X$ is a pair of points in $\P^1$.
In this case, the right-hand side of the formula in Theorem~\ref{main_thm0}
vanishes by~\e_ref{hypergeom_prp_e2}.
This is exactly as one would expect, since there are no positive-degree maps 
from a curve to a point.\\

\noindent
If $n\!=\!3$, $X$ is a plane cubic, i.e.~a $2$-torus embedded as a degree-$3$ curve in $\P^2$.
Thus, its degree~$d$ GW-invariant is zero unless $d$ is divisible by~$3$.
Furthermore, its genus~$1$ degree~$3r$ GW-invariant is the number of $r$-fold (unramified)
covers of a torus by a torus divided by $r$, the order of the automorphism group of
each such cover.
Since the number~$\si_r$ of such covers is given by~\e_ref{toruscovers_e}, it follows that 
\begin{equation}\label{n3concl_e1}
\sum_{d=1}^{\i}e^{dT}\GW_{1,0}^X(d;1)
= -\sum_{d=1}^{\i}\ln\big(1\!-\!e^{3dT}\big).
\end{equation}
On the other hand, Theorem~\ref{main_thm0} gives
\begin{equation}\label{n3concl_e2}
\sum_{d=1}^{\i}e^{dT}\GW_{1,0}^{0;X}(d;1)
=\frac{1}{8}(T\!-\!t)-\frac{1}{24}\ln\big(1\!-\!27e^t)-\frac{1}{2}\ln I_0(t),
\end{equation}
where $I_0(t)$ and $T$ are given~by:
$$I_0(t)=1+\sum_{d=1}^{\i}e^{dt}\frac{(3d)!}{(d!)^3}\,, \qquad
T=t+\frac{1}{I_0(t)}\sum_{d=1}^{\i}e^{dt}
\bigg(\frac{(3d)!}{(d!)^3}\sum_{r=d+1}^{3d}\!\frac{3}{r}\bigg).$$
Since the standard and reduced (genus~$1$) invariants of a (complex) curves
are the same if no descendant classes are involved, 
\begin{equation}\label{n3concl_e3}
\frac{1}{8}(T\!-\!t)-\frac{1}{24}\ln\big(1\!-\!27e^t)-\frac{1}{2}\ln I_0(t)
=-\sum_{d=1}^{\i}\ln\big(1\!-\!e^{3dT}\big)
\end{equation}
by \e_ref{n3concl_e1} and~\e_ref{n3concl_e2}.
We do not see a direct proof of \e_ref{n3concl_e3} at this point.\\

\noindent
If $n\!=\!4$, $X$ is a quartic surface in $\P^3$, i.e.~a $K3$.
All its GW-invariants are known to be zero.
With $n\!=\!4$, we find that the two coefficients of $(T\!-\!t)$ in Theorem~\ref{main_thm0}
add up to zero; the same is the case for the two coefficients of $\ln I_{0,0}(t)$.
Thus, the sum of the terms on the first two lines of the right-hand side in the formula
of Theorem~\ref{main_thm0} is zero.
The remaining term~is
$$\frac{4}{24}\cdot\frac{1}{2!}\bigg\{\frac{d}{dw}\bigg\}^2
\bigg(\ln\sum_{q=0}^{\i}J_q(t)w^q\bigg)\bigg|_{w=0}
=\frac{1}{6}\bigg(J_2(t)-\frac{1}{2}J_1(t)^2\bigg).$$
Here are $J_1(t)$ and $J_2(t)$ are the $n\!=\!4$ analogues of the functions
in~Subsection~\ref{mirsym_subs}:
$$\sum_{q=0}^{\i} I_q(t)w^q \equiv e^{wt}\sum_{d=0}^{\i}e^{dt}
\frac{\prod_{r=1}^{r=4d}(4w\!+\!r)}{\prod_{r=1}^{r=d}((w\!+\!r)^4\!-\!w^4)},
\qquad J_q(t)=I_q(t)\big/I_0(t).$$
We note that 
\begin{equation}\label{n4concl_e}\begin{split}
\bigg(J_2(t)-\frac{1}{2}J_1(t)^2\bigg)'&
=J_1'(t)\bigg(\frac{J_2'(t)}{J_1'(t)}-J_1(t)\bigg);\\
\bigg(\frac{J_2'(t)}{J_1'(t)}-J_1(t)\bigg)'&=
\bigg(\frac{J_2'(t)}{J_1'(t)}\bigg)'-J_1'(t)\equiv I_{2,2}(t)-I_{1,1}(t)=0.
\end{split}\end{equation}
The last equality above holds by the $(n,p)\!=\!(4,1)$ case of~\e_ref{hypergeom_prp_e3}.
Since $J_2(t)-\frac{1}{2}J_1(t)^2$ is a power series in $e^t$ with no $e^t$-constant
term, \e_ref{n4concl_e} implies that it is zero as expected.\footnote{For 
a surface $X$, the standard and reduced (genus~$1$) GW-invariants
are the same if no descendant classes are involved.}\\

\noindent
The $n\!=\!5$ case of Theorem~\ref{main_thm0} implies Theorem~\ref{bcov_thm}.
In this case, the power series $I_{0,0}(t)$ and $I_{1,1}(t)$ in $e^t$ in the statement of 
Theorem~\ref{main_thm0} are $I_0(t)$ and $J_1'(t)$ in the notation 
of Subsection~\ref{mirsym_subs}.
Thus, the sum of the terms on the first two lines of the right-hand side in the formula
of Theorem~\ref{main_thm0} is precisely the right-hand side of~\e_ref{bcov_thm}
divided by~$2$.
The remaining term in Theorem~\ref{main_thm0} is a sum of two terms, one of which  
(the one corresponding to $p\!=\!2$) is easily seen to be zero.
The other term~is
\begin{equation}\label{n5concl_e}\begin{split}
\frac{5}{24}\cdot\frac{1}{3!}\cdot
\bigg\{\frac{d}{dw}\bigg\}^3
\bigg(\ln\sum_{q=0}^{\i}J_q(t)w^q\bigg)\bigg|_{w=0}
&=\frac{5}{24}\bigg(J_3(t)-J_1(t)J_2(t)+\frac{1}{3}J_1(t)^3\bigg)\\
&=-\frac{1}{12}\sum_{d=1}^{\i}N_{0,d}e^{dT},
\end{split}\end{equation}
with $J_1(t)$, $J_2(t)$, and $J_3(t)$ as in~Subsection~\ref{mirsym_subs}.
The last equality in~\e_ref{n5concl_e} is immediate from~\e_ref{g0ms_e}.
Theorem~\ref{bcov_thm} thus follows from Theorem~\ref{main_thm} and~\e_ref{genus1_e3}.\\

\noindent
If $n\!=\!6$, $X$ is a sextic fourfold in $\P^5$.
Theorem~\ref{main_thm0} in this case gives
\begin{equation*}\label{n6concl_e}\begin{split}
&\sum_{d=1}^{\i}e^{dT}\GW_{1,0}^{0;X}(d;1)
=-\frac{35}{2}\big(T\!-\!t\big)+\frac{423}{4}\ln I_0(t)
-\ln J_1'(t)-\frac{1}{24}\ln(1\!-\!6^6e^t)\\
&\qquad
+\frac{6}{24}\cdot 15\cdot\frac{1}{2!}\bigg\{\frac{d}{dw}\bigg\}^2
\bigg(\ln\sum_{q=0}^{\i}J_q(t)w^q\bigg)\bigg|_{w=0}
+\frac{6}{24}\cdot \frac{1}{4!}\bigg\{\frac{d}{dw}\bigg\}^4
\bigg(\ln\sum_{q=0}^{\i}J_q(t)w^q\bigg)\bigg|_{w=0}.
\end{split}\end{equation*}
The terms on the second line above arise from 
the last term in the formula of Theorem~\ref{main_thm0}.
In the $n\!=\!3,4,5$ cases, the latter is 
$$\GW_{1,k}^{0;X}(d;1)-\GW_{1,k}^X(d;1).$$
We show in~\cite{g1diff} that this is the case for all $n$.
For $n\!=\!6$, Theorem~\ref{main_thm0} would then give
$$\sum_{d=1}^{\i}e^{dT}\GW_{1,0}^X(d;1)
=-\frac{35}{2}\big(T\!-\!t\big)+\frac{423}{4}\ln I_0(t)
-\ln J_1'(t)-\frac{1}{24}\ln(1\!-\!6^6e^t),$$
confirming the mirror symmetry prediction of~\cite[Section~6.1]{KlPa}.\footnote{The variable 
$t$ in~\cite[(46)]{KlPa} is not the same as the variable $t$ in this paper; 
see Appendix~\ref{comparison_app}.}

\section{Equivariant Cohomology and Stable Maps} 
\label{introlocal_sec}

\subsection{Definitions and Notation}
\label{equivcoh_subs}

\noindent
This subsection reviews the notion of equivariant cohomology and sets up related notation
that will be used throughout the rest of the paper.
For the most part, our notation agrees with \cite[Chapters~29,30]{MirSym};
the main difference is that we work with $\P^{n-1}$ instead of~$\P^n$.\\

\noindent
We denote by $\T$ the $n$-torus $(\C^*)^n$ (or $(S^1)^n$).
It acts freely on $E\T\!=\!(\C^{\i})^n\!-\!0$ (or $(S^{\i})^n$):
$$\big(e^{\I\th_1},\ldots,e^{\I\th_n}\big)\cdot (z_1,\ldots,z_n)
=\big(e^{\I\th_1}z_1,\ldots,e^{\I\th_n}z_n\big).$$
Thus, the classifying space for $\T$ and its group cohomology are given by
$$\qquad B\T\equiv E\T/\T=(\P^{\i})^n  \qquad\hbox{and}\qquad
H_{\T}^*\equiv H^*(B\T;\Q)=\Q[\al_1,\ldots,\al_n],$$
where $\al_i\!=\!\pi_i^*c_1(\ga^*)$ if
$$\pi_i\!: (\P^{\i})^n\lra\P^{\i} \qquad\hbox{and}\qquad
\ga\lra\P^{\i}$$
are the projection onto the $i$-th component and the tautological line bundle, respectively.
Denote by $\H_{\T}^*$ the field of fractions of $H_{\T}^*$:
$$\H_{\T}^*=\Q_{\al}\equiv \Q(\al_1,\ldots,\al_n).$$\\

\noindent
A representation $\rho$ of $\T$, i.e.~a linear action of $\T$ on $\C^k$,
induces a vector bundle over~$B\T$:
$$V_{\rho}\equiv E{\T}\times_{\T}\C^k.$$
If $\rho$ is one-dimensional, we will call
$$c_1(V_{\rho}^*)=-c_1(V_{\rho})\in H_{\T}^*\subset \H_{\T}^*$$ 
the \sf{weight} of $\rho$.
For example, $\al_i$ is the weight of representation
\begin{equation}\label{ifactorrep_e}
\pi_i\!: \T\lra\C^*, \qquad 
\big(e^{\I\th_1},\ldots,e^{\I\th_n}\big)\cdot z =  e^{\I\th_i}z.
\end{equation}
More generally, if a representation $\rho$ of $\T$ on $\C^k$ splits 
into one-dimensional representations with weights $\be_1,\ldots,\be_k$,
we will call $\be_1,\ldots,\be_k$ the \sf{weights} of~$\rho$.
In such a case,
\begin{equation}\label{weightspord_e}
e(V_{\rho}^*)=\be_1\cdot\ldots\cdot\be_k.
\end{equation}
We will call the representation $\rho$ of $\T$ on $\C^n$ with weights $\al_1,\ldots,\al_n$
the \sf{standard representation} of~$\T$.\\

\noindent
If $\T$ acts on a topological space $M$, let
$$H_{\T}^*(M)\equiv H^*(BM;\Q), \qquad\hbox{where}\qquad BM=E\T\!\times_{\T}\!M,$$
denote the corresponding \sf{equivariant cohomology} of $M$.
The projection map $BM\!\lra\!B\T$ induces an action of $H_{\T}^*$ on $H_{\T}^*(M)$.
Let
$$\H_{\T}^*(M)=H_{\T}^*(M)\otimes_{H_{\T}^*}\H_{\T}^*.$$
If the $\T$-action on $M$ lifts to an action on a (complex) vector bundle $V\!\lra\!M$,
then 
$$BV\equiv E{\T}\!\times_{\T}\!V$$
is a vector bundle over $BM$.
Let
$$\E(V)\equiv e(BV)\in H_{\T}^*(M)\subset \H_{\T}^*(M)$$
denote the \sf{equivariant euler class of} $V$.\\

\noindent
Throughout the paper we work with the standard action of $\T$ on $\P^{n-1}$,
i.e.~the action induced by the standard action $\rho$ of $\T$ on~$\C^n$:
$$\big(e^{\I\th_1},\ldots,e^{\I\th_n}\big)\cdot [z_1,\ldots,z_n] 
=\big[e^{\I\th_1}z_1,\ldots,e^{\I\th_n}z_n\big].$$
Since $B\P^{n-1}=\P V_{\rho}$, 
$$H_{\T}^*(\P^{n-1})\equiv H^*\big(\P V_{\rho};\Q\big)
= \Q[x,\al_1,\ldots,\al_n]
\big/\big(x^n\!+\!c_1(V_{\rho})x^{n-1}\!+\!\ldots\!+\!c_n(V_{\rho})\big),$$
where $x\!=\!c_1(\ti\ga^*)$ and $\ti\ga\!\lra\!\P V_{\rho}$ 
is the tautological line bundle.
Since 
$$c(V_{\rho})=(1-\al_1)\ldots(1-\al_n),$$
it follows that 
\begin{equation}\label{pncoh_e}\begin{split}
H_{\T}^*(\P^{n-1}) &= \Q[x,\al_1,\ldots,\al_n]\big/(x\!-\!\al_1)\ldots(x\!-\!\al_n),\\
\H_{\T}^*(\P^{n-1}) &= \Q_{\al}[x]\big/(x\!-\!\al_1)\ldots(x\!-\!\al_n).
\end{split}\end{equation}\\

\noindent
The standard action of $\T$ on $\P^{n-1}$ has $n$~fixed points:
$$P_1=[1,0,\ldots,0], \qquad P_2=[0,1,0,\ldots,0], 
\quad\ldots\quad P_n=[0,\ldots,0,1].$$
For each $i\!=\!1,2,\ldots,n$, let
\begin{equation}\label{phidfn_e}
\phi_i= \prod_{k\neq i}(x\!-\!\al_k) \in H_{\T}^*(\P^{n-1}).
\end{equation}
By equation~\e_ref{phiprop_e} below, $\phi_i$ is the equivariant Poincare dual of $P_i$.
We also note that $\ti\ga|_{BP_i}\!=\!V_{\pi_i}$, where $\pi_i$ is as in~\e_ref{ifactorrep_e}.
Thus, the restriction map on the equivariant cohomology induced by the inclusion 
$P_i\!\lra\!\P^{n-1}$ is given~by
\begin{equation}\label{restrmap_e}
H_{\T}^*(\P^{n-1})=\Q[x,\al_1,\ldots,\al_n]\big/\prod_{k=1}^{k=n}(x\!-\!\al_k)
\lra H_{\T}^*(P_i)=\Q[\al_1,\ldots,\al_n], \qquad x\lra\al_i.
\end{equation}
By~\e_ref{restrmap_e},
\begin{equation}\label{uniquecond_e}
\eta=0 \in H_{\T}^*(\P^{n-1}) \qquad\Llra\qquad
\eta|_{P_i}=0\in H_{\T}^* ~~\forall~i=1,2,\ldots,n.
\end{equation}\\

\noindent
The tautological line bundle $\ga_{n-1}\!\lra\!\P^{n-1}$ is a subbundle of 
$\P^{n-1}\!\times\!\C^n$ preserved by the diagonal action of~$\T$.
Thus, the action of $\T$ on $\P^{n-1}$ naturally lifts to an action 
on~$\ga_{n-1}$ and
\begin{equation}\label{garestr_e}
\E\big(\ga_{n-1}^*\big)\big|_{P_i}=\al_i \qquad\forall~i=1,2,\ldots,n.
\end{equation}
The $\T$-action on $\P^{n-1}$ also has a natural lift to the vector bundle
$T\P^{n-1}\!\lra\!\P^{n-1}$ so that  there is a short exact sequence
$$0\lra \ga_{n-1}^*\otimes\ga_{n-1}
\lra \ga_{n-1}^*\otimes\big(\P^{n-1}\!\times\!\C^n\big) 
\lra T\P^{n-1} \lra0$$
of $\T$-equivariant vector bundles on $\P^{n-1}$.
By~\e_ref{weightspord_e}, \e_ref{garestr_e}, and~\e_ref{phidfn_e}, 
\begin{equation}\label{tangrestr_e}
\E\big(T\P^{n-1}\big)\big|_{P_i}=\prod_{k\neq i}(\al_i\!-\!\al_k)
=\phi_i|_{P_i}
\qquad\forall~i=1,2,\ldots,n.
\end{equation}\\

\noindent
If $\T$ acts smoothly on a smooth compact oriented manifold $M$, there is a well-defined
integration-along-the-fiber homomorphism
$$\int_M\!: H_{\T}^*(M)\lra H_{\T}^*$$
for the fiber bundle $BM\!\lra\!B\T$.
The classical localization theorem of~\cite{ABo} relates it to 
integration along the fixed locus of the $\T$-action.
The latter is a union of smooth compact orientable manifolds~$F$;
$\T$ acts on the normal bundle $\N F$ of each~$F$.
Once an orientation of $F$ is chosen, there is a well-defined 
integration-along-the-fiber homomorphism
$$\int_F\!: H_{\T}^*(F)\lra H_{\T}^*.$$
The localization theorem states that 
\begin{equation}\label{ABothm_e}
\int_M\eta = \sum_F\int_F\frac{\eta|_F}{\E(\N F)} \in \H^*_{\T}
\qquad\forall~\eta\in H_{\T}^*(M),
\end{equation}
where the sum is taken over all components $F$ of the fixed locus of $\T$. 
Part of the statement of~\e_ref{ABothm_e} is that $\E(\N F)$
is invertible in~$\H_{\T}^*(F)$.
In the case of the standard action of $\T$ on~$\P^{n-1}$, \e_ref{ABothm_e} implies that 
\begin{equation}\label{phiprop_e}
\eta|_{P_i}=\int_{\P^{n-1}}\eta\phi_i \in\H_{\T}^*
\qquad\qquad\forall~\eta\!\in\!\H_{\T}^*(\P^{n-1}),~i=1,2,\ldots,n;
\end{equation}
see also \e_ref{tangrestr_e}.\\

\noindent
Finally, if $f\!:M\!\lra\!M'$ is a $\T$-equivariant map between two compact oriented 
manifolds, there is a well-defined pushforward homomorphism
$$f_*\!: H_{\T}^*(M) \lra H_{\T}^*(M').$$
It is characterized by the property that 
\begin{equation}\label{pushdfn_e}
\int_{M'}(f_*\eta)\,\eta'=\int_M\eta\,(f^*\eta') 
\qquad~\forall~\eta\!\in\! H_{\T}^*(M),\, \eta'\!\in\! H_{\T}^*(M').
\end{equation}
The homomorphism $\int_M$ of the previous paragraph corresponds to $M'$ being a point.
It is immediate from~\e_ref{pushdfn_e} that 
\begin{equation}\label{pushprop_e}
f_*\big(\eta\,(f^*\eta')\big)=(f_*\eta)\,\eta' 
\qquad~\forall~\eta\!\in\! H_{\T}^*(M),\, \eta'\!\in\! H_{\T}^*(M').
\end{equation}

\subsection{Setup for Localization Computation on $\wt\M_{1,1}^0(\P^{n-1},d)$}
\label{local_subs}

\noindent
The standard $\T$-action on $\P^{n-1}$ (as well as any other action) induces 
$\T$-actions on moduli spaces of stable maps $\ov\M_{g,k}(\P^{n-1},d)$ 
by composition on the right,
$$h\cdot[\cC,f]=[\cC,h\circ f]  \qquad\forall~h\in\T, 
~[\cC,f]\in\ov\M_{g,k}(\P^{n-1},d),$$
and lifts to an action on $\wt\M_{1,k}^0(\P^{n-1},d)$.
All the evaluation maps,
$$\ev_i\!:\ov\M_{g,k}(\P^{n-1},d),\wt\M_{1,k}^0(\P^{n-1},d)\lra\P^{n-1}, 
\quad [\cC,y_1,\ldots,y_k,f]\lra f(y_i), \qquad i=1,2,\ldots,k,$$
are $\T$-equivariant.
These actions lift naturally to the tautological tangent line bundles
$$L_1,\ldots,L_k\lra \ov\M_{g,k}(\P^{n-1},d);$$
see \cite[Section~25.2]{MirSym}.
Let
$$\psi_i\equiv c_1(L_i^*)\in\H_{\T}^*\big(\ov\M_{g,k}(\P^{n-1},d)\big)$$
denote the equivariant $\psi$-class.\\

\noindent
Via the natural lift of the $\T$-action to $\ga_{n-1}\!\lra\!\P^{n-1}$
described in Subsection~\ref{equivcoh_subs}, 
the $\T$-actions on $\ov\M_{g,k}(\P^{n-1},d)$ and $\wt\M_{1,k}^0(\P^{n-1},d)$ 
lift to $\T$-actions on the sheafs $\pi_*\ev^*\O_{\P^{n-1}}(a)$ and the vector bundle
$$\V_1\lra \wt\M_{1,1}^0(\P^{n-1},d)$$
introduced in Subsection~\ref{appr_subs}.
We denote by 
$$\V_0\lra\ov\M_{0,k}(\P^{n-1},d)$$ 
the vector bundle corresponding to the locally free sheaf 
$\pi_*\ev^*\O_{\P^{n-1}}(n)$ on $\ov\M_{0,k}(\P^{n-1},d)$.
Let 
$$\cL=\ga_{n-1}^{*\otimes n}\lra \P^{n-1}$$
be the vector bundle corresponding to the locally free sheaf 
$\O_{\P^{n-1}}(n)\!\lra\!\P^{n-1}$.
For $g\!=\!0,1$, the equivariant bundle map
$$\wt\ev_1\!: \V_g\lra\ev_1^*\cL, \qquad
[\cC,y_1,\ldots,y_k,f,\xi]\lra\big[\xi(y_1)\big],$$
is surjective.
Thus,
$$\V_0'\equiv\ker\wt\ev_1\lra\ov\M_{0,k}(\P^{n-1},d)
\qquad\hbox{and}\qquad
\V_1'\equiv\ker\wt\ev_1\lra\wt\M_{1,1}^0(\P^{n-1},d)$$
are equivariant vector bundles.\footnote{In \cite[Chapters 29,30]{MirSym}, 
the analogues of $\V_0$ and $\V_0'$ over $\ov\M_{0,2}(\P^n,d)$ are denoted by
$E_{0,d}$ and $E_{0,d}'$, respectively.
However, $E_{0,d}'$ is the kernel of the evaluation map at the second marked point.}
Furthermore,
\begin{equation}\label{eulerclass_e}
\E(\V_g)=\E\big(\ev_1^*\cL\big)\,\E(\V_g')=n\,\ev_1^*(x)\,\E(\V_g'), 
\end{equation}
where $x\!\in\!H_{\T}^2(\P^{n-1})$ is the equivariant hyperplane class
as in Subsection~\ref{equivcoh_subs}.\\

\noindent
We denote by $\al$ the tuple $(\al_1,\ldots,\al_n)$. 
With $\ev_{1,d}$ denoting the evaluation map on $\wt\M_{1,1}^0(\P^{n-1},d)$,
let
$$\F(\al,x,u)=\sum_{d=1}^{\i}u^d(\ev_{1,d*}\E(\V_1)\big)\in 
\big(H_{\T}^{n-2}(\P^{n-1})\big)[[u]].$$
By~\e_ref{pncoh_e}, 
\begin{equation}\label{pushclass_e}
\F(\al,x,u)=\F_0(u)x^{n-2}+\F_1(\al,u)x^{n-3}+\ldots+\F_{n-2}(\al,u)x^0,
\end{equation}
for some power series $\F_0(u)$ in $u$ and degree~$p$ homogeneous $\al$-polynomials 
$$\F_p(\al,u)\in\Q[[u]]\big[\al_1,\ldots,\al_n\big].$$
These polynomials must be symmetric in $\al_1,\ldots,\al_n$.
Note that by~\e_ref{derivgen_e1} and \e_ref{pushprop_e}, 
\begin{equation}\label{derivgen_e2}
\frac{d}{dT}\sum_{d=1}^{\i}e^{dT}\GW_{1,0}^{0;X}(d;1)=\F_0(e^T).
\end{equation}
Thus, our aim is to determine the power series $\F_0(u)$ in $u$ 
defined by~\e_ref{pushclass_e}.\\

\noindent
By~\e_ref{restrmap_e}, \e_ref{phiprop_e}, and~\e_ref{pushdfn_e}, 
\begin{equation}\label{Frestr_e}\begin{split}
\F(\al,\al_i,u) &=\F(\al,x,u)\big|_{P_i}
=\sum_{d=1}^{\i}u^d\int_{\P^{n-1}}\big(\ev_{1,d*}\E(\V_1)\big)\phi_i\\
&=\sum_{d=1}^{\i}u^d\int_{\wt\M_{1,1}^0(\P^{n-1},d)}\E(\V_1)\ev_1^*\phi_i
\end{split}\end{equation}
for each $i\!=\!1,2,\ldots,n$.
By~\e_ref{uniquecond_e}, the power series $\F_0(u)$ is completely 
determined by 
$$\F(\al,\al_1,u),\ldots,\F(\al,\al_n,u)\in\Q_{\al}[[u]],$$
where $\Q_{\al}$ is the ring $\Q(\al_1,\ldots,\al_n)$
of rational fractions in $\al_1,\ldots,\al_n$.
We will apply the localization formula~\e_ref{ABothm_e} to the last 
expression in~\e_ref{Frestr_e}.
In order to do so, we need to describe the fixed loci of the $\T$-action 
on $\wt\M_{1,1}^0(\P^{n-1},d)$ and for each fixed locus $F$ the corresponding triple $(F,\E(\V_1)\ev_1^*\phi_i|_F,\N F)$ or another triple $(F',\eta,\N')$ such that 
$$\int_{F'}\frac{\eta}{\E(\N')}=\int_F\frac{\E(\V_1)\ev_1^*\phi_i|_F}{\E(\N F)}.$$
In one case in Subsection~\ref{localdata_subs2}, choosing such a replacement 
turns out to be advantageous.\\

\noindent
An element $[\cC,f]$ of $\ov\M_{1,1}(\P^{n-1},d)$ is the equivalence class of
a nodal genus-one curve $\cC$ with one marked point and 
a stable degree-$d$ map $f\!:\cC\!\lra\!\P^{n-1}$.
We denote by
$$\M_{1,1}^{\eff}(\P^{n-1},d)\subset \ov\M_{1,1}^0(\P^{n-1},d)$$
the open subset of $\ov\M_{1,1}(\P^{n-1},d)$ consisting of the stable maps
$[\cC,f]$ that are not constant on the \sf{principal}, 
genus-one, component~$\cC_P$ of~$\cC$.\footnote{The connected curve~$\cC$ is nodal and 
has arithmetic genus one. 
Thus, either one of the component of~$\cC$ is a smooth torus or 
$\cC$ contains a circle of one or more spheres (each irreducible component
is a $\P^1$ with exactly two nodes). 
In the first case, $\cC_P$ is the smooth torus; in the second, $\cC_P$ is 
the circle of spheres.}
By definition, $\M_{1,1}^{\eff}(\P^{n-1},d)$ is dense in $\ov\M_{1,1}^0(\P^{n-1},d)$.
Let 
$$\partial\ov\M_{1,1}^0(\P^{n-1},d)
=\ov\M_{1,1}^0(\P^{n-1},d)-\M_{1,1}^{\eff}(\P^{n-1},d).$$
The desingularization $\wt\M_{1,1}^0(\P^{n-1},d)$ of 
$\ov\M_{1,1}^0(\P^{n-1},d)$ is obtained by blowing up along subvarieties
contained in $\partial\ov\M_{1,1}^0(\P^{n-1},d)$;
see \cite[Subsection~1.2]{VaZ}.
Thus,
$$\M_{1,1}^{\eff}(\P^{n-1},d)\subset \wt\M_{1,1}^0(\P^{n-1},d)$$
is a dense open subset. Let
$$\partial\wt\M_{1,1}^0(\P^{n-1},d)
=\wt\M_{1,1}^0(\P^{n-1},d)-\M_{1,1}^{\eff}(\P^{n-1},d).$$\\

\noindent
Since each of the fixed loci of the $\T$-action on $\ov\M_{1,1}^0(\P^{n-1},d)$
is contained either in $\M_{1,1}^{\eff}(\P^{n-1},d)$ or in
$\partial\ov\M_{1,1}^0(\P^{n-1},d)$,
each of the fixed loci of the $\T$-action on $\wt\M_{1,1}^0(\P^{n-1},d)$
is contained either in $\M_{1,1}^{\eff}(\P^{n-1},d)$ or in
$\partial\wt\M_{1,1}^0(\P^{n-1},d)$.
Furthermore, the fixed loci contained in $\M_{1,1}^{\eff}(\P^{n-1},d)$
and the corresponding triples $(F,\eta,\N F)$ as in~\e_ref{ABothm_e}
are the same for the $\T$-actions on $\ov\M_{1,1}^0(\P^{n-1},d)$
and $\wt\M_{1,1}^0(\P^{n-1},d)$.
These loci and their total contribution to~\e_ref{Frestr_e} are described 
in Subsection~\ref{localdata_subs1}.
The fixed loci contained in $\partial\wt\M_{1,1}^0(\P^{n-1},d)$ and their total
contribution to~\e_ref{Frestr_e} are described in Subsection~\ref{localdata_subs2} 
based on \cite[Subsection~1.4]{VaZ}.\footnote{We will not describe 
$\wt\M_{1,1}^0(\P^{n-1},d)$ in this paper as this is not necessary.}\\

\noindent
Many expressions throughout the paper involve residues of rational functions
in a complex variable~$\hb$.
If $f\!=\!f(\hb)$ is a rational function in $\hb$ and $\hb_0\!\in\!S^2$, 
we denote by $\fR_{\hb=\hb_0}f(\hb)$ the residue of $f(\hb)d\hb$ at $\hb\!=\!\hb_0$:
$$\fR_{\hb=\hb_0}f(\hb) = \frac{1}{2\pi\I}\oint f(\hb)d\hb,$$
where the integral is taken over a positively oriented loop around $\hb\!=\!\hb_0$
containing no other singular points of $f$.
With this definition,
$$\fR_{\hb=\i}f(\hb)=-\fR_{w=0}\big\{w^{-2}f(w^{-1})\big\}.$$
If $f$ involves variables other than $\hb$, 
$\fR_{\hb=\hb_0}f(\hb)$ will be a function  of such variables.
If $f$ is a power series in $u$ with coefficients that are rational
functions in~$\hb$ and possibly other variables, denote by 
$\fR_{\hb=\hb_0}f(\hb)$ the  power series in~$u$ obtained by 
replacing each of the coefficients by its residue at $\hb\!=\!\hb_0$.
If $\hb_1,\ldots,\hb_k$ is a collection of points in~$S^2$, let
$$\fR_{\hb=\hb_1,\ldots,\hb_k}f(\hb) =\sum_{i=1}^{i=k} \fR_{\hb=\hb_i}f.$$
Finally, we will denote by $\bar\Z^+$ the set of non-negative integers and by $[n]$,
whenever $n\!\in\!\bar\Z^+$, the set of positive integers not exceeding~$n$:
$$\bar\Z^+\equiv\big\{0,1,2,\ldots\big\}, \qquad
[n]=\big\{1,2,\ldots,n\big\}.$$

\subsection{Contributions from Fixed Loci, I}
\label{localdata_subs1}

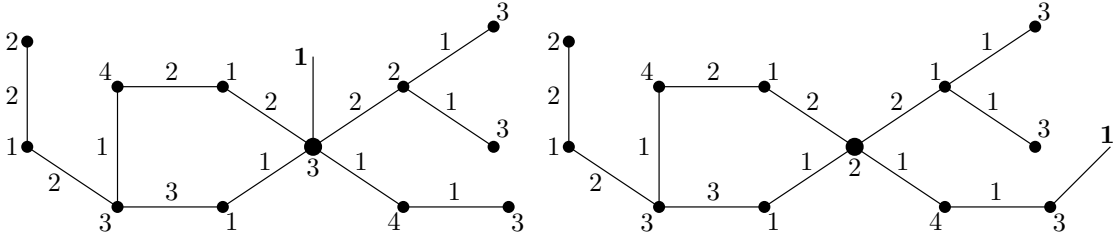
\begin{figure}
\begin{pspicture}(-.3,-1.1)(10,2)
\psset{unit=.4cm}
\pscircle*(12,0){.3}\rput(12,-.7){\smsize{$3$}}
\psline[linewidth=.04](12,0)(12,3)\rput(11.6,3){\smsize{$\bf 1$}}
\psline[linewidth=.04](12,0)(9,2)\pscircle*(9,2){.2}
\rput(10.6,1.5){\smsize{$2$}}\rput(9.3,2.5){\smsize{$1$}}
\psline[linewidth=.04](12,0)(9,-2)\pscircle*(9,-2){.2}
\rput(10.4,-.5){\smsize{$1$}}\rput(9.3,-2.5){\smsize{$1$}}
\psline[linewidth=.04](5.5,2)(9,2)\pscircle*(5.5,2){.2}
\rput(7.3,2.5){\smsize{$2$}}\rput(5.1,2.5){\smsize{$4$}}
\psline[linewidth=.04](5.5,-2)(9,-2)\pscircle*(5.5,-2){.2}
\rput(7.3,-1.5){\smsize{$3$}}\rput(5.1,-2.5){\smsize{$3$}}
\psline[linewidth=.04](5.5,-2)(5.5,2)\rput(5,0){\smsize{$1$}}
\psline[linewidth=.04](2.5,0)(5.5,-2)\pscircle*(2.5,0){.2}
\rput(3.4,-1.2){\smsize{$2$}}\rput(2,0){\smsize{$1$}}
\psline[linewidth=.04](2.5,0)(2.5,3.5)\pscircle*(2.5,3.5){.2}
\rput(2,1.7){\smsize{$2$}}\rput(2,3.5){\smsize{$2$}}
\psline[linewidth=.04](12,0)(15,2)\pscircle*(15,2){.2}
\rput(13.4,1.5){\smsize{$2$}}\rput(14.7,2.5){\smsize{$2$}}
\psline[linewidth=.04](12,0)(15,-2)\pscircle*(15,-2){.2}
\rput(13.6,-.5){\smsize{$1$}}\rput(14.7,-2.5){\smsize{$4$}}
\psline[linewidth=.04](18,4)(15,2)\pscircle*(18,4){.2}
\rput(16.4,3.5){\smsize{$1$}}\rput(18.3,4.5){\smsize{$3$}}
\psline[linewidth=.04](18,0)(15,2)\pscircle*(18,0){.2}
\rput(16.6,1.5){\smsize{$1$}}\rput(18.3,.6){\smsize{$3$}}
\psline[linewidth=.04](18.5,-2)(15,-2)\pscircle*(18.5,-2){.2}
\rput(16.7,-1.5){\smsize{$1$}}\rput(18.8,-2.5){\smsize{$3$}}
\pscircle*(30,0){.3}\rput(30,-.7){\smsize{$2$}}
\psline[linewidth=.04](30,0)(27,2)\pscircle*(27,2){.2}
\rput(28.6,1.5){\smsize{$2$}}\rput(27.3,2.5){\smsize{$1$}}
\psline[linewidth=.04](30,0)(27,-2)\pscircle*(27,-2){.2}
\rput(28.4,-.5){\smsize{$1$}}\rput(27.3,-2.5){\smsize{$1$}}
\psline[linewidth=.04](23.5,2)(27,2)\pscircle*(23.5,2){.2}
\rput(25.3,2.5){\smsize{$2$}}\rput(23.1,2.5){\smsize{$4$}}
\psline[linewidth=.04](23.5,-2)(27,-2)\pscircle*(23.5,-2){.2}
\rput(25.3,-1.5){\smsize{$3$}}\rput(23.1,-2.5){\smsize{$3$}}
\psline[linewidth=.04](23.5,-2)(23.5,2)\rput(23,0){\smsize{$1$}}
\psline[linewidth=.04](20.5,0)(23.5,-2)\pscircle*(20.5,0){.2}
\rput(21.4,-1.2){\smsize{$2$}}\rput(20,0){\smsize{$1$}}
\psline[linewidth=.04](20.5,0)(20.5,3.5)\pscircle*(20.5,3.5){.2}
\rput(20,1.7){\smsize{$2$}}\rput(20,3.5){\smsize{$2$}}
\psline[linewidth=.04](30,0)(33,2)\pscircle*(33,2){.2}
\rput(31.4,1.5){\smsize{$2$}}\rput(32.7,2.5){\smsize{$1$}}
\psline[linewidth=.04](30,0)(33,-2)\pscircle*(33,-2){.2}
\rput(31.6,-.5){\smsize{$1$}}\rput(32.7,-2.5){\smsize{$4$}}
\psline[linewidth=.04](36,4)(33,2)\pscircle*(36,4){.2}
\rput(34.4,3.5){\smsize{$1$}}\rput(36.3,4.5){\smsize{$3$}}
\psline[linewidth=.04](36,0)(33,2)\pscircle*(36,0){.2}
\rput(34.6,1.5){\smsize{$1$}}\rput(36.3,.6){\smsize{$3$}}
\psline[linewidth=.04](36.5,-2)(33,-2)\pscircle*(36.5,-2){.2}
\rput(34.7,-1.5){\smsize{$1$}}\rput(36.8,-2.5){\smsize{$3$}}
\psline[linewidth=.04](38.5,0)(36.5,-2)\rput(38.4,.4){\smsize{$\bf 1$}}
\end{pspicture}
\caption{A decorated graph of type $A_3$ and a decorated graph of type $\ti{A}_{32}$}
\label{loopgraph_fig}
\end{figure}

\noindent
As described in detail in \cite[Section~27.3]{MirSym},
the fixed loci of the $\T$-action on $\ov\M_{g,k}(\P^{n-1},d)$ 
are indexed by \sf{decorated graphs}.
A \sf{graph} consists of a set~$\Ver$ of \sf{vertices} and 
a collection $\Edg$ of \sf{edges}, i.e.~of two-element subsets 
of~$\Ver$.\footnote{If $g\!=\!0$, $\Edg$ can be taken to be a subset of
the set $\Sym^2(\Ver)$ of two-element subsets of~$\Ver$.
For $g\!\ge\!1$, $\Edg$ should be viewed as a map from a finite set $\Dom(\Edg)$
to $\Sym^2(\Ver)$; this map may not be injective (i.e.~there can be multiple
edges connecting a pair of vertices).
In the latter case, $e\!\in\!\Edg$ will mean that $e$ is an element of $\Dom(\Edg)$;
if $v\!\in\!\Ver$ and $e\!\in\!\Edg$, $v\!\in\!e$ will mean that $v$ 
is an element of the image of~$e$ in $\Sym^2(\Ver)$;
a map from $\Edg$ will mean a map from $\Dom(\Edg)$.}
In Figure~\ref{loopgraph_fig}, the vertices are represented by dots,
while each edge $\{v_1,v_2\}$ is shown as the line segment between $v_1$ and~$v_2$.
For the purposes of describing the fixed loci of 
$\ov\M_{0,k}(\P^{n-1},d)$ and $\M_{1,k}^{\eff}(\P^{n-1},d)$, 
it is sufficient to define a \sf{decorated graph} as a tuple
\begin{equation}\label{decortgraphdfn_e}
\Ga = \big(\Ver,\Edg;\mu,\d,\eta\big),
\end{equation}
where $(\Ver,\Edg)$ is a graph and
$$\mu\!:\Ver\lra [n], \qquad
\d\!: \Edg\lra\Z^+, \quad\hbox{and}\quad \eta\!:[k]\lra\Ver$$
are maps such that
\begin{equation}\label{decorgraphcond_e}
\mu(v_1)\neq\mu(v_2)  \qquad\hbox{if}\quad \{v_1,v_2\}\in\Edg.
\end{equation}
In Figure~\ref{loopgraph_fig}, the value of the map $\mu$ on each vertex
is indicated by the number next to the vertex.
Similarly, the value of the map $\d$ on each edge is indicated by the number next to the edge.
The only element of the~set $[k]\!=\![1]$ is shown in bold face.
It is linked by a line segment to its image under~$\eta$.
By~\e_ref{decorgraphcond_e}, no two consecutive vertex labels are the same.\\

\begin{figure}
\begin{pspicture}(-.3,-1.2)(10,2)
\psset{unit=.4cm}
\psline[linewidth=.04](15,1)(18,1)\rput(18.4,1){\smsize{$\bf 2$}}
\pscircle*(15,1){.2}\rput(15,1.6){\smsize{$2$}}
\psline[linewidth=.04](15,-1)(18,-1)\rput(18.4,-1){\smsize{$\bf 1$}}
\pscircle*(15,-1){.2}\rput(15,-1.6){\smsize{$2$}}
\psline[linewidth=.04](15,1)(12,2)\pscircle*(12,2){.2}
\rput(13.6,2){\smsize{$2$}}\rput(12.3,2.5){\smsize{$1$}}
\psline[linewidth=.04](15,-1)(12,-2)\pscircle*(12,-2){.2}
\rput(13.4,-1){\smsize{$1$}}\rput(12.3,-2.5){\smsize{$1$}}
\psline[linewidth=.04](8.5,2)(12,2)\pscircle*(8.5,2){.2}
\rput(10.3,2.5){\smsize{$2$}}\rput(8.1,2.5){\smsize{$4$}}
\psline[linewidth=.04](8.5,-2)(12,-2)\pscircle*(8.5,-2){.2}
\rput(10.3,-1.5){\smsize{$3$}}\rput(8.1,-2.5){\smsize{$3$}}
\psline[linewidth=.04](8.5,-2)(8.5,2)\rput(8,0){\smsize{$1$}}
\psline[linewidth=.04](5.5,0)(8.5,-2)\pscircle*(5.5,0){.2}
\rput(6.4,-1.2){\smsize{$2$}}\rput(5,0){\smsize{$1$}}
\psline[linewidth=.04](5.5,0)(5.5,3.5)\pscircle*(5.5,3.5){.2}
\rput(5,1.7){\smsize{$2$}}\rput(5,3.5){\smsize{$2$}}
\rput(11.8,0){\lgsize{$\bf{\Ga_{\pm}}$}}
\psline[linewidth=.04](28,1)(25,1)\rput(24.6,1){\smsize{$\bf 1$}}
\pscircle*(28,1){.2}\rput(28,1.6){\smsize{$2$}}
\psline[linewidth=.04](28,-1)(25,-1)\rput(24.6,-1){\smsize{$\bf 2$}}
\pscircle*(28,-1){.2}\rput(28,-1.6){\smsize{$2$}}
\psline[linewidth=.04](28,1)(31,2)\pscircle*(31,2){.2}
\rput(29.4,1.9){\smsize{$2$}}\rput(30.7,2.5){\smsize{$2$}}
\psline[linewidth=.04](28,-1)(31,-2)\pscircle*(31,-2){.2}
\rput(29.6,-1){\smsize{$1$}}\rput(30.7,-2.5){\smsize{$4$}}
\psline[linewidth=.04](34,4)(31,2)\pscircle*(34,4){.2}
\rput(32.4,3.5){\smsize{$1$}}\rput(34.3,4.5){\smsize{$3$}}
\psline[linewidth=.04](34,0)(31,2)\pscircle*(34,0){.2}
\rput(32.6,1.5){\smsize{$1$}}\rput(34.3,.6){\smsize{$3$}}
\psline[linewidth=.04](34.5,-2)(31,-2)\pscircle*(34.5,-2){.2}
\rput(32.7,-1.5){\smsize{$1$}}\rput(34.8,-2.5){\smsize{$3$}}
\psline[linewidth=.04](36.5,0)(34.5,-2)\rput(36.4,.4){\smsize{$\bf 1$}}
\rput(28,3){\lgsize{$\bf{\Ga_1}$}}\rput(28,-3){\lgsize{$\bf{\Ga_0}$}}
\end{pspicture}
\caption{The strands of the second graph in Figure~\ref{loopgraph_fig}}
\label{loopgraph_fig2}
\end{figure}
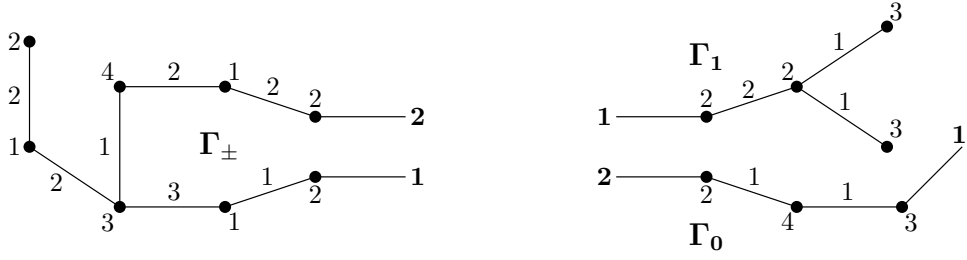

\noindent
A graph $(\Ver,\Edg)$ is a \sf{tree} if it contains no~\sf{loops}, 
i.e.~the set $\Edg$ contains no subset of the~form
$$\big\{\{v_1,v_2\},\{v_2,v_3\},\ldots,\{v_N,v_1\}\big\},
\qquad v_1,\ldots,v_N\!\in\!\Ver,~ N\!\ge\!1.$$
For example, the graphs in Figure~\ref{loopgraph_fig2} are trees,
while those in Figure~\ref{loopgraph_fig} contain one loop each.
Via the construction of the next paragraph, decorated trees describe 
the fixed loci of $\ov\M_{0,k}(\P^{n-1},d)$, while decorated graphs 
with exactly one loop describe the fixed loci of $\M_{1,k}^{\eff}(\P^{n-1},d)$.\\

\noindent
The fixed locus $\cZ_{\Ga}$ of $\ov\M_{g,k}(\P^{n-1},d)$ corresponding to a 
decorated graph $\Ga$ consists of the stable maps~$f$ from a 
genus~$g$ nodal curve $\cC_f$ with $k$ marked points into~$\P^{n-1}$ 
that satisfy the following conditions.
The components of $\cC_f$ on which the map~$f$ is not constant
are rational and correspond to the edges of~$\Ga$.
Furthermore, if $e\!=\!\{v_1,v_2\}$ is an edge,
the restriction of~$f$ to the component $\cC_{f,e}$ corresponding to~$e$
is a degree~$\d(e)$ cover of the line 
$$\P^1_{\mu(v_1),\mu(v_2)}\subset\P^{n-1}$$
passing through the fixed points $P_{\mu(v_1)}$ and $P_{\mu(v_2)}$.
The map $f|_{\cC_{f,e}}$ is ramified only over $P_{\mu(v_1)}$ and~$P_{\mu(v_2)}$.
In particular, $f|_{\cC_{f,e}}$ is unique up to isomorphism.
The remaining, contracted, components of $\cC_f$ are rational 
and indexed by the vertices $v\!\in\!\Ver$ such that 
$$\val(v)\equiv \big|\{e\!\in\!\Edg\!: v\!\in\!e\}\big|
+ \big|\{i\!\in\![k]\!: \eta(i)\!=\!v\}\big| \ge 3.$$
The map $f$ takes such a component $\cC_{f,v}$ to the fixed point~$P_{\mu(v)}$.
Thus,
$$\cZ_{\Ga}\approx \ov\cM_{\Ga}\!\equiv\!\prod_{v\in\Ver}\!\!\ov\cM_{0,\val(v)},$$
where $\ov\cM_{g',l}$ denotes the moduli space
of stable genus $g'$ curves with $l$ marked points.
For the purposes of this definition, $\ov\cM_{0,1}$ and $\ov\cM_{0,2}$
are one-point spaces.
For example, in the case of the first diagram in Figure~\ref{loopgraph_fig},
$$\cZ_{\Ga}\approx
\ov\cM_{\Ga} \!\equiv \ov\cM_{0,5} \!\times\! \ov\cM_{0,3}^2
\!\times\! \ov\cM_{0,2}^5 \!\times\! \ov\cM_{0,1}^4 \approx \ov\cM_{0,5}$$
is a fixed locus\footnote{after dividing by the appropriate automorphism 
group; see \cite[Section 27.3]{MirSym}}
in $\M_{1,1}^{\eff}(\P^{n-1},19)$, with $n\!\ge\!4$.
Since $n$ is fixed throughout the main computation in the paper,
each graph $\Ga$ completely determines the ambient moduli space containing 
the fixed locus~$\cZ_{\Ga}$; it will be denoted by~$\ov\M_{\Ga}$.\\

\noindent
Suppose $\Ga$ is a decorated graph as in~\e_ref{decortgraphdfn_e} and has exactly one loop.
By~\e_ref{phidfn_e} and~\e_ref{restrmap_e},
$$\ev_1^*\phi_i\big|_{\cZ_{\Ga}}
=\prod_{k\neq i}\big(\al_{\mu(\eta(1))}-\al_k\big)
=\de_{i,\mu(\eta(1))}\prod_{k\neq i}(\al_i-\al_k),$$
where $\de_{i,\mu(\eta(1))}$ is the Kronecker delta function.
Thus, by~\e_ref{ABothm_e}, $\Ga$ does not contribute to~\e_ref{Frestr_e} 
unless $\mu(\eta(1))\!=\!i$, i.e.~the marked point of the map is taken 
to the point $P_i\!\in\!\P^{n-1}$.
There are two types of  graphs that do (or may) contribute to~\e_ref{Frestr_e}; 
they will be called \sf{$A_i$} and \sf{$\ti{A}_{ij}$-types}.
In a graph of the $A_i$-type, the marked point~$1$ is attached to some vertex 
$v_0\!\in\!\Ver$ that lies inside of the loop and is labeled~$i$.
In a graph of the $\ti{A}_{ij}$-type, the marked point~$1$ is attached to 
a vertex that lies outside of the loop and is still labeled~$i$, while 
the vertex $v_0$ of the loop which is the closest to the marked point is 
labeled by some $j\!\in\![n]$.
This vertex is thus mapped to the point $P_j\!\in\!\P^{n-1}$.
Examples of graphs of the two types, with $i\!=\!3$ and $j\!=\!2$,
are depicted in Figure~\ref{loopgraph_fig}.\\

\noindent
Whether a graph $\Ga$ is of type $A_i$ or $\ti{A}_{ij}$, 
it contains a distinguished vertex~$v_0$;
it is indicated with a thick dot in Figure~\ref{loopgraph_fig}.
If we break $\Ga$ at~$v_0$, keeping a copy of~$v_0$ on each of the edges of 
$\Ga$ containing~$v_0$, and cap off each of the ``loose'' ends with a marked point
attached to~$v_0$, we obtain several decorated trees,
which will be called the \sf{strands of~$\Ga$}.
If $e_-,e_+\!\in\!\Edg$ are the two edges in the loop in $\Ga$ joined at~$v_0$,
the strands are naturally indexed by the set
$$\ov\Edg(v_0)\equiv\big\{e\!\in\!\Edg\!:v_0\!\in\!e\big\}\big/e_-\!\sim\!e_+$$
of edges leaving $v_0$, with $e_-$ and $e_+$ identified.\footnote{By 
\e_ref{decorgraphcond_e}, $e_-\!\neq\!e_+$, i.e.~there is no edge from vertex
to itself.}
The distinguished strand $\Ga_{e_{\pm}}$ with two marked points 
arising from the loop of~$\Ga$  will be denoted by~$\Ga_{\pm}$.
There are also $m\!\ge\!0$ strands, $\Ga_1,\ldots,\Ga_m$,
each of which has exactly one marked point.
Finally, if $\Ga$ is of type~$\ti{A}_{ij}$, there is also a second distinguished
strand with two marked points that contains the marked point~$1$ of~$\Ga$.
This strand will be denoted by~$\Ga_0$. 
The strands of the second graph in Figure~\ref{loopgraph_fig}
are shown in Figure~\ref{loopgraph_fig2}.\\

\noindent
The strands of a one-loop decorated graph $\Ga$ correspond to fixed loci of 
the $\T$-action on $\ov\M_{0,k}(\P^{n-1},d)$, with $k\!=\!1,2$ and $d\!\in\!\Z^+$. 
Furthermore, 
\begin{equation}\label{Alocusdecomp_e}
\cZ_{\Ga}\approx\ov\cM_{0,\val(v_0)}\times\cZ_{\Ga;B}
\equiv\ov\cM_{0,\val(v_0)}\times\prod_{e\in\ov\Edg(v_0)}\!\!\!\!\!\!\cZ_{\Ga_e},
\end{equation}
up to a quotient by a finite group, where $B$ stands for the bubble components.
If
$$\pi_P,\pi_B\!:\cZ_{\Ga}\lra\ov\cM_{0,\val(v_0)},\cZ_{\Ga;B}
\qquad\hbox{and}\qquad \pi_e\!: \cZ_{\Ga;B}\lra\cZ_{\Ga_e}$$ 
are the projection maps, then
\begin{equation}\label{Abndl_e}
\V_1|_{\cZ_{\Ga}}\approx \pi_B^*\bigg(\bigoplus_{e\in\ov\Edg(v_0)}
\!\!\!\!\!\pi_e^*\V_0'\bigg) \qquad\Lra\qquad
\E(\V_1)\big|_{\cZ_{\Ga}}=
\pi_B^*\bigg(\prod_{e\in\ov\Edg(v_0)}\!\!\!\!\!\!\pi_e^*\E(\V_0')\bigg).
\end{equation}\\

\noindent
``Most'' of the normal bundle to $\cZ_{\Ga}$ in $\wt\M_{1,1}^0(\P^{n-1},d)$ and
$\ov\M_{1,1}(\P^{n-1},d)$, as described in \cite[Section~27.4]{MirSym},
also comes from the components of~$\cZ_{\Ga}$ in the following sense.
The marked points on $\ov\cM_{0,\val(v_0)}$ are naturally indexed by the~set
$$\Edg(v_0)\equiv\big\{e\!\in\!\Edg\!:v_0\!\in\!e\big\}$$
of the edges leaving $v_0$, along with $1$ if $\Ga$ is of 
type~$A_i$.
For each $e\!\in\!\Edg(v_0)$, let 
$$L_e'\lra\ov\cM_{1,|\val(v_0)|} \qquad\hbox{and}\qquad
\hb_e\equiv c_1(L_e')\in H^2\big(\ov\cM_{1,|\val(v_0)|})$$
be the tautological tangent line bundle at the marked point 
corresponding to $e$ and its first chern class, respectively.
Analogously, let
$$L_{e_-},L_{e_+}\lra\cZ_{\Ga_{\pm}} \qquad\hbox{and}\qquad
L_e\lra\cZ_{\Ga_e},~~e\in\Edg(v_0)\!-\!\{e_-,e_+\},$$
be the restrictions to $\cZ_{\Ga_{\pm}}$ and $\cZ_{\Ga_e}$ of
the tautological line bundles on $\ov\M_{\Ga_{\pm}}$ and $\ov\M_{\Ga_v}$
at the marked points corresponding to the edges leaving~$v_0$.
Let
$$\psi_e=c_1(L_e^*)$$
be the corresponding $\psi$-classes.
The normal bundle of $\wt\cZ_{\Ga}$ in $\wt\M_{1,k}^0(\P^{n-1},d)$
is then given~by
$$\frac{\N\cZ_{\Ga}}{T_{\mu(v_0)}\P^{n-1}} =  
\pi_B^* \bigg(\bigoplus_{e\in\ov\Edg(v_0)}\!\!\!\!\!\!\pi_e^*\N\cZ_{\Ga_e}\bigg)
\oplus\bigoplus_{e\in\Edg(v_0)}\!\!
\frac{\pi_P^*L_e'\!\otimes\!\pi_B^*\pi_e^*L_e}{\pi_B^*\pi_e^*T_{\mu(v_0)}\P^{n-1}},$$
where $\N\cZ_{\Ga_e}\!\lra\!\cZ_{\Ga_e}$ is the normal bundle of $\cZ_{\Ga_e}$
in $\ov\M_{\Ga_e}$.
Thus,
\begin{equation}\label{Anormalbndl_e}
\frac{\E(T_{\mu(v_0)}\P^{n-1})}{\E(\N\cZ_{\Ga})}
 =  \prod_{e\in\ov\Edg(v_0)}\!\!\! \frac{1}{\E(\N\cZ_{\Ga_e})} \,
\prod_{e\in\Edg(v_0)}\!\!\!\!\frac{\E(T_{\mu(v_0)}\P^{n-1})}{\hb_e\!-\!\psi_e},
\end{equation}
where we omit the pullback maps $\pi_P^*$, $\pi_B^*$, and~$\pi_e^*$.\\

\noindent
By~\e_ref{Abndl_e} and~\e_ref{Anormalbndl_e}, 
$$\frac{\E(\V_1)\ev_1^*\phi_i|_{\cZ_{\Ga}}}{\E(\N\cZ_{\Ga})}\in H^*(\cZ_{\Ga})$$
splits into factors coming from the strands of $\Ga$ after integrating
over the first factor in~\e_ref{Alocusdecomp_e}.
These factors are the contributions of $\cZ_{\Ga_e}$ to integrals
over~$\ov\M_{\Ga_e}$ involving $\E(\V_0)$ and~$\psi_e$.
We will next describe the total contribution of all graphs of 
types $A_i$ and $\ti{A}_{ij}$ to~\e_ref{Frestr_e} in terms of such integrals.\\

\noindent
For all $i,j\!=\!1,2,\ldots$, let
\begin{alignat}{1}
\label{Z1ptdfn_e}
\cZ_i^*(\hb,u) &\equiv \sum_{d=1}^{\i}
\int_{\ov\M_{0,2}(\P^{n-1},d)}\frac{\E(\V_0')}{\hb\!-\!\psi_1}\ev_1^*\phi_i;\\
\label{Z1pt2dfn_e}
\cZ_{ij}^*(\hb,u) &\equiv \hb^{-1}\sum_{d=1}^{\i}
\int_{\ov\M_{0,2}(\P^{n-1},d)}\frac{\E(\V_0')}{\hb\!-\!\psi_1}
\ev_1^*\phi_i\ev_2^*\phi_j;\\
\label{Z2ptdfn_e}
\wt\cZ_{ij}^*(\hb_1,\hb_2,u) &\equiv \frac{1}{2\hb_1\hb_2}\sum_{d=1}^{\i}
\int_{\ov\M_{0,2}(\P^{n-1},d)}
\frac{\E(\V_0')}{(\hb_1\!-\!\psi_1)(\hb_2\!-\!\psi_2)}
\ev_1^*\phi_i\ev_2^*\phi_j.
\end{alignat}
These generating functions have been explicitly computed; 
see Subsection~\ref{g0form_subs}.
For the moment, we simply note that 
$$\cZ_i^*,\cZ_{ij}^*\in\Q_{\al}(\hb)[[u]]
\qquad\hbox{and}\qquad \wt\cZ_{ij}^*\in\Q_{\al}(\hb_1,\hb_2)[[u]].$$
Thus, the $\hb$-residues of these power series are well-defined.
So is
\begin{equation}\label{etadfn_e}
\eta_i(u)=\fR_{\hb=0}\Big\{\ln\big(1\!+\!\cZ_i^*(\hb,u)\big)\Big\}
\in\Q_{\al}[[u]],
\end{equation}
since the degree-zero term of the power series $\cZ_i^*(\hb,u)$ is $0$.
Let
\begin{equation}\label{Phi0dfn_e}
\Phi_0(\al_i,u)=\fR_{\hb=0}\Big\{
\hb^{-1}e^{-\eta_i(u)/\hb}\big(1\!+\!\cZ_i^*(\hb,u)\big)\Big\}
\in\Q_{\al}[[u]].
\end{equation}
By Lemma~\ref{Zreg_lmm}, $e^{-\eta_i(u)/\hb}\big(1\!+\!\cZ_i^*(\hb,u)\big)$ 
is in fact holomorphic at $\hb\!=\!0$ and 
thus $\Phi_0(\al_i,u)$ is simply the value of this power series at $\hb\!=\!0$.
Note that the degree-zero term of $\Phi_0(\al_i,u)$ is~$1$.\\

\noindent
{\it Remark:} The star in $\cZ_i^*$, $\cZ_{ij}^*$, and $\wt\cZ_{ij}^*$ 
indicates that these power series are obtained by removing the $u$-constant term
from certain natural power series $\cZ_i$, $\cZ_{ij}$, and $\wt\cZ_{ij}$;
see \cite[Chapter~29]{MirSym} and~\cite[Section~1.1]{bcov0}.

\begin{prp}
\label{contr_prp1}
(i) The total contribution $\cA_i(u)$ to~\e_ref{Frestr_e} from all graphs of type $A_i$ 
is given~by
\begin{equation}\label{contr_prp1e1}
\cA_i(u)=\frac{1}{\Phi_0(\al_i,u)}\fR_{\hb_1=0}\Big\{\fR_{\hb_2=0}\Big\{
e^{-\eta_i(u)/\hb_1}e^{-\eta_i(u)/\hb_2}\wt\cZ_{ii}^*(\hb_1,\hb_2,u)\Big\}\Big\};
\end{equation}
(ii) The total contribution $\ti\cA_{ij}(u)$ to~\e_ref{Frestr_e} from all graphs 
of type $\ti{A}_{ij}$ is
\begin{equation}\label{contr_prp1e2}
\ti\cA_{ij}(u)=\frac{\cA_j(u)}{\prod_{k\neq j}(\al_j\!-\!\al_k)}
\fR_{\hb=0}\Big\{e^{-\eta_j(u)/\hb}\cZ_{ji}^*(\hb,u)\Big\}.
\end{equation}\\
\end{prp}

\noindent
It is fairly straightforward to express $\cA_i(u)$ and $\ti\cA_{ij}(u)$
in terms of sums of products of the residues of $\cZ_i^*$, $\cZ_{ij}^*$, 
and $\wt\cZ_{ij}^*$ at every possible $\hb\!\in\!\C^*$.
A straightforward application of the Residue Theorem on~$S^2$
then reduces the resulting expressions to sums of products of the residues of $\cZ_i^*$, $\cZ_{ij}^*$, and $\wt\cZ_{ij}^*$ at $\hb\!=\!0$.
However, the products will have either $m\!+\!2$ or $m\!+\!3$ factors,
where $m$ is the number of strands of $\Ga$ with one marked point,
and must be summed over all possible~$m$. 
We are able to sum over~$m$ because $e^{-\eta_i(u)/\hb}\big(1\!+\!\cZ_i^*(\hb,u)\big)$
turns out to be holomorphic at~$\hb\!=\!0$; see Lemmas~\ref{regular_lmm} and~\ref{Zreg_lmm}.
Proposition~\ref{contr_prp1} is proved in Subsection~\ref{localprppf_subs}.

\subsection{Contributions from Fixed Loci, II}
\label{localdata_subs2}

\noindent
In this subsection we describe the contribution to~\e_ref{Frestr_e} from
the fixed loci of the $\T$-action on $\wt\M_{1,1}^0(\P^{n-1},d)$ that are contained
in $\partial\wt\M_{1,1}^0(\P^{n-1},d)$.
We begin by reviewing the description of such loci and their normal bundles given in
\cite[Subsection~1.4]{VaZ}.\\

\begin{figure}
\begin{pspicture}(0,-1.5)(10,2)
\psset{unit=.4cm}
\pscircle*(10,0){.3}\rput(10,.7){\smsize{$2$}}
\psline[linewidth=.1](10,0)(7,2)\pscircle*(7,2){.2}
\rput(8.6,1.5){\smsize{$2$}}\rput(7.3,2.5){\smsize{$3$}}
\psline[linewidth=.1](10,0)(6,0)\pscircle*(6,0){.2}
\rput(8.5,.4){\smsize{$2$}}\rput(6.3,.6){\smsize{$3$}}
\psline[linewidth=.04](6,0)(3,2)\pscircle*(3,2){.2}
\rput(4.6,1.4){\smsize{$3$}}\rput(3.4,2.5){\smsize{$2$}}
\psline[linewidth=.04](6,0)(3,-2)\pscircle*(3,-2){.2}
\rput(4.1,-.8){\smsize{$1$}}\rput(3.2,-1.3){\smsize{$1$}}
\psline[linewidth=.1](10,0)(7,-1.5)\pscircle*(7,-1.5){.2}
\rput(7.7,-.7){\smsize{$2$}}\rput(6.8,-.9){\smsize{$3$}}
\psline[linewidth=.04](10,0)(8,-3)\pscircle*(8,-3){.2}
\rput(8.5,-1.6){\smsize{$2$}}\rput(7.5,-3.3){\smsize{$1$}}
\psline[linewidth=.04](10,0)(12,-2)\pscircle*(12,-2){.2}
\rput(11,-1.5){\smsize{$3$}}\rput(12.5,-2.2){\smsize{$3$}}
\psline[linewidth=.05,linestyle=dashed](10,0)(14,0)\pscircle*(14,0){.2}
\psline[linewidth=.04](14,0)(17,-1)\pscircle*(17,-1){.2}
\rput(16.2,-.3){\smsize{$2$}}\rput(17.5,-1){\smsize{$1$}}
\psline[linewidth=.04](14,0)(17,1)\pscircle*(17,1){.2}
\rput(15.6,.9){\smsize{$3$}}\rput(17.2,1.6){\smsize{$3$}}
\psline[linewidth=.04](17,1)(20,2)\pscircle*(20,2){.2}
\rput(18.6,1.9){\smsize{$1$}}\rput(20.2,2.6){\smsize{$1$}}
\psline[linewidth=.05,linestyle=dashed](10,0)(13,2)\pscircle*(13,2){.2}
\psline[linewidth=.04](13,2)(16,3)\pscircle*(16,3){.2}
\rput(14.6,3){\smsize{$1$}}\rput(16.2,3.6){\smsize{$3$}}
\psline[linewidth=.04](13,2)(14,4)\rput(14.3,4.4){\smsize{$\bf 1$}}
\rput(7.5,0.5){\smsize{$e_1$}}\rput(11.5,1.7){\smsize{$e_2$}}
\pscircle*(30,-2){.2}\rput(29.8,-2.6){\smsize{$2$}}
\psline[linewidth=.05](30,-2)(34,-2)\pscircle*(34,-2){.2}
\rput(32.5,-2.4){\smsize{$1$}}\rput(34,-2.6){\smsize{$3$}}
\psline[linewidth=.05](30,-2)(32,-4)\rput(32.4,-4.1){\smsize{$\bf 1$}}
\psline[linewidth=.05](30,-2)(32,0)\rput(32.4,-.1){\smsize{$\bf 2$}}
\rput(27.2,-2.1){$\Ga_{e_2}\!=$}
\pscircle*(30,2){.2}\rput(29.8,1.4){\smsize{$2$}}
\psline[linewidth=.05](30,2)(34,2)\pscircle*(34,2){.2}
\rput(32.5,1.6){\smsize{$2$}}\rput(34,1.4){\smsize{$3$}}
\psline[linewidth=.05](34,2)(37,3)\pscircle*(37,3){.2}
\rput(35,2.8){\smsize{$3$}}\rput(37.5,3){\smsize{$2$}}
\psline[linewidth=.05](34,2)(37,1)\pscircle*(37,1){.2}
\rput(35.5,1){\smsize{$1$}}\rput(37.5,1){\smsize{$1$}}
\psline[linewidth=.05](30,2)(32,3.5)\rput(32.4,3.4){\smsize{$\bf 1$}}
\rput(27.2,1.9){$\Ga_{e_1}\!=$}
\end{pspicture}
\caption{A refined decorated rooted tree and some of its strands}
\label{reftree_fig}
\end{figure}
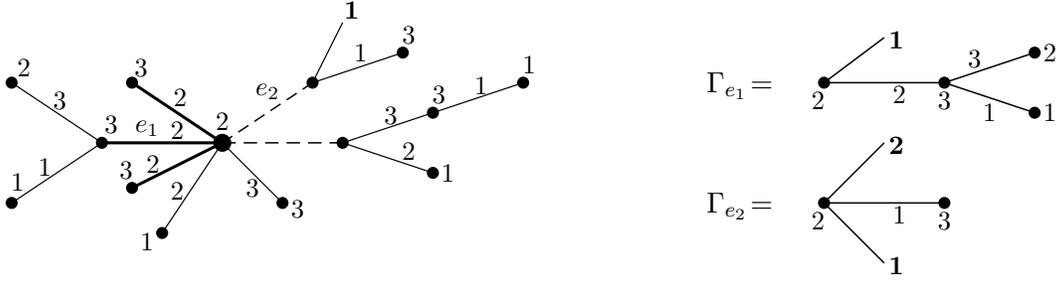

\noindent
A \sf{rooted tree} is a tree (i.e.~a graph with no loops) 
with a distinguished vertex.
A tuple 
$$(\Ver,\Edg,v_0;\Ver_+,\Ver_0)$$
is a \sf{refined rooted tree} if $(\Ver,\Edg,v_0)$ is a rooted tree,
i.e.~$v_0$ is the distinguished vertex of the tree~$(\Ver,\Edg)$, and 
$$\Ver_+,\Ver_0\subset\Ver\!-\!\{v_0\}, \quad \Ver_+\neq\eset, \quad
\Ver_+\!\cap\!\Ver_0=\eset, \quad
\{v_0,v\}\!\in\!\Edg ~\forall\, v\!\in\!\Ver_+\!\cup\!\Ver_0.$$
Given such a refined rooted tree, we put
$$\Edg_+=\big\{\{v_0,v\}\!: v\!\in\!\Ver_+\big\} \qquad\hbox{and}\qquad
\Edg_0=\big\{\{v_0,v\}\!: v\!\in\!\Ver_0\big\}.$$
In the first diagram of Figure~\ref{reftree_fig}, the distinguished vertex $v_0$
is indicated by the thick dot.
The elements of $\Edg_+$ and $\Edg_0$ are shown
as the thick solid lines and the thin dashed lines, respectively.\\

\noindent
A \sf{refined decorated rooted tree} is a tuple
\begin{equation}\label{treedfn_e}
\Ga = \big(\Ver,\Edg,v_0;\Ver_+,\Ver_0;\mu,\d,\eta \big),
\end{equation}
where $(\Ver,\Edg;v_0;\Ver_+,\Ver_0)$ is a refined rooted tree
and
$$\mu\!:\Ver\!-\!\Ver_0\lra[n], \qquad
\d\!: \Edg\!-\!\Edg_0\lra\Z^+, \quad\hbox{and}\quad \eta\!: [k]\lra\Ver$$
are maps such that\\
${}\quad$ (i) $\mu(v_1)\!=\!\mu(v_2)$ and $\d(\{v_0,v_1\})\!=\!\d(\{v_0,v_2\})$
for all $v_1,v_2\!\in\!\Ver_+$;\\
${}\quad$ (ii) if $v_1\!\in\!\Ver_+$, $v_2\!\in\!\Ver\!-\!\Ver_0\!-\!\Ver_+$,
and $\{v_0,v_2\}\!\in\!\Edg$, then 
\begin{equation}\label{notthickedges_e}
\mu(v_1)\neq\mu(v_2) \qquad\hbox{or}\qquad
\d(\{v_0,v_1\})\!\neq\!\d(\{v_0,v_2\});
\end{equation}
${}\quad$ (iii) if $\{v_1,v_2\}\!\in\!\Edg$ and $v_2\!\not\in\!\Ver_0\!\cup\!\{v_0\}$,
then 
$$\mu(v_2)\!\neq\!\mu(v_1)  \quad\hbox{if}~~~ v_1\!\not\in\!\Ver_0
\qquad\hbox{and}\qquad 
\mu(v_2)\!\neq\!\mu(v_0) \quad\hbox{if}~~~ v_1\!\in\!\Ver_0;$$
${}\quad$ (iv) if $v_1\!\in\!\Ver_0$, then $\{v_1,v_2\}\!\in\!\Edg$ for some
$v_2\!\in\!\Ver\!-\!\{v_0\}$ and $\val(v_1)\!\ge\!3$.\\
In Figure~\ref{reftree_fig}, the value of the map $\mu$ on each vertex,
not in $\Ver_0$, is indicated by the number next to the vertex.
Similarly, the value of the map $\d$ on each edge,
not in $\Edg_0$, is indicated by the number next to the edge.
The elements of the~set $[k]\!=\![1]$ are shown in bold face.
Each of them is linked by a line segment to its image under~$\eta$.
The first condition above implies that all of the thick edges have the same labels,
and so do their vertices, other than the root~$v_0$.
By the second condition, the set of thick edges is 
a maximal set of edges leaving~$v_0$ which satisfies the first condition.
By the third condition, no two consecutive vertex labels are the same.
The final condition implies that there are at least two solid lines, at least one of which
is an edge, leaving from every vertex which is connected to the root by a dashed line.\\

\noindent
{\it Remark:} In \cite[Subsection~1.4]{VaZ},
refined decorated rooted trees are required to satisfy a fifth condition,
$\sum_{e\in\Edg_+}\!\!\d(e)\ge2$.
If this condition is not satisfied (i.e.~$\d(e)\!=\!1$ for the unique element $e\!\in\!\Edg_+$),
the locus $\ti\cZ_{\Ga}$ corresponding to~$\Ga$ via the construction below
will not contribute to~\e_ref{Frestr_e}; see also the next remark.\\

\noindent
Let $\Ga$ be a refined decorated rooted tree as in~\e_ref{treedfn_e}.
Breaking $\Ga$ at $v_0$, we split $\Ga$ into pieces $\Ga_e'$ indexed by the set
$$\Edg(v_0)\equiv \big\{e\!\in\!\Edg\!: v_0\!\in\!e\big\}$$
of the edges in $\Ga$ leaving $v_0$.
If $e\!\not\in\!\Edg_0$, we will keep the vertex $v_0$ and the edge $e$
and cap $\Ga_e'$ off with a new marked point attached to~$v_0$, 
just as in Subsection~\ref{localdata_subs1}.
If $e\!=\!\{v_0,v\}$ with $v\!\in\!\Ver_0$, we remove $v_0$ and $e$ from~$\Ga_e'$,
cap it off with a new marked point attached to~$v$, and assign 
the $\mu$-value of~$v_0$ to~$v$.
In either case, we denote the resulting decorated tree by $\Ga_e$ and 
call it a \sf{strand} of~$\Ga$; see Figure~\ref{reftree_fig} for two examples.
If $v\!\in\!\Ver_+$, let
$$\mu(\Ga)=\mu(v) \qquad\hbox{and}\qquad \d(\Ga)=\d\big(\{v_0,v\}\big).$$
By the requirement (i) on $\Ga$, $\mu(\Ga)$ and $\d(\Ga)$ do not depend 
on the choice of $v\!\in\!\Ver_+$.\\

\noindent
Via the construction of Subsection~\ref{localdata_subs1}, each strand $\Ga_e$ of $\Ga$
determines a $\T$-fixed locus $\cZ_{\Ga_e}$ in a moduli space $\ov\M_{\Ga_e}$
of genus~$0$ stable maps.
Let 
$$\cZ_{\Ga;B}=\prod_{e\in\Edg(v_0)}\!\!\!\!\!\!\cZ_{\Ga_e},$$
where $B$ stands for the ``bubble'' components.
Denote by
$$\pi_e\!: \cZ_{\Ga;B}\!\lra\!\cZ_{\Ga_e} \qquad\hbox{and}\qquad
L_e\lra\cZ_{\Ga_e}$$ 
the natural projection map and the restriction  to~$\cZ_{\Ga_e}$ of 
the tautological tangent line bundle on $\ov\M_{\Ga_e}$ for 
the marked point corresponding to the attachment at~$v_0$.
Let
\begin{gather*}
L_{\Ga}=\pi_e^*L_e\lra\cZ_{\Ga;B} ~~\hbox{if}~e\!\in\!\Edg_+, 
\qquad \psi_{\Ga}=c_1\big(L_{\Ga}^*\big)\in H^2(\cZ_{\Ga;B}),\\
F_{\Ga;B}=\bigoplus_{e\in\Edg_+}\!\!\!\pi_e^*L_e,\quad
F_{\Ga;B}^c=\bigoplus_{e\in\Edg(v_0)-\Edg_+}\!\!\!\!\!\!\!\!\!\!\!\!\pi_e^*L_e,
\quad \wt\cZ_{\Ga;B}=\P F_{\Ga;B}.
\end{gather*}
By the requirement (i) on $\Ga$, $L_{\Ga}$ is well-defined as a $\T$-equivariant line
bundle and  
$$\wt\cZ_{\Ga;B}\approx \cZ_{\Ga;B}\times \P^{|\Edg_+|-1}.$$
Let
$$\pi_1,\pi_2\!: \wt\cZ_{\Ga}\lra\cZ_{\Ga},\P^{|\Edg_+|-1}$$
be the two projection maps.
Up to a quotient by a finite group, the fixed locus of the $\T$-action on 
$\wt\M_{1,k}^0(\Pn,d)$ corresponding to $\Ga$
(or its equivalent for our purposes; see the next remark) is
\begin{equation}\label{Blocusdecomp_e}
\wt\cZ_{\Ga}\equiv \wt\cM_{1,|\val(v_0)|}\!\times\ti\cZ_{\Ga;B}
\approx \wt\cM_{1,|\val(v_0)|}\!\times \cZ_{\Ga;B}\times \P^{|\Ver_+|-1},
\end{equation}
where $\wt\cM_{1,|\val(v_0)|}$ is a certain blowup of the moduli space
of genus-one curves with $|\val(v_0)|$ marked points
constructed in \cite[Subsection~2.3]{VaZ}.\footnote{\label{curvebl_ftnt}
The blowups $\wt\cM_{1,(I,J)}$ of $\ov\cM_{1,N}$ constructed
in \cite{VaZ} are indexed by ordered partitions $(I,J)$ of~$[N]$.
The only cases encountered as components of the fixed loci of $\wt\M_{1,1}^0(\P^{n-1},d)$
are $|J|\!=\!0,1$.
In these two cases, the blowups are the same, 
as are the universal tangent line bundles~$\L$ appearing in the following paragraph.} 
The only property of $\wt\cM_{1,|\val(v_0)|}$ relevant for the purposes of
this paper is~\e_ref{Lnum_e} below.
Denote by
$$\pi_P,\pi_B\!:\wt\cZ_{\Ga}\lra \wt\cM_{1,|\val(v_0)|},\ti\cZ_{\Ga;B}$$
the projection maps.\\

\noindent
The normal bundle of $\wt\cZ_{\Ga}$ in $\ov\M_{1,k}^0(\P^{n-1},d)$ (or its equivalent) 
is given~by
\begin{equation}\label{Bnormbdl_e1}\begin{split}
\N\wt\cZ_{\Ga}=\pi_B^*\Big( \big(\pi_1^*\N\cZ_{\Ga;B}\oplus 
\pi_1^*(L_{\Ga}^*\!\otimes\!F_{\Ga;B}^c)\!\otimes\!\pi_2^*\ga^*\big)
\big/\big(\pi_1^*(L_{\Ga}^*\!\otimes\!T_{\mu(v_0)}\P^{n-1})\!\otimes\!\pi_2^*\ga^*
\big)\Big)&\\
\oplus\pi_P^*\L\!\otimes\!\pi_B^*\big(\pi_1^*L_{\Ga}\!\otimes\!\pi_2^*\ga\big)&,
\end{split}\end{equation}
where $\ga\!\lra\!\P^{|\Ver_+|-1}$ is the tautological line bundle,
\begin{equation}\label{Bnormbdl_e3}
\frac{\E\big(\N\cZ_{\Ga;B}\big)}{\E(T_{\mu(v_0)}\P^{n-1})}
=\prod_{e\in\Edg(v_0)}\! 
\pi_e^*\bigg(\frac{\E\big(\N\cZ_{\Ga_e}\big)}{\E(T_{\mu(v_0)}\P^{n-1})}\bigg)
\in\H_{\T}^*\big(\cZ_{\Ga;B}\big),\footnotemark
\end{equation}
and $\L\!\lra\!\wt\cM_{1,|\val(v_0)|}$ is the universal tangent line bundle
constructed in  \cite[Subsection~2.3]{VaZ}.
\footnotetext{\label{tuplemod_ftnt}The vector bundle
$\N\cZ_{\Ga;B}\!\lra\!\cZ_{\Ga;B}$ is the normal bundle of $\cZ_{\Ga;B}$ 
in the moduli space $\ov\M_{\Ga;B}$ of $|\Edg(v_0)|$ tuples 
of genus-zero stable maps that agree at a distinguished marked point of 
each element of the tuple.}
The only property of this line bundle relevant for our purposes~is
\begin{equation}\label{Lnum_e}
\int_{\wt\cM_{1,|\val(v_0)|}}\!\!\!\ti\psi^{\,|\val(v_0)|}=\frac{(|\val(v_0)|\!-\!1)!}{24}
\qquad\hbox{if}~
k\!=\!1,\footnote{This assumption on $k$ implies that $|J|\!=\!0,1$; see
Footnote~\ref{curvebl_ftnt}.}
\end{equation}
where $\ti\psi\!=\!c_1(\L^*)$ is the \sf{universal $\psi$-class};
see \cite[Corollary~1.2]{g1desing2}.\\

\noindent
The final piece of localization data we need to recall from~\cite{VaZ} is that
\begin{equation}\label{Bbdl_e1}
\E(\V_1')\big|_{\wt\cZ_{\Ga}}= \pi_B^*\big(\V_{\Ga;B}'/
\pi_1^*(L_{\Ga}^*\!\otimes\!\cL_{\mu(v_0)})\!\otimes\!\pi_2^*\ga^*\big),
\end{equation}
where 
\begin{equation}\label{Bbdl_e3}
\V_{\Ga;B}'=\prod_{e\in\Edg(v_0)}\! \pi_e^*\E(\V_0').
\footnote{The vector bundle $\V_{\Ga;B}'$ is the analogue of the vector bundle $\V_0'$ 
for the moduli space $\ov\M_{\Ga;B}$; see Footnote~\ref{tuplemod_ftnt}.
It is obtained by pulling $\cL$ back to the universal curve, then pushing down,
and then taking the kernel of a natural evaluation map.}
\end{equation}\\

\noindent
{\it Remark:} If $\d(\Ga)\!\ge\!2$, the pair $(\wt\cZ_{\Ga},\N\wt\cZ_{\Ga})$
described above is precisely the fixed locus corresponding to~$\Ga$ and its normal
bundle as described in  \cite[Subsection~1.4]{VaZ}.
If $\d(\Ga)\!=\!1$, the actual fixed locus corresponding to~$\Ga$ has 
$\P^{|\Ver_+|-2}$ instead of $\P^{|\Ver_+|-1}$ as the last factor in~\e_ref{Blocusdecomp_e}
and $F_{\Ga;B}^c$ has an extra component of~$L_{\Ga}$.
Thus, by~\e_ref{Bnormbdl_e1}, the expression
$$\int_{\wt\cZ_{\Ga}}\frac{(\E(\V_1)\ev_1^*\phi_i)|_{\cZ_{\Ga}}}{\E(\N\wt\cZ_{\Ga})}$$
as described above agrees with the correct one, as the extra dimension in 
the last factor in~\e_ref{Blocusdecomp_e} is canceled by 
the extra factor of~$c_1(\ga^*)$ in the integrand.\footnote{If $\d(\Ga)\!=\!1$, 
$L_{\Ga}$ is a direct summand in $T_{\mu(v_0)}\P^{n-1}$. 
The extra factor of $c_1(\ga^*)$ in the integrand comes from 
the direct summand $\ga^*$ of the bundle $\pi_1^*(L_{\Ga}^*\!\otimes\!T_{\mu(v_0)}\P^{n-1})\!\otimes\!\pi_2^*\ga^*$
in~\e_ref{Bnormbdl_e1}. 
This summand is canceled by the extra summand of $L_{\Ga}$ in $F_{\Ga;B}^c$ in~\cite{VaZ}.}\\

\begin{figure}
\begin{pspicture}(0,-1.3)(10,2)
\psset{unit=.4cm}
\pscircle*(10,0){.3}\rput(10,.7){\smsize{$2$}}
\psline[linewidth=.1](10,0)(7,2)\pscircle*(7,2){.2}
\rput(8.6,1.5){\smsize{$2$}}\rput(7.3,2.5){\smsize{$3$}}
\psline[linewidth=.1](10,0)(6,0)\pscircle*(6,0){.2}
\rput(7.7,.4){\smsize{$2$}}\rput(6.3,.6){\smsize{$3$}}
\psline[linewidth=.04](6,0)(3,2)\pscircle*(3,2){.2}
\rput(4.6,1.4){\smsize{$3$}}\rput(3.4,2.5){\smsize{$2$}}
\psline[linewidth=.04](6,0)(3,-2)\pscircle*(3,-2){.2}
\rput(4.1,-.8){\smsize{$1$}}\rput(3.2,-1.3){\smsize{$1$}}
\psline[linewidth=.1](10,0)(7,-1.5)\pscircle*(7,-1.5){.2}
\rput(7.7,-.7){\smsize{$2$}}\rput(6.8,-.9){\smsize{$3$}}
\psline[linewidth=.04](10,0)(8,-3)\pscircle*(8,-3){.2}
\rput(8.5,-1.6){\smsize{$2$}}\rput(7.5,-3.3){\smsize{$1$}}
\psline[linewidth=.04](10,0)(12,-2)\pscircle*(12,-2){.2}
\rput(11,-1.5){\smsize{$3$}}\rput(12.5,-2.2){\smsize{$3$}}
\psline[linewidth=.05,linestyle=dashed](10,0)(14,0)\pscircle*(14,0){.2}
\psline[linewidth=.04](14,0)(17,-1)\pscircle*(17,-1){.2}
\rput(16.2,-.3){\smsize{$2$}}\rput(17.5,-1){\smsize{$1$}}
\psline[linewidth=.04](14,0)(17,1)\pscircle*(17,1){.2}
\rput(15.6,.9){\smsize{$3$}}\rput(17.2,1.6){\smsize{$3$}}
\psline[linewidth=.04](17,1)(20,2)\pscircle*(20,2){.2}
\rput(18.6,1.9){\smsize{$1$}}\rput(20.2,2.6){\smsize{$1$}}
\psline[linewidth=.04](10,0)(13,2)\pscircle*(13,2){.2}
\rput(11.3,1.4){\smsize{$1$}}\rput(13.2,2.7){\smsize{$3$}}
\psline[linewidth=.04](10,0)(10,-2.7)\rput(10,-3.2){\smsize{$\bf 1$}}
\pscircle*(29,0){.3}\rput(29,.7){\smsize{$2$}}
\psline[linewidth=.1](29,0)(26,2)\pscircle*(26,2){.2}
\rput(27.6,1.5){\smsize{$2$}}\rput(26.3,2.5){\smsize{$3$}}
\psline[linewidth=.1](29,0)(25,0)\pscircle*(25,0){.2}
\rput(26.7,.4){\smsize{$2$}}\rput(25.3,.6){\smsize{$3$}}
\psline[linewidth=.04](25,0)(22,2)\pscircle*(22,2){.2}
\rput(23.6,1.4){\smsize{$3$}}\rput(22.4,2.5){\smsize{$2$}}
\psline[linewidth=.04](25,0)(22,-2)\pscircle*(22,-2){.2}
\rput(23.1,-.8){\smsize{$1$}}\rput(22.2,-1.3){\smsize{$1$}}
\psline[linewidth=.1](29,0)(26,-1.5)\pscircle*(26,-1.5){.2}
\rput(26.7,-.7){\smsize{$2$}}\rput(25.8,-.9){\smsize{$3$}}
\psline[linewidth=.04](29,0)(27,-3)\pscircle*(27,-3){.2}
\rput(27.5,-1.6){\smsize{$2$}}\rput(26.5,-3.3){\smsize{$1$}}
\psline[linewidth=.04](29,0)(31,-2)\pscircle*(31,-2){.2}
\rput(30,-1.5){\smsize{$3$}}\rput(31.5,-2.2){\smsize{$3$}}
\psline[linewidth=.05,linestyle=dashed](29,0)(33,0)\pscircle*(33,0){.2}
\psline[linewidth=.04](33,0)(36,-1)\pscircle*(36,-1){.2}
\rput(35.2,-.3){\smsize{$2$}}\rput(36.5,-1){\smsize{$1$}}
\psline[linewidth=.04](33,0)(36,1)\pscircle*(36,1){.2}
\rput(34.6,.9){\smsize{$3$}}\rput(36.2,1.6){\smsize{$3$}}
\psline[linewidth=.04](36,1)(39,2)\pscircle*(39,2){.2}
\rput(37.6,1.9){\smsize{$1$}}\rput(39.2,2.6){\smsize{$1$}}
\psline[linewidth=.05,linestyle=dashed](29,0)(32,2)\pscircle*(32,2){.2}
\psline[linewidth=.04](32,2)(35,3)\pscircle*(35,3){.2}
\rput(33.6,3){\smsize{$1$}}\rput(35.2,3.6){\smsize{$3$}}
\psline[linewidth=.04](32,2)(33,4)\rput(33.3,4.4){\smsize{$\bf 1$}}
\end{pspicture}
\caption{Refined decorated rooted trees of types $B_2$ and $\ti{B}_{22}$}
\label{treegraps_fig}
\end{figure}
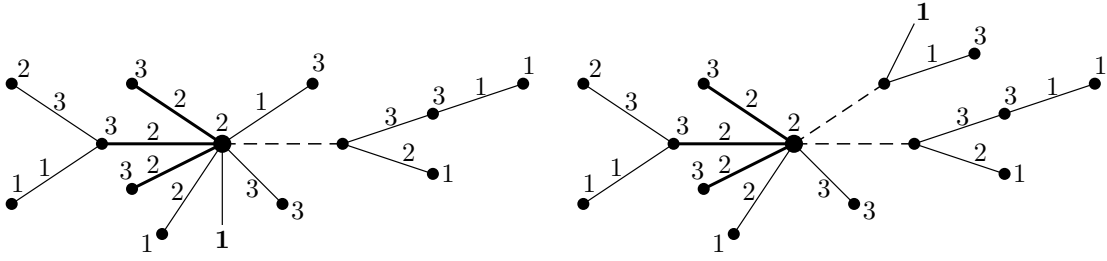

\noindent
We now consider the refined decorated rooted trees $\Ga$ 
as in~\e_ref{treedfn_e} that contribute to~\e_ref{Frestr_e}.
As in Subsection~\ref{localdata_subs1}, $\Ga$ does not contribute to~\e_ref{Frestr_e}
unless $\mu(\eta(1))\!=\!i$, i.e.~the marked point is mapped to $P_i\!\in\!\P^{n-1}$.
Similarly to Subsection~\ref{localdata_subs1}, we group all graphs that 
contribute to~\e_ref{Frestr_e} into two types: $B_i$ and~$\ti{B}_{ij}$.
In the graphs of type~$B_i$, $\eta(1)\!=\!v_0$, i.e.~the marked point~$1$
lies on the principal contracted component of the domain of the maps
(and is mapped to~$P_i$).
In the graphs of type~$B_{ij}$, $\eta(1)\!\neq\!v_0$, i.e.~the marked point~$1$
lies on one of the strands of~$\Ga$, while $\mu(v_0)\!=\!j$;
see Figure~\ref{treegraps_fig} for examples.
A graph of type~$B_i$ has $m\!\ge\!1$ strands with one marked point.
On the other hand, a graph of type $\ti{B}_{ij}$ has a distinguished strand
with two marked points and $m\!\ge\!0$ strands with one marked point.
In either case, the first factor in~\e_ref{Blocusdecomp_e} is $\wt\cM_{1,m+1}$.\\

\noindent
By~\e_ref{Bbdl_e1}, \e_ref{Bbdl_e3}, and~\e_ref{eulerclass_e},
\begin{equation}\label{Bbndl_e}
\E(\V_1)|_{\wt\cZ_{\Ga}}
=n\al_{\mu(v_0)}\prod_{e\in\Edg(v_0)}\!\!\!\!\!\!\!\pi_e^*\E(\V_0')
\Big/\big(n\al_{\mu(v_0)}\!+\!\psi_{\Ga}\!+\!\la\big),
\end{equation}
where $\la\!=\!c_1(\ga^*)$ and we omit the pullback maps $\pi_P^*$ and $\pi_B^*$.
By~\e_ref{Bnormbdl_e1} and~\e_ref{Bnormbdl_e3}, 
\begin{equation}\label{Bnormalbndl_e}
\frac{\E(T_{\mu(v_0)}\P^{n-1})}{\E\big(\N\wt\cZ_{\Ga}\big)}
=-\frac{\E(\ga^*\!\otimes\!L_{\Ga}^*\!\otimes\!T_{\mu(v_0)}\P^{n-1})}
{\psi_{\Ga}\!+\!\ti\psi\!+\la}
\!\!\!\!\!\prod_{e\in\Edg(v_0)}\!\frac{\E(T_{\mu(v_0)}\P^{n-1})}{\E\big(\N\cZ_{\Ga_e}\big)}
\prod_{e\in\Edg(v_0)-\Edg_+}\!\!\frac{1}{\psi_{\Ga}\!-\!\psi_e\!+\!\la}.
\end{equation}
Thus,
$$\frac{\E(\V_1)\ev_1^*\phi_i|_{\wt\cZ_{\Ga}}}{\E(\N\wt\cZ_{\Ga})}
\in H^*(\wt\cZ_{\Ga})$$
splits into factors coming from the strands of $\Ga$ after integrating
over the first and last factors in~\e_ref{Blocusdecomp_e}.
These factors are the contributions of $\cZ_{\Ga_e}$ to  integrals
over~$\ov\M_{\Ga_e}$ involving $\E(\V_0)$ and~$\psi_e$.

\begin{prp}
\label{contr_prp2}
(i) The total contribution $\cB_i(u)$ to~\e_ref{Frestr_e} from all graphs of type $B_i$ is
\begin{equation}\label{contr_prp2e1}
\cB_i(u)=\frac{n\al_i}{24} \fR_{\hb=0,\i,-n\al_i}\bigg\{
\frac{\prod_{k=1}^{k=n}(\al_i\!-\!\al_k\!+\!\hb)}{(n\al_i\!+\!\hb)\hb^3}
\frac{\cZ_i^*(\hb,u)}{1\!+\!\cZ_i^*(z,u)}\bigg\}.
\end{equation}
(ii) The total contribution $\ti\cB_{ij}(u)$ to~\e_ref{Frestr_e} from all graphs of type 
$\ti{B}_{ij}$ is
\begin{equation}\label{contr_prp2e2}
\ti\cB_{ij}(u)=-\frac{1}{\prod_{k\neq j}(\al_j\!-\!\al_k)}
\frac{n\al_j}{24} \fR_{\hb=0,\i,-n\al_j}\bigg\{
\frac{\prod_{k=1}^{k=n}(\al_j\!-\!\al_k\!+\!\hb)}{(n\al_j\!+\!\hb)\hb^2}
\frac{\cZ_{ji}^*(\hb,u)}{1\!+\!\cZ_j^*(z,u)}\bigg\}.
\end{equation}\\
\end{prp}

\noindent
The proof of this proposition turns out to be quite a bit simpler than that of
Proposition~\ref{contr_prp1}.
At an early stage in the computation, Lemma~\ref{prodres_lmm} reduces a sum of $m$~products 
of residues to the residue of a product.
The resulting products  sum over~$m$ to the functions of $\hb$ appearing 
in the statement of Proposition~\ref{contr_prp2}.
The residues of these functions on $S^2\!-\!\{0,\i,-n\al_j\}$ are summed,
using the Residue Theorem on~$S^2$, to get $\cB_i(u)$ and~$\cB_{ij}(u)$.
Proposition~\ref{contr_prp2} is proved in Subsection~\ref{localprppf_subs}.

\section{Localization Computations} 
\label{mainlocal_sec}

\noindent
In Subsection~\ref{regul_subs} we give two equivalent characterizations
of a property of power series in rational functions that reduces
infinite summations involving certain products of residues of such
power series to simple expressions.
In Subsection~\ref{Zreg_subs}, we dig deeper into Givental's proof of~\e_ref{g0ms_e}
to show that a certain generating function for genus~$0$ GW-invariants satisfies
this property.
We use these observations to prove Proposition~\ref{contr_prp1} in 
Subsection~\ref{localprppf_subs}, along with Proposition~\ref{contr_prp2}.

\subsection{Regularizable Power Series in Rational Functions}
\label{regul_subs}

\begin{dfn}
\label{regular_dfn}
A power series $\cZ^*\!=\!\cZ^*(\hb,u)\!\in\!\Q_{\al}(\hb)[[u]]$ 
is \sf{regularizable at} $\hb\!=\!0$ if there exist power series
$$\eta\!=\!\eta(u)\in\Q_{\al}[[u]] \qquad\hbox{and}\qquad
\bar\cZ^*\!=\!\bar\cZ^*(\hb,u)\in\Q_{\al}(\hb)[[u]]$$
with no degree-zero term such that $\bar\cZ^*$ is regular at $\hb\!=\!0$ and
\begin{equation}\label{regdfn_e}
1+\cZ^*(\hb,u)=e^{\eta(u)/\hb}\big(1+\bar\cZ^*(\hb,u)\big).
\end{equation}\\
\end{dfn}

\noindent
If $\cZ^*$ is regularizable at $\hb\!=\!0$, $\cZ^*$ has no degree-zero term and
the \sf{regularizing pair} $(\eta,\bar\cZ^*)$ is unique.
It is determined~by
\begin{equation}\label{regprp_e}
\eta(u)=\fR_{\hb=0}\Big\{\ln\big(1\!+\!\cZ^*(\hb,u)\big)\Big\}.
\end{equation}
The logarithm above is a well-defined power series in $u$,
since $\cZ^*(\hb,u)$ has no degree-zero term.

\begin{lmm}
\label{regular_lmm}
Suppose $\cZ^*\!=\!\cZ^*(\hb,u)\!\in\!\Q_{\al}(\hb)[[u]]$ has no degree-zero term.\\
(i) The power series $\cZ^*$ is regularizable at $\hb\!=\!0$ if and only if 
for every $a\!\ge\!0$
\begin{equation}\label{regprp_e1}
\sum_{m=2}^{\i} \frac{1}{m(m\!-\!1)}
\sum_{\underset{a_l\ge0}{\sum\limits_{l=1}\limits^{l=m}\!\!a_l=m-2-a}}
\!\!\!\!\Bigg(\prod_{l=1}^{l=m}\frac{(-1)^{a_l}}{a_l!}
\fR_{\hb=0}\Big\{\hb^{-a_l}\cZ^*(\hb,u)\Big\}\Bigg)
=a!\,\fR_{\hb=0}\Big\{\hb^{a+1}\cZ^*(\hb,u)\Big\}.
\end{equation}
(ii) If $(\eta,\bar\cZ^*)$ is the regularizing pair for $\cZ^*$ at $\hb\!=\!0$, 
then for every $a\!\ge\!0$
\begin{equation}\label{regprp_e2}
\sum_{m=0}^{\i} \sum_{\underset{a_l\ge0}{\sum\limits_{l=1}\limits^{l=m}\!\!a_l=m-a}}
\!\!\!\Bigg(\prod_{l=1}^{l=m}\frac{(-1)^{a_l}}{a_l!}
\fR_{\hb=0}\Big\{\hb^{-a_l}\cZ^*(\hb,u)\Big\}\Bigg)
=\frac{\eta(u)^a}{1+\bar\cZ^*(0,u)}.
\end{equation}\\
\end{lmm}

\noindent
Note that the sums on the left-hand sides of~\e_ref{regprp_e1} and~\e_ref{regprp_e2} 
are finite in each $u$-degree, since $\cZ^*$ has no degree-zero term.\\

\noindent
Suppose $(\eta,\bar\cZ^*)$ is the regularizing pair for $\cZ^*$ at $\hb\!=\!0$.
Let
\begin{equation}\label{bZexp_e}
\bar\cZ^*(\hb,u)=\sum_{q=0}^{\i}C_q(u)\hb^q
\end{equation}
be the Taylor series expansion for $\bar\cZ^*(\hb,u)$ at $\hb\!=\!0$.
If $a\!\in\!\Z$, then
\begin{equation}\label{Rdfn_e}
\fR_{\hb=0}\big\{\hb^a \cZ^*(\hb,u)\big\}
=\sum_{\underset{p,q\ge 0}{p-q=1+a}}\!\frac{\eta(u)^p}{p!}C_q(u)
+\begin{cases}
\frac{\eta(u)^{a+1}}{(a+1)!},&\hbox{if}~a\!\ge\!0;\\
0,&\hbox{otherwise}.\end{cases}
\end{equation}
The identities~\e_ref{regprp_e1} and~\e_ref{regprp_e2} follow from~\e_ref{Rdfn_e} 
by a fairly direct computation; see Appendix~\ref{comb_app}.\\

\noindent
It remains to show that if $\cZ^*$ satisfies~\e_ref{regprp_e1} for all $a\!\ge\!0$,
then $\cZ^*$ admits a regularization.
Since $\cZ^*\!\in\!\Q_{\al}(\hb)[[u]]$, we can expand $\cZ^*$ at $\hb\!=\!0$ as
$$\cZ^*(\hb,u) = \sum_{d=1}^{\i}\sum_{q=-N_d}^{\i}\!\!\ti{C}_{q;d}\hb^qu^d
\equiv\sum_{q\in\Z}\ti{C}_q(u)\hb^q, \qquad\hbox{where}\qquad
\ti{C}_q(u)=\fR_{\hb=0}\Big\{\hb^{-q-1}\cZ^*(\hb,u)\Big\}.$$
{\it Claim:} There exists $\eta\!\in\!\Q_{\al}[[u]]$ such that 
$$\fR_{\hb=0}\bar\cZ^*(\hb,u)=0, \qquad\hbox{where}\qquad
1+\cZ^*(\hb,u)=e^{\eta(u)/\hb}\big(1+\bar\cZ^*(\hb,u)\big).$$
Since $\bar\cZ^*\!\in\!\Q_{\al}(\hb)[[u]]$, 
we can expand $\bar\cZ^*$ at $\hb\!=\!0$ as
$$\bar\cZ^*(\hb,u) = 
\sum_{d=1}^{\i}\sum_{q=-\bar{N}_d}^{\i} \!\!\!C_{q;d}\hb^qu^d
\equiv\sum_{q\in\Z}C_q(u)\hb^q, \qquad\hbox{where}\qquad
C_q(u)=\fR_{\hb=0}\Big\{\hb^{-q-1}\bar\cZ^*(\hb,u)\Big\}.$$
By assumption on $\eta$, $C_{-1}(u)\!=\!0$.
Let 
$$\bar\cY^*(\hb,u)=\sum_{q=0}^{\i}C_q(u)\hb^q\in\Q_{\al}(\hb)[[u]]
\qquad\hbox{and}\qquad
1+\cY^*(\hb,u)=e^{\eta(u)/\hb}\big(1+\bar\cY^*(\hb,u)\big).$$
Since $C_{-1}(u)\!=\!0$, 
\begin{equation}\label{regprppf_e5}
\fR_{\hb=0}\Big\{\hb^{-a}\cY^*(\hb,u)\Big\}
=\fR_{\hb=0}\Big\{\hb^{-a}\cZ^*(\hb,u)\Big\} \qquad\forall\,a\!\ge\!0.
\end{equation}
Since $\bar\cY^*$ is holomorphic at~$\hb\!=\!0$, $\cY^*$ satisfies~\e_ref{regprp_e1},
with $\cZ^*$ replaced by~$\cY^*$.
Thus, for all $a\!\ge\!0$ 
\begin{equation*}\begin{split}
a!\,\fR_{\hb=0}\Big\{\hb^{a+1}\cY^*(\hb,u)\Big\}
&=\sum_{m=2}^{\i} \frac{1}{m(m\!-\!1)}
\sum_{\underset{a_l\ge0}{\sum\limits_{l=1}\limits^{l=m}\!\!a_l=m-2-a}}
\!\!\!\!\Bigg(\prod_{l=1}^{l=m}\frac{(-1)^{a_l}}{a_l!}
\fR_{\hb=0}\Big\{\hb^{-a_l}\cY^*(\hb,u)\Big\}\Bigg)\\
&=\sum_{m=2}^{\i} \frac{1}{m(m\!-\!1)}
\sum_{\underset{a_l\ge0}{\sum\limits_{l=1}\limits^{l=m}\!\!a_l=m-2-a}}
\!\!\!\!\Bigg(\prod_{l=1}^{l=m}\frac{(-1)^{a_l}}{a_l!}
\fR_{\hb=0}\Big\{\hb^{-a_l}\cZ^*(\hb,u)\Big\}\Bigg)\\
&=a!\,\fR_{\hb=0}\Big\{\hb^{a+1}\cZ^*(\hb,u)\Big\}=a!\ti{C}_{-a-2}.
\end{split}\end{equation*}
Along with~\e_ref{regprppf_e5}, this implies that $\cZ^*\!=\!\cY^*$.
Thus, $\bar\cZ^*\!=\!\bar\cY^*$ is holomorphic at $\hb\!=\!0$.\\

\noindent
{\it Proof of the Claim:} The required property of $\eta$ is equivalent to
$$\eta=\sum_{q=0}^{\i}\frac{(-\eta)^q}{q!}\ti{C}_{q-1}.$$
Since $C_q\!\in\!\Q_{\al}[[u]]$ has no degree-zero term,
this equation has a unique power-series solution $\eta\!\in\!\Q_{\al}[[u]]$
with no degree-zero term.\footnote{Take $\eta_0\!=\!\ti{C}_{-1}$, 
$\eta_p\!=\!\sum_{q=0}^p\frac{(-\eta_{p-1})^q}{q!}\ti{C}_{q-1}$ for $p\!\ge\!1$.
The sequence $\eta_0,\eta_1,\ldots\!\in\!\Q_{\al}[[u]]$ converges, since
it is constant in degree $d$ after the $d$-th term.}\\

\noindent
{\it Remark:} The identity~\e_ref{regprp_e2} is valid as long as the residue of 
$\bar\cZ^*$ at~$\hb$ vanishes, but $\bar\cZ^*$ is not necessarily holomorphic at $\hb\!=\!0$.
In such a case, $\bar\cZ^*(0,u)$ must be replaced by 
$\fR_{\hb=0}\big\{\hb^{-1}\bar\cZ^*(\hb,u)\big\}$ on the right-hand side of~\e_ref{regprp_e2}.
By the proof of the claim, this residue is completely determined by~$\cZ^*$.
The assumption that $\cZ^*$ is regularizable at $\hb\!=\!0$ allows us to
compute the sum in~\e_ref{regprp_e2} explicitly.

\subsection{Regularizability of GW Generating Functions}
\label{Zreg_subs}

\begin{lmm}
\label{Zreg_lmm}
The power series $\cZ_i^*\!=\!\cZ_i^*(\hb,u)\!\in\!\Q_{\al}(\hb)[[u]]$ 
defined in~\e_ref{Z1ptdfn_e} is regularizable at $\hb\!=\!0$.
\end{lmm}

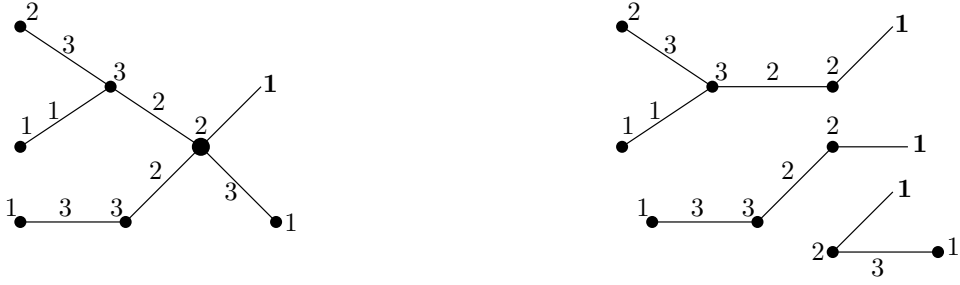
\begin{figure}
\begin{pspicture}(-0.5,-1.5)(10,2)
\psset{unit=.4cm}
\pscircle*(10,0){.3}\rput(10,.7){\smsize{$2$}}
\psline[linewidth=.04](10,0)(7,2)\pscircle*(7,2){.2}
\rput(8.6,1.5){\smsize{$2$}}\rput(7.3,2.5){\smsize{$3$}}
\psline[linewidth=.04](7,2)(4,4)\pscircle*(4,4){.2}
\rput(5.6,3.4){\smsize{$3$}}\rput(4.4,4.5){\smsize{$2$}}
\psline[linewidth=.04](7,2)(4,0)\pscircle*(4,0){.2}
\rput(5.1,1.2){\smsize{$1$}}\rput(4.2,.7){\smsize{$1$}}
\psline[linewidth=.04](10,0)(7.5,-2.5)\pscircle*(7.5,-2.5){.2}
\rput(8.5,-.8){\smsize{$2$}}\rput(7.2,-2){\smsize{$3$}}
\psline[linewidth=.04](7.5,-2.5)(4,-2.5)\pscircle*(4,-2.5){.2}
\rput(5.5,-2){\smsize{$3$}}\rput(3.7,-2){\smsize{$1$}}
\psline[linewidth=.04](10,0)(12.5,-2.5)\pscircle*(12.5,-2.5){.2}
\rput(11,-1.5){\smsize{$3$}}\rput(13,-2.5){\smsize{$1$}}
\psline[linewidth=.04](10,0)(12,2)\rput(12.3,2.1){\smsize{$\bf 1$}}
\pscircle*(31,2){.2}\rput(31,2.7){\smsize{$2$}}
\psline[linewidth=.04](31,2)(27,2)\pscircle*(27,2){.2}
\rput(29,2.5){\smsize{$2$}}\rput(27.3,2.5){\smsize{$3$}}
\psline[linewidth=.04](27,2)(24,4)\pscircle*(24,4){.2}
\rput(25.6,3.4){\smsize{$3$}}\rput(24.4,4.5){\smsize{$2$}}
\psline[linewidth=.04](27,2)(24,0)\pscircle*(24,0){.2}
\rput(25.1,1.2){\smsize{$1$}}\rput(24.2,.7){\smsize{$1$}}
\psline[linewidth=.04](31,2)(33,4)\rput(33.3,4.1){\smsize{$\bf 1$}}
\psline[linewidth=.04](31,0)(28.5,-2.5)\pscircle*(28.5,-2.5){.2}
\rput(29.5,-.8){\smsize{$2$}}\rput(28.2,-2){\smsize{$3$}}
\psline[linewidth=.04](28.5,-2.5)(25,-2.5)\pscircle*(25,-2.5){.2}
\rput(26.5,-2){\smsize{$3$}}\rput(24.7,-2){\smsize{$1$}}
\pscircle*(31,0){.2}\rput(31,.7){\smsize{$2$}}
\psline[linewidth=.04](31,0)(33.5,0)\rput(33.9,.1){\smsize{$\bf 1$}}
\psline[linewidth=.04](31,-3.5)(34.5,-3.5)\pscircle*(34.5,-3.5){.2}
\rput(32.5,-4){\smsize{$3$}}\rput(35,-3.3){\smsize{$1$}}
\pscircle*(31,-3.5){.2}\rput(30.5,-3.5){\smsize{$2$}}
\psline[linewidth=.04](31,-3.5)(33,-1.5)\rput(33.4,-1.4){\smsize{$\bf 1$}}
\end{pspicture}
\caption{A decorated tree contributing to $Q_i(\hb,u)$, with $i\!=\!2$, 
and its strands}
\label{Zgraph_fig}
\end{figure}

\noindent
{\it Proof:} We will verify that $\cZ^*\!\equiv\cZ_i^*$ satisfies \e_ref{regprp_e1} 
for all $a\!\ge\!0$.
By the string relation (see \cite[Section 26.3]{MirSym}),
\begin{equation}\label{Zreg_e1}
\cZ_i'^*(\hb,u)\equiv
\int_{\ov\M_{0,1}(\P^{n-1},d)}\frac{\E(\V_0')}{\hb\!-\!\psi_1}\ev_1^*\phi_i
=\hb\,\cZ_i^*(\hb,u).
\end{equation}
By the same argument as in the proof of \cite[Lemma~30.11]{MirSym},
\begin{gather}\label{Zreg_e3}
\cZ_i'^*(\hb,u)=Q_i(\hb,u)+ \sum_{d=1}^{\i}\sum_{j\neq i}
\frac{1}{\hb-\frac{\al_j-\al_i}{d}}\fR_{\hb=\frac{\al_j-\al_i}{d}}
\big\{\cZ_i'^*(\hb,u)\big\}\\
\hbox{for some}\qquad Q_i\in\Q_{\al}[\hb,\hb^{-1}]\big[\big[u\big]\big].\notag
\footnotemark
\end{gather}
\footnotetext{The statement of \cite[Lemma 30.11]{MirSym} is made for 
a renormalized version of the power series $\cZ_i^*(\hb,u)$ and is in fact 
sharper.}The middle 
term in~\e_ref{Zreg_e3} is the sum of contributions to $\cZ_i'^*(\hb,u)$
from the graphs $\Ga$ with one marked point such that the marked point is attached
to a vertex $v_0$ of valence at least~$3$.
Similarly to Subsections~\ref{localdata_subs1} and~\ref{localdata_subs2},
the vertex $v_0$ of $\Ga$ must then be labeled~$i$.
An example of such a graph is shown in Figure~\ref{Zgraph_fig}.\\

\noindent
Let $\Ga$ be a decorated tree with one marked point as in~\e_ref{decortgraphdfn_e} 
that contributes to $Q_i(\hb,u)$, i.e.
$$k=1, \quad \mu(v_0)=i, \quad \val(v_0)\ge3, 
\qquad\hbox{where}\quad v_0=\eta(1).$$
As in Subsections~\ref{localdata_subs1} and~\ref{localdata_subs2},
we break $\Ga$ into strands $\Ga_v$ indexed by the set
$$\Edg(v_0)=\big\{e\!\in\!\Edg\!: v_0\!\in\!e\big\}$$
of the edges leaving from $v_0$.
In this case, there are 
$$m\equiv \big|\Edg(v_0)\big|\ge 2$$
strands, each with exactly one marked point.\\

\noindent
The fixed locus $\cZ_{\Ga}$ corresponding to~$\Ga$, 
the restriction of $\E(\V_0')$ to~$\cZ_{\Ga}$, and 
the euler class of the normal bundle of~$\cZ_{\Ga}$ are given~by
\begin{equation}\label{Zreg_e5}\begin{split}
&\cZ_{\Ga}=\ov\cM_{0,|\val(v_0)|}\times\prod_{e\in\Edg(v_0)}\!\!\!\!\!\!\cZ_{\Ga_e},
\qquad \E(\V_0')=\prod_{e\in\Edg(v_0)}\!\!\!\!\!\!\pi_e^*\E(\V_0'),\\
&\frac{\E(T_{\mu(v_0)}\P^{n-1})}{\E(\N\cZ_{\Ga})}=
\prod_{e\in\Edg(v_0)}\!\!\bigg(\frac{\E(T_{\mu(v_0)}\P^{n-1})}
{\E(\N\cZ_{\Ga_e})\,(\hb_e\!-\!\pi_e^*\psi_1)}\bigg),
\end{split}\end{equation} 
where $\hb_e\!\equiv\!c_1(L_e')\!\in\!H^*(\ov\cM_{0,|\val(v_0)|})$ is
the first chern class of the universal tangent line bundle for the marked point
corresponding to the edge~$e$.
By \cite[Section~27.2]{MirSym}, if $e\!=\!\{v_0,v_e\}$, then
\begin{equation}\label{psiweigh_e}
\psi_1|_{\cZ_{\Ga_e}}=\frac{\al_{\mu(v_e)}\!-\!\al_{\mu(v_0)}}{\d(e)}
=\frac{\al_{\mu(v_e)}\!-\!\al_i}{\d(e)}.
\end{equation}
Thus, by \e_ref{Zreg_e5} and~\e_ref{tangrestr_e},
\begin{alignat}{1}
&\int_{\cZ_{\Ga}}\!\!\frac{\E(\V_0')\ev_1^*\phi_i}{(\hb\!-\!\psi_1)\E(\N\cZ_{\Ga})}
=\sum_{a=0}^{\i} \hb^{-(a+1)}(-1)^m 
 \!\!\!\!\! \sum_{\underset{e\in\Edg(v_0)}{a_e\ge0}} \!\!\!\!\! 
 \Bigg\{ \int_{\ov\cM_{0,|\val(v_0)|}} \psi_1^a 
\!\!\!\!\!\!\prod_{~~e\in\Edg(v_0)}\!\!\!\!\!\!\!\!\! \hb_e^{a_e} \notag\\
\label{Zreg_e6}
&\hspace{2.5in}
\times\!\!\! \prod_{e\in\Edg(v_0)}\!\!\! \Bigg(
\bigg(\frac{\al_{\mu(v_e)}\!-\!\al_i}{\d(e)}\bigg)^{-(a_e+1)}
\!\!\!\int_{\cZ_{\Ga_e}}\!\!\!\frac{\E(\V_0')\ev_1^*\phi_i}{\E(\N\cZ_{\Ga_e})}\Bigg)\Bigg\}\\
&\quad
=\sum_{a+\!\!\!\!\!\!\!\!
\underset{a_e\ge0}{\sum\limits_{e\in\Edg(v_0)}\!\!\!\!\!\!\!\!\!a_e=m-2}}\!\!
\binom{m\!-\!2}{a,(a_e)_{e\in\Edg(v_0)}}(-1)^a\hb^{-(a+1)}
\!\!\!\!\!\!\prod_{e\in\Edg(v_0)}\!\!\! \Bigg(
\bigg(\frac{\al_{\mu(v_e)}\!-\!\al_i}{\d(e)}\bigg)^{-(a_e+1)}
\!\!\!\int_{\cZ_{\Ga_e}}\!\!\!\frac{\E(\V_0')\ev_1^*\phi_i}{\E(\N\cZ_{\Ga_e})}
\Bigg).\notag
\end{alignat}
The first equality holds after dividing the expressions on the right-hand side
by the order of the appropriate groups of symmetries; see~\cite[Section 27.3]{MirSym}.
This group is taken into account in the next paragraph.\\

\noindent
We now sum up \e_ref{Zreg_e6} over all possibilities for $\Ga$.
By~\e_ref{Zreg_e3} and its proof, for every $j\!\in\![n]\!-\!i$ and $d\!\in\!\Z^+$,
\begin{equation}\label{Zreg_e7}\begin{split}
\sum_{\underset{\mu(v_e)=j,\d(e)=d}{\Ga_e}}\!\!\!\!\!\!\!\!
\bigg(\frac{\al_j\!-\!\al_i}{d}\bigg)^{-(a_e+1)}
\!\!\!\int_{\cZ_{\Ga_e}}\!\!\!\frac{\E(\V_0')\ev_1^*\phi_i}{\E(\N\cZ_{\Ga_e})}
&=\fR_{\hb=\frac{\al_j-\al_i}{d}}\Big\{\hb^{-(a_e+1)}\cZ_i'^*(\hb,u)\Big\}\\
&=\fR_{\hb=\frac{\al_j-\al_i}{d}}\Big\{\hb^{-a_e}\cZ_i^*(\hb,u)\Big\}.
\footnotemark
\end{split}\end{equation}
\footnotetext{The proof is the same as the proof of~\e_ref{Z1pt_lmm_e1} 
in \cite[Chapter~30]{MirSym}.
If a graph $\Ga$ contributes a factor of $\hb\!-\!(\al_j\!-\!\al_i)/d$
to the denominator of $\cZ_i'^*$ (or of $\hb^{-(a_e+1)}\cZ_i'^*$), then 
the marked point $1$ is attached to a vertex of $\Ga$ of valence~$2$.
Furthermore, if $e$ is the unique element of $\Edg(v_0)$, then
$\mu(v_e)\!=\!j$ and $\d(e)\!=\!d_v$; 
see \cite[Section 30.1]{MirSym} for more details.
The second equality in~\e_ref{Zreg_e7} is immediate from~\e_ref{Zreg_e1}.}
Since $\hb^{-a_e}\cZ_i^*(\hb,u)$ has no residue at $\hb\!=\!\i$ by Lemma~\ref{Z1pt_lmm},
\begin{equation}\label{Zreg_e8}\begin{split}
\sum_{\Ga_e}
\bigg(\frac{\al_j\!-\!\al_i}{d}\bigg)^{-(a_e+1)}
\!\!\!\int_{\cZ_{\Ga_e}}\!\!\!\frac{\E(\V_0')\ev_1^*\phi_i}{\E(\N\cZ_{\Ga_e})}
&=\sum_{d=1}^{\i}\sum_{j\neq i}
\fR_{\hb=\frac{\al_j-\al_i}{d}}\Big\{\hb^{-a_e}\cZ_i^*(\hb,u)\Big\}\\
&=-\fR_{\hb=0}\Big\{\hb^{-a_e}\cZ_i^*(\hb,u)\Big\}
\end{split}\end{equation}
by~\e_ref{Zreg_e7}, the Residue Theorem on~$S^2$, and Lemma~\ref{Z1pt_lmm}.
By~\e_ref{Zreg_e6} and~\e_ref{Zreg_e8},
\begin{equation}\label{Zreg_e9}\begin{split}
&\sum_{\underset{|\Edg(v_0)|=m}{\Ga}}\!\!\!\!\!\!
\int_{\cZ_{\Ga}}\!\!\frac{\E(\V_0')\ev_1^*\phi_i}{\E(\N\cZ_{\Ga})}\\
&\qquad
=\sum_{a=0}^{\i}\hb^{-(a+1)}\Bigg\{
(-1)^{m-a} \!\!\!\!\!\!\!\!\!\!
\sum_{\underset{a_e\ge0}{\sum\limits_{e\in\Edg(v_0)}\!\!\!\!\!\!\!\!\!a_e=m-2-a}}
\!\! \binom{m\!-\!2}{a,(a_e)_{e\in\Edg(v_0)}}
\prod_{e\in\Edg(v_0)}\!\!\!\!\!\! \fR_{\hb=0}\Big\{\hb^{-a_e}\cZ_i^*(\hb,u)\Big\}\Bigg\}.
\end{split}\end{equation}
Taking into account the group of symmetries, i.e.~$S_m$ in the case of \e_ref{Zreg_e9},
and summing up over all possible values of $|\Edg(v_0)|$, i.e.~$m\!\ge\!2$, we obtain
\begin{equation}\label{Zreg_e10}
Q_i(\hb,u)=\sum_{a=0}^{\i}\frac{\hb^{-(a+1)}}{a!}
\sum_{m=2}^{\i} \frac{1}{m(m\!-\!1)}
\sum_{\underset{a_l\ge0}{\sum\limits_{l=1}\limits^{l=m}\!\!a_l=m-2-a}}
\!\!\!\!\Bigg(\prod_{l=1}^{l=m}\frac{(-1)^{a_l}}{a_l!}
\fR_{\hb=0}\Big\{\hb^{-a_l}\cZ_i^*(\hb,u)\Big\}\Bigg).
\end{equation}\\

\noindent
On the other hand, by the Residue Theorem on $S^2$, Lemma~\ref{Z1pt_lmm},
and~\e_ref{Zreg_e1},
\begin{equation*}\begin{split}
\sum_{d=1}^{\i}\sum_{j\neq i}
\frac{1}{\hb-\frac{\al_j-\al_i}{d}}\fR_{\hb=\frac{\al_j-\al_i}{d}}
\big\{\cZ_i'^*(\hb,u)\big\}
&=\sum_{d=1}^{\i}\sum_{j\neq i}\fR_{z=\frac{\al_j-\al_i}{d}}
\Big\{\frac{1}{\hb\!-\!z}\cZ_i'^*(z,u)\Big\}\\
&=-\fR_{z=\hb,0}\Big\{\frac{1}{\hb\!-\!z}\cZ_i'^*(z,u)\Big\}\\
&=\cZ_i'^*(\hb,u)-\sum_{a=0}^{\i}\hb^{-(a+1)}
\fR_{z=0}\Big\{z^a\cZ_i'^*(z,u)\Big\}\\
&=\cZ_i'^*(\hb,u)-\sum_{a=0}^{\i}\hb^{-(a+1)}
\fR_{\hb=0}\Big\{\hb^{a+1}\cZ_i^*(\hb,u)\Big\}.
\end{split}\end{equation*}
Comparing with \e_ref{Zreg_e3}, we conclude that  
\begin{equation}\label{Zreg_e12}
Q_i(\hb,u)=\sum_{a=0}^{\i}\hb^{-(a+1)}
\fR_{\hb=0}\Big\{\hb^{a+1}\cZ_i^*(\hb,u)\Big\}.
\end{equation}
By~\e_ref{Zreg_e10} and~\e_ref{Zreg_e12}, $\cZ_i^*$ satisfies \e_ref{regprp_e1}
for all $a\!\ge\!1$.

\subsection{Proofs of Propositions~\ref{contr_prp1} and~\ref{contr_prp2}}
\label{localprppf_subs}

\noindent
In this subsection we prove Propositions~\ref{contr_prp1} and~\ref{contr_prp2}.
An argument nearly identical to part of the proof of Lemma~\ref{Zreg_subs}
leads to a long expression like~\e_ref{Zreg_e6}.
In the proof of the first proposition, the Residue Theorem on~$S^2$ reduces it
to the form~\e_ref{Zreg_e9}. We then use the second statement of Lemma~\ref{regular_lmm}
to deal with the infinite summation.
The situation in Proposition~\ref{contr_prp2} is a bit different,
as the possible strands $\Ga_e$ of $\Ga$ are not mutually independent
due the requirements~(i) and~(ii) on~$\Ga$ in Subsection~\ref{localdata_subs2}.
In this case, we will use Lemma~\ref{prodres_lmm} to reduce a product of
residues to the residue of a product; the resulting products are readily summable.
The Residue Theorem on~$S^2$ is used at the last~step.
As the proof of Lemma~\ref{prodres_lmm} is completely straightforward,
we relegate it to Appendix~\ref{comb_app}.\\

\noindent
If $f\!=\!f(\la)$ is holomorphic at $\la\!=\!0$ and $m\!\ge\!0$, let
$$\cD_{\la}^mf=\frac{1}{m!}\frac{d^m}{d\la^m}f(\la)\Big|_{\la=0}.$$

\begin{lmm}
\label{prodres_lmm}
If $(f_e\!=\!f_e(\la))_{e\in E}$ is a finite collection of functions 
with at most a simple pole at $\la\!=\!0$, then
\begin{equation*}\begin{split}
\fR_{\la=0}\bigg\{\prod_{e\in E}\! f_e(\la)\bigg\}
=\sum_{E_+\subset E}\!\!
\Bigg\{\prod_{e\in E_+}\!\!\fR_{\la=0}\big\{f_e(\la)\big\}
\cdot \cD_{\la}^{|E_+|-1} \bigg(\prod_{e\not\in E_+}\!\!\!\Big(f_e(\la)
\!-\!\la^{-1}\fR_{\la=0}\big\{f_e(\la)\big\}\Big)\bigg) \Bigg\}.
\end{split}\end{equation*}\\
\end{lmm}

\noindent
In the case of (i) of Proposition~\ref{contr_prp1}, 
equations~\e_ref{Abndl_e} and~\e_ref{Anormalbndl_e}
describe the splitting of the integrand corresponding to each fixed locus $\cZ_{\Ga}$
in the sense of~\e_ref{ABothm_e} between the strands of~$\Ga$.
Summing over all possible strands as in Subsection~\ref{Zreg_subs},
we obtain
\begin{equation}\label{contrprp1_e1}\begin{split}
\cA_i(u)=\sum_{m=0}^{\i}
\sum_{\underset{a_{\pm},a_l\ge0}{a_-+a_++\sum\limits_{l=1}\limits^{l=m}\!\!a_l=m}}
\!\!\!\Bigg(\frac{(-1)^{a_-+a_+}}{a_-!a_+!}\fR_{\hb_1=0}\fR_{\hb_2=0}\Big\{
\hb_1^{-a_-}\hb_2^{-a_+}\ti\cZ_{ii}^*(\hb_1,\hb_2,u)\Big\}&\\
\times\prod_{l=1}^{l=m}\frac{(-1)^{a_l}}{a_l!}
\fR_{\hb=0}\Big\{\hb^{-a_l}\cZ_i^*(\hb,u)\Big\}&\Bigg).
\end{split}\end{equation}
This is the analogue of~\e_ref{Zreg_e10}.
In this case, we use Lemma~\ref{Z2pt_lmm}, in addition to Lemma~\ref{Z1pt_lmm} 
and the Residue Theorem on~$S^2$, to obtain~\e_ref{contrprp1_e1} from 
the analogue of~\e_ref{Zreg_e6} for~$\cA_i(u)$.
The sum is taken over every possible number $m$ of strands with 
one marked point.
The ends of the distinguished strand $\Ga_{\pm}$ are ordered, 
accounting for the factor of $1/2$ in~\e_ref{Z2ptdfn_e}.
By Lemma~\ref{Zreg_lmm}, $\cZ_i^*$ is regularizable at $\hb\!=\!0$.
Therefore, (ii) of Lemma~\ref{regular_lmm} reduces the right-hand side
of~\e_ref{contrprp1_e1} to the right-hand side of~\e_ref{contr_prp1e1}.\\

\noindent
In the case of (ii) of Proposition~\ref{contr_prp1}, 
the analogue of~\e_ref{contrprp1_e1} is easily seen to be
\begin{equation*}\begin{split}
\ti{A}_{ij}(u)=\frac{1}{\prod\limits_{k\neq j}(\al_j\!-\!\al_k)}\sum_{m=0}^{\i}
\sum_{\underset{a_0,a_{\pm},a_l\ge0}{a_-+a_++a_0+\sum\limits_{l=1}\limits^{l=m}\!\!a_l=m}}
\!\!\!\Bigg(\frac{(-1)^{a_-+a_+}}{a_-!a_+!}\fR_{\hb_1=0}\fR_{\hb_2=0}\Big\{
\hb_1^{-a_-}\hb_2^{-a_+}\ti\cZ_{jj}^*(\hb_1,\hb_2,u)\Big\}&\\
\times\frac{(-1)^{a_0}}{a_0!}\fR_{\hb=0}
\Big\{\hb^{-a_0}\cZ_{ji}^*(\hb,u)\Big\}
\cdot\prod_{l=1}^{l=m}\frac{(-1)^{a_l}}{a_l!}
\fR_{\hb=0}\Big\{\hb^{-a_l}\cZ_j^*(\hb,u)\Big\}&\Bigg).
\end{split}\end{equation*}
In this case, Lemma~\ref{Z1pt2_lmm} is used in addition to 
Lemmas~\ref{Z1pt_lmm} and~\ref{Z2pt_lmm} and the Residue Theorem on~$S^2$.
Lemmas~\ref{regular_lmm} and~\ref{Zreg_lmm} reduce the right-hand side
of the above expression to the right-hand side of~\e_ref{contr_prp1e2}.\\

\noindent
In the case of (i) of Proposition~\ref{contr_prp2}, 
\e_ref{Bbndl_e} and~\e_ref{Bnormalbndl_e} describe the splitting of 
the integrand for the fixed locus $\wt\cZ_{\Ga}$ between the strands of~$\Ga$.
Let
$$\Edg_-=\Edg(v_0)-\Edg_+, \qquad m_+=\big|\Edg_+|, \qquad
\Psi(\hb,\ti\psi)= -\frac{\prod_{k\neq i}(\al_i\!-\!\al_k\!+\!\hb)}
{(n\al_i\!+\!\hb)(\hb\!+\!\ti\psi)}\,.$$
The analogue of~\e_ref{Zreg_e6} is then
\begin{equation}\label{contrprp2_e1}\begin{split}
&\int_{\wt\cZ_{\Ga}}\!\!\frac{\E(\V_1')\ev_1^*\phi_i}{\E(\N\wt\cZ_{\Ga})}
=\int_{\wt\cM_{1,|\val(v_0)|}\times\P^{m_+-1}}\Bigg\{\Psi(\hb,\ti\psi)
\prod_{e\in\Edg_+}\!\!\! \bigg(\int_{\cZ_{\Ga_e}}\!\!
\frac{\E(\V_0')\ev_1^*\phi_i}{\E\big(\N\cZ_{\Ga_e}\big)}\bigg)\\
&\hspace{2.7in}
\times\prod_{e\in\Edg_-}\!\!\! \bigg(\int_{\cZ_{\Ga_e}}\!\!
\frac{\E(\V_0')\ev_1^*\phi_i}{(\hb\!-\!\psi_e)\E\big(\N\cZ_{\Ga_e}\big)}
\bigg)\Bigg\}_{\hb=\psi_{\Ga}+\la}\,.
\end{split}\end{equation}\\

\noindent
We now sum \e_ref{contrprp2_e1} over all possibilities for $\Ga_e$
with $e\!\in\!\Edg_-$.
In contrast to the three cases encountered above, $\Ga_e$ can be any graph with 
one marked point such that $\mu(1)\!=\!i$, as long as the edge leaving 
the vertex $\mu(1)$ is not labeled~$\d(\Ga)$ or its other end is not 
labeled~$\mu(\Ga)$.
This restriction is due to~\e_ref{notthickedges_e}.
By~\e_ref{ABothm_e} applied to the function $\cZ_i'^*(\hb,u)$ defined in~\e_ref{Zreg_e1},
the sum over such graphs $\Ga_e$ is 
\begin{equation}\label{contrprp2_e2}\begin{split}
&\sum_{\underset{(\mu(v_e),\d(e))\neq(\mu(\Ga),\d(\Ga))}{\Ga_e}} \!\!\!\!\!
\int_{\cZ_{\Ga_e}}
\frac{\E(\V_0')\ev_1^*\phi_i}{(\hb\!-\!\psi_e)\E\big(\N\cZ_{\Ga_e}\big)}\\
&\hspace{1.5in}
=\cZ_i'^*(\hb,u) -\sum_{\underset{(\mu(v_e),\d(e))
=(\mu(\Ga),\d(\Ga))}{\Ga_e}}\!\!\!\!
\int_{\cZ_{\Ga_e}}\frac{\E(\V_0')
\ev_1^*\phi_i}{(\hb\!-\!\psi_e)\E\big(\N\cZ_{\Ga_e}\big)}\\
&\hspace{1.5in}
=\cZ_i'^*(\hb,u)-
\frac{1}{\hb\!-\!\psi_{\Ga}}\,\fR_{z=\psi_{\Ga}}\big\{\cZ_i'^*(z,u)\big\},
\end{split}\end{equation}
since $\psi_{\Ga}=(\al_{\mu(\Ga)}\!-\!\al_i)/\d(\Ga)$.
The last equality uses~\e_ref{Z1pt_lmm_e1} and~\e_ref{Zreg_e1}.
On the other hand, summing over all possibilities for $\Ga_e$ with $e\!\in\!\Edg_+$,
without changing $(\mu(\Ga),\d(\Ga))$, we obtain
\begin{equation}\label{contrprp2_e3}\begin{split}
&\sum_{\underset{(\mu(v_e),\d(e))=(\mu(\Ga),\d(\Ga))}{\Ga_e}} \!\!\!\!\!
\int_{\cZ_{\Ga_e}}\frac{\E(\V_0')\ev_1^*\phi_i}{\E\big(\N\cZ_{\Ga_e}\big)}
=\fR_{z=\psi_{\Ga}}\big\{\cZ_i'^*(z,u)\big\},
\end{split}\end{equation}
also by~\e_ref{Z1pt_lmm_e1} and~\e_ref{Zreg_e1}.\\

\noindent
By~\e_ref{contrprp2_e2} and~\e_ref{contrprp2_e3}, the sum of the terms in~\e_ref{contrprp2_e1}
with $\Edg(v_0)$, $\Edg_+$, and $(\mu(\Ga),\d(\Ga))$ fixed is
\begin{equation}\label{contrprp2_e4}\begin{split}
&\int_{\wt\cM_{1,|\val(v_0)|}\times\P^{m_+-1}}\Bigg\{
\prod_{e\in\Edg_+}\!\!\! \Big(\fR_{z=\psi_{\Ga}}\big\{\cZ_i'^*(z,u)\big\}\Big)
\Psi(\hb,\ti\psi)\\
&\hspace{1.8in}
\times\prod_{e\in\Edg_-}\!\!\!\Big(\cZ_i'^*(\hb,u)-\frac{1}{\hb\!-\!\psi_{\Ga}}\,
\fR_{z=\psi_{\Ga}}\big\{\cZ_i'^*(z,u)\big\}\Big) \Bigg\}_{\hb=\psi_{\Ga}+\la}\\
&=\int_{\wt\cM_{1,|\val(v_0)|}}\Bigg\{
\prod_{e\in\Edg_+}\!\!\! \Big(\fR_{z=\psi_{\Ga}}\big\{\cZ_i'^*(z,u)\big\}\Big)\\
&\hspace{0.9in} \times
\cD_{\la}^{m_+-1}\bigg(\Psi(\psi_{\Ga}\!+\!\la,\ti\psi)
\prod_{e\in\Edg_-}\!\!\!\!\!\! \Big(\cZ_i'^*(\psi_{\Ga}\!+\!\la,u)-
\la^{-1} \fR_{z=\psi_{\Ga}}\big\{\cZ_i'^*(z,u)\big\}\Big)\bigg)\Bigg\}.
\end{split}\end{equation}
By the last expression in~\e_ref{contrprp2_e4} and Lemmas~\ref{prodres_lmm}
and~\ref{Z1pt_lmm}, the sum of the terms in~\e_ref{contrprp2_e1} with 
only $m\!\equiv\!|\Edg(v_0)|$ and $(\mu(\Ga),\d(\Ga))$ fixed is
\begin{equation}\label{contrprp2_e5}\begin{split}
&\int_{\wt\cM_{1,|\val(v_0)|}}\fR_{z=\psi_{\Ga}}\bigg(
\Psi(z,\ti\psi)\prod_{e\in\Edg(v_0)}\!\!\!\!\!\cZ_i'^*(z,u)\bigg)\\
&\hspace{1in}
=\frac{(-1)^m m!}{24}\fR_{z=\frac{\al_{\mu(\Ga)}-\al_i}{\d(\Ga)}}
\bigg\{z^{-(m+2)}\frac{\prod_{k\neq i}(\al_i\!-\!\al_k\!+\!z)}{(n\al_i\!+\!z)}
\cZ_i'^*(z,u)^m\bigg\}\\
&\hspace{1in}
=\frac{(-1)^m m!}{24}\fR_{z=\frac{\al_{\mu(\Ga)}-\al_i}{\d(\Ga)}}
\bigg\{\frac{\prod_{k\neq i}(\al_i\!-\!\al_k\!+\!z)}{(n\al_i\!+\!z)z^2}
\cZ_i^*(z,u)^m\bigg\}.
\end{split}\end{equation}
The first equality above follows from~\e_ref{Lnum_e}, while the second 
one follows  from~\e_ref{Zreg_e1}.\\

\noindent
Finally, taking into account \e_ref{eulerclass_e} and the group of symmetries, i.e.~$S_m$,
and summing \e_ref{contrprp2_e5} over all possible numbers of one-pointed 
strands, i.e.~$m\!\ge\!1$, and all possible values of~$(\mu(\Ga),\d(\Ga))$,
we obtain
\begin{equation}\label{contrprp2_e6}\begin{split}
\cB_i(u)&=\frac{n\al_i}{24}\sum_{j\neq i}\sum_{d=1}^{\i}\sum_{m=1}^{\i}
(-1)^m\fR_{\hb=\frac{\al_j-\al_i}{d}}
\bigg\{\frac{\prod_{k\neq i}(\al_i\!-\!\al_k\!+\!\hb)}{(n\al_i\!+\!\hb)\hb^2}
\cZ_i^*(\hb,u)^m\bigg\}\\
&=-\frac{n\al_i}{24}\sum_{j\neq i}\sum_{d=1}^{\i}
\fR_{\hb=\frac{\al_j-\al_i}{d}}
\bigg\{\frac{\prod_{k=1}^{k=n}(\al_i\!-\!\al_k\!+\!\hb)}{(n\al_i\!+\!\hb)\hb^3}
\frac{\cZ_i^*(\hb,u)}{1\!+\!\cZ_i^*(\hb,u)}\bigg\}.
\end{split}\end{equation}
The first claim of Proposition~\ref{contr_prp2} now follows from the Residue Theorem on~$S^2$
and Lemma~\ref{Z1pt_lmm}.\\

\noindent
The proof of (ii) of Proposition~\ref{contr_prp2} is nearly identical.
In the starting equation~\e_ref{contrprp2_e1}, $i$ is replaced by~$j$ 
and the entire expression is divided by $\prod_{k\neq j}(\al_j\!-\!\al_k)$.
One of the strands now has two marked points 
and  there are $m\!\ge\!0$ other strands.
The two-pointed strand can appear as an element of $\Edg_+$ as well as of~$\Edg_-$
and contributes to $\cZ_{ji}^*(\hb,u)$ instead of~$\cZ_j^*(\hb,u)$.
The number $|\val(v_0)|$ is still $m\!+\!1$.
Therefore, \e_ref{contrprp2_e6} becomes
\begin{equation*}\begin{split}
\ti\cB_{ij}(u)&=\frac{1}{\prod\limits_{k\neq j}(\al_j\!-\!\al_k)}
\frac{n\al_j}{24}\sum_{l\neq j}\sum_{d=1}^{\i}\sum_{m=0}^{\i}
(-1)^m\fR_{\hb=\frac{\al_l-\al_j}{d}}
\bigg\{\frac{\prod_{k\neq j}(\al_j\!-\!\al_k\!+\!\hb)}{(n\al_j\!+\!\hb)\hb}
\cZ_j^*(\hb,u)^m\cZ_{ji}^*(\hb,u)\bigg\}\\
&=\frac{1}{\prod\limits_{k\neq j}(\al_j\!-\!\al_k)}
\frac{n\al_j}{24}\sum_{l\neq j}\sum_{d=1}^{\i} \fR_{\hb=\frac{\al_l-\al_j}{d}}
\bigg\{\frac{\prod_{k=1}^{k=n}(\al_j\!-\!\al_k\!+\!\hb)}{(n\al_j\!+\!\hb)\hb^2}
\frac{\cZ_{ji}^*(\hb,u)}{1\!+\!\cZ_j^*(\hb,u)}\bigg\}.
\end{split}\end{equation*}
The second claim of Proposition~\ref{contr_prp2} now follows from the Residue Theorem,
along with Lemmas~\ref{Z1pt_lmm} and~\ref{Z1pt2_lmm}.\\

\noindent
{\it Remark:} In (i) and (ii) of Proposition~\ref{contr_prp1},
$|\val(v_0)|\!=\!m\!+\!3$, where $m$ is the number of one-pointed strands.
Of the extra $3$, $2$ comes from the distinguished  strand~$\Ga_{\pm}$.
The remaining $1$ comes from the marked point $1$ that lies on the contracted 
component $\cC_{f,v_0}$ corresponding to the vertex~$v_0$ in~(i)
and from the second distinguished strand $\Ga_0$ in~(ii).
In (i) and (ii) of Proposition~\ref{contr_prp2},
$|\val(v_0)|\!=\!m\!+\!1$, where $m$ is again the number of one-pointed strands.
The extra $1$ comes from the marked point $1$ that lies on $\cC_{f,v_0}$ 
in~(i) and from the two-pointed strand in~(ii).

\section{Algebraic Computations}
\label{alg_sec}

\noindent
In this section we use Lemmas~\ref{algebr_lmm1}-\ref{Z2pt_lmm} to deduce the main theorem 
of this paper, Theorem~\ref{main_thm} below, from Propositions~\ref{contr_prp1},
\ref{contr_prp2}, \ref{hypergeom_prp}, and~\ref{hypergeom_prp2}.
Theorem~\ref{main_thm} expresses the contribution from each of 
the two types of $\T$-fixed loci in $\wt\M_{1,1}^0(\P^{n-1},d)$ 
to the generating function $\F_0(u)$ for the reduced genus~$1$ GW-invariants of 
a degree~$n$ hypersurface in~$\P^{n-1}$ in terms of hypergeometric series.
Along with~\e_ref{derivgen_e2}, it immediately implies Theorem~\ref{main_thm0}.
Propositions~\ref{hypergeom_prp} and~\ref{hypergeom_prp2} describe the structure
of the function $R\!=\!R(w,t)$ defined in~\e_ref{Rdfn_e} 
at $w\!=\!0$ and $w\!=\!\infty$, respectively; they are proved in~\cite{ZaZ}.
Lemma~\ref{algebr_lmm1} serves a tool for extracting the non-equivariant part
of an equivariant cohomology class, while Lemmas~\ref{Z1pt_lmm}-\ref{Z2pt_lmm}
recall the relevant information about genus~$0$ generating functions.
Let
\begin{equation}\label{mudfn_e0}
\mu(e^t)\equiv\fR_{\hb=0}\big\{\ln R(\hb^{-1},t)\big\}-t.
\end{equation}

\begin{mythm}
\label{main_thm}
The generating function $\F_0(u)$ defined in \e_ref{pushclass_e} is given~by
$$\F_0(e^T)=\frac{d}{dT}\big(\ti{A}(e^t)+\ti{B}(e^t)\big),$$
where $T$ and $t$ are related by the mirror map \e_ref{mirmap_e} and
\begin{alignat}{1}
\label{main_thm_e1}
\ti{A}(e^t)&=\frac{1}{2}\Bigg(\frac{(n\!-\!2)(n\!+\!1)}{24}\mu(e^t)
-\frac{(n\!-\!2)(3n\!-\!5)}{24}\ln\big(1\!-\!n^ne^t\big)
-\sum_{p=0}^{n-3}\binom{n\!-\!1\!-\!p}{2}\ln I_{p,p}(t)\Bigg)\\
&=\frac{(n\!-\!2)(n\!+\!1)}{48}\mu(e^t)-
\begin{cases}
\frac{n+1}{48}\ln\big(1\!-\!n^ne^t\big)+
\sum_{p=0}^{(n-3)/2}\frac{(n-1-2p)^2}{8}\ln I_{p,p}(t),
&\hbox{if}~2\!\not|n;\\
\frac{n-2}{48}\ln\big(1\!-\!n^ne^t\big)+
\sum_{p=0}^{(n-4)/2}\frac{(n-2p)(n-2-2p)}{8}\ln I_{p,p}(t),
&\hbox{if}~2|n;
\end{cases}\notag\\
\label{main_thm_e2}
\ti{B}(e^t)&=\bigg(\frac{(n\!-\!2)(n\!+\!1)}{48}+\frac{1-(1\!-\!n)^n}{24n^2}\bigg)(T\!-\!t)
-\frac{(n\!-\!2)(n\!+\!1)}{48}\mu(e^t)+\frac{1}{24}\ln\big(1-n^ne^t\big)\\
&\hspace{1.2in}
+\frac{n^2\!-\!1+(1\!-\!n)^n}{24n}\ln I_{0,0}(t)
+\frac{n}{24}\sum_{p=2}^{n-2}\bigg(\cD_w^{n-2-p}\frac{(1\!+\!w)^n}{(1\!+\!nw)}\bigg)
\big(\cD_w^p\ln\bar{R}(w,t)\big)\notag
\end{alignat}
are $T$-integrals of the contributions of the effective fixed loci of
Subsection~\ref{localdata_subs1} and of the boundary fixed loci of 
Subsection~\ref{localdata_subs2}, respectively.
\end{mythm}

\noindent
The next two propositions, which are proved in~\cite{ZaZ}, are used in the proof of Theorem~\ref{main_thm} in Subsection~\ref{mainthmpf_subs}.

\begin{prp}
\label{hypergeom_prp}
(i) There exist $\ti{I}_{p,r}\!\in\!\Q[[e^t]]$ for 
$p,r\in\!\bar\Z^+$ with $r\!\ge\!p$ such~that 
$$I_{p,q}(t)=\sum_{r=p}^{r=q}\frac{t^{q-r}}{(q\!-\!r)!}\ti{I}_{p,r}(t)
\qquad\forall~p,q\!\in\!\bar\Z^+ ~\hbox{with}~q\!\ge\!p$$
and the constant term of $\ti{I}_{p,r}$ is $1$ for $r\!=\!p$
and is $0$ for $r\!>\!p$.\\
(ii) The power series $I_{p,p}$ in $e^t$ with $p\!=\!0,1,\ldots,n\!-\!1$ satisfy
\begin{alignat}{1}
\label{hypergeom_prp_e1}
& I_{0,0}(t)I_{1,1}(t)\ldots I_{n-1,n-1}(t)\big(1\!-\!n^ne^t\big)=1,\\
\label{hypergeom_prp_e2}
& I_{0,0}(t)^{n-1}I_{1,1}(t)^{n-2}\ldots I_{n-1,n-1}(t)^0
\big(1\!-\!n^ne^t\big)^{(n-1)/2}=1,\\
\label{hypergeom_prp_e3}
& I_{p,p}(t)=I_{n-1-p,n-1-p}(t) \qquad\forall~p\!=\!0,1,\ldots,n\!-\!1.
\end{alignat}
\end{prp}

\begin{prp}
\label{hypergeom_prp2}
For all $n\!\in\!\Z^+$,
\begin{equation}\label{hypergeom_prp2e1}
\mu(e^t) =\int_0^{e^t}\!\frac{(1\!-\!n^nu)^{-1/n}-1}{u}du
\in e^t\cdot\Q[[e^t]].
\end{equation}
The coefficients of the power series 
\begin{equation}\label{Qdfn_e}
Q(\hb,e^t)\equiv e^{-(t+\mu(e^t))/\hb}R(\hb^{-1},t)
\in\Q(\hb)[[e^t]]
\end{equation}
are holomorphic at $\hb\!=\!0$.
Furthermore,
\begin{alignat}{1}
\label{hypergeom_prp2e2}
\Phi_0(e^t)&\equiv Q(0,e^t)=\big(1\!-\!n^ne^t\big)^{-1/n};\\
\label{hypergeom_prp2e3}
\Phi_1(e^t)&\equiv \frac{d}{d\hb}Q(\hb,e^t)\Big|_{\hb=0}
=\frac{(n\!-\!2)(n\!+\!1)}{24n}
\Big(\big(1\!-\!n^ne^t\big)^{-1/n}-\big(1\!-\!n^ne^t\big)^{-1}\Big).
\end{alignat}\\
\end{prp}

\noindent
Any one of the identities in~\e_ref{hypergeom_prp_e1}-\e_ref{hypergeom_prp_e3}
is implied by the others. We state them all for convenience.
A simple algorithm for determining all coefficients of the expansion of~$Q$
at~$\hb\!=\!0$ is provided by \cite[Theorem~1.5]{ZaZ};
these may be needed for computing higher-genus GW-invariants of
projective CY-hypersurfaces.

\subsection{Linear Independence in Symmetric Rational Functions}
\label{symmalg_subs}

\noindent
In this subsection we prove a lemma showing that most terms appearing in our 
computation of $\F(\al,x,u)$ can be ignored if our only aim is to determine~$\F_0(u)$.\\

\noindent
For each $p\!\in\![n]$, let $\si_p$ be the $p$-th elementary symmetric 
polynomial in $\al_1,\ldots,\al_n$.
Denote by
$$\Q[\al]^{S_n}\equiv\Q[\al_1,\ldots,\al_n]^{S_n}\subset\Q [\al_1,\ldots,\al_n]$$
the subspace of symmetric polynomials, by 
$$\cI\subset \Q[\al]^{S_n}$$ 
the ideal generated by $\si_1,\ldots,\si_{n-1}$, and by
$$\ti\Q[\al]^{S_n} \equiv 
\Q[\al_1,\ldots,\al_n]_{\lr{\al_j,(\al_j-\al_k)|j\neq k}}^{S_n} \subset \Q_{\al}$$
the subalgebra of symmetric rational functions in $\al_1,\ldots,\al_n$
whose denominators are products of $\al_j$ and $(\al_j\!-\!\al_k)$ with $j\!\neq\!k$.
For each $i\!=\!1,\ldots,n$, let 
$$\ti\Q_i[\al]^{S_{n-1}} \equiv 
\Q[\al_1,\ldots,\al_n]_{\lr{\al_i,(\al_i-\al_k)|k\neq i}}^{S_{n-1}}
\subset \Q_{\al}$$
be the subalgebra consisting of rational functions symmetric in $\{\al_k\!:k\!\neq\!i\}$
and with denominators that are  products of $\al_i$ and 
$(\al_i\!-\!\al_k)$ with $k\!\neq\!i$.
Let
\begin{equation}\label{congspacedfn_e}
\K_i\equiv\Span\big\{\cI\cdot\ti\Q_i[\al]^{S_{n-1}},
\{1,\al_i,\ldots,\al_i^{n-3},\al_i^{n-1}\}\cdot\ti\Q[\al]^{S_n}\big\}
\end{equation}
be the linear span (over $\Q$) of $\cI\cdot\ti\Q_i[\al]^{S_{n-1}}$ and 
$\al_i^p\cdot\ti\Q[\al]^{S_n}$ with $p\!=\!0,1,\ldots,n\!-\!3,n\!-\!1$.
If $f,g\!\in\!\Q_{\al}$, we will write
\begin{equation}\label{cong_e1}
f\cong_i g \qquad\hbox{if}\qquad f-g\in\K_i.
\end{equation}

\begin{lmm}
\label{algebr_lmm1}
(i) The ideal $\cI$ does not contain the product of 
any powers of $\si_n$ and
$$D\equiv \prod_{j\neq k}(\al_j\!-\!\al_k).$$
(ii) If $n\!\ge\!2$, the linear span of $\al_i^{n-2}$ is disjoint from $\K_i$:
$$\Span\big\{\al_i^{n-2}\big\}\cap \K_i=\{0\}\subset \Q_{\al}.$$\\
\end{lmm}

\noindent
{\it Proof:} (i) Suppose $m\!\in\bar\Z^+$. 
If $\al_1,\ldots,\al_n$ are the $n$ distinct roots of the polynomial $x^n\!-\!1$,
then
$$\si_1(\al_1,\ldots,\al_n),\ldots,\si_{n-1}(\al_1,\ldots,\al_n)=0,
\qquad\hbox{but}\qquad
\si_n(\al_1,\ldots,\al_n)^mD(\al_1,\ldots,\al_n)^m\neq0.$$
Thus, $\si_n^mD^m\!\not\in\!\cI$.\\

\noindent
(ii) Every polynomial in $\al_1,\ldots,\al_n$ which is symmetric in
$\{\al_k\!:k\!\neq\!i\}$ can be written as a polynomial in 
$\al_i$ with coefficients in~$\Q[\al]^{S_n}$.
Thus, suppose
\begin{gather}\label{indepen_e1}
\al_i^{n-1}=\frac{\sum_{r=0}^{r=N}\al_i^rf_r}
{\al_i^m\prod_{k\neq i}(\al_i\!-\!\al_k)^m}
+\frac{\sum_{p=0}^{n-2}\al_i^p g_p+\al_i^ng_n}{\si_n^m\prod_{j\neq k}(\al_j\!-\!\al_k)^m},\\
\hbox{where}\qquad m\!\in\!\bar\Z^+,~N=n(m\!+\!1)\!-\!2,~
f_r\in\cI,~g_p\!\in\!\Q[\al]^{S_n}. \notag
\end{gather}
Multiplying out the denominators, we obtain
$$\al_i^{n-1}\si_n^mD^m=\sum_{r=0}^{r=\ti{N}}\al_i^rF_r
+\sum_{p=0}^{n-2}\al_i^p g_p+\al_i^ng_n,
\qquad\hbox{where}\qquad \ti{N}=n(mn\!+\!1)\!-\!2,~F_r\in\cI.$$
It follows that 
\begin{equation}\label{indepen_e3}\begin{split}
&\sum_{i=1}^{i=n}\fR_{x=\al_i}
\bigg\{\frac{x^{n-1}}{\prod_{k=1}^{k=n}(x\!-\!\al_k)}\si_n^mD^m\bigg\}\\
&\qquad\qquad\qquad
=\sum_{i=1}^{i=n}\fR_{x=\al_i}\Bigg\{\frac{1}{\prod_{k=1}^{k=n}(x\!-\!\al_k)}
\bigg(\sum_{r=0}^{r=\ti{N}}x^pF_r+\sum_{p=0}^{n-2}x^pg_p+x^ng_n\bigg)\Bigg\}.
\end{split}\end{equation}
By the Residue Theorem on $S^2$, both sides of~\e_ref{indepen_e3} 
equal to the negative of the residues of the corresponding one-forms at $x\!=\!\i$.
Thus,
\begin{equation}\label{indepen_e4}
1\cdot \si_n^mD^m=\sum_{r=n-1}^{r=\ti{N}}\!\!\!\ti{F}_r+\si_1g_n,
\qquad\hbox{where}\qquad \ti{F}_r\in\cI.
\end{equation}
Note that $\ti{F}_r\!\in\!\cI$ because $F_r\!\in\!\cI$.
Thus, \e_ref{indepen_e4} contradicts the first statement of the lemma,
and \e_ref{indepen_e1} cannot hold.\\

\noindent
{\it Remark:} The proof of (ii) of Lemma~\ref{algebr_lmm1} shows that its statement
remains valid if $n\!-\!2$ is replaced by any $p\!=\!0,1,\ldots,n\!-\!1$,
in the definition of $\K_i$ and in the statement of~(ii).

\subsection{The Genus Zero Generating Functions}
\label{g0form_subs}

\noindent
By Subsections~\ref{local_subs}-\ref{localdata_subs2}, 
the generating function $\F(\al,\al_i,u)$ for the reduced genus~$1$ GW-invariants 
of a degree~$n$ hypersurface in $\P^{n-1}$ is given~by
\begin{equation}\label{sumofcontr_e}
\F(\al,\al_i,u)=\cA_i(u)+\sum_{j=1}^{j=n}\ti\cA_{ij}(u)
+\cB_i(u)+\sum_{j=1}^{j=n}\ti\cB_{ij}(u);
\end{equation}
see~\e_ref{Frestr_e} for the definition of~$\F(\al,\al_i,u)$.
Propositions~\ref{contr_prp1} and~\ref{contr_prp2} express
the four terms on the right-hand side of~\e_ref{sumofcontr_e}
in terms of the generating functions for genus~$0$ GW-invariants defined 
in~\e_ref{Z1ptdfn_e}-\e_ref{Z2ptdfn_e}.
These functions have been previously computed in terms of hypergeometric series.
We describe them in this subsection.
For the rest of Section~\ref{alg_sec}, we assume that $n\!\ge\!2$.\\

\noindent
Let
\begin{equation}\label{Ydfn_e}
\cY(\hb,x,e^t) =\frac{1}{I_{0,0}(t)}\sum_{d=0}^{\i}e^{dt}
\frac{\prod_{r=1}^{r=nd}(nx\!+\!r\hb)}
{\prod_{r=1}^{r=d}\big(\prod_{k=1}^{k=n}(x\!-\!\al_k\!+\!r\hb)
-\prod_{k=1}^{k=n}(x\!-\!\al_k)\big)}.
\end{equation}
If $p\!\in\!\bar\Z^+$, let 
\begin{equation}\label{Tderivdfn_e}
\cY_p(\hbar,x,e^t)=e^{-xt/\hb}
\bigg\{\frac{\hb}{I_{p,p}(t)}\frac{d}{dt}\bigg\}\ldots
\bigg\{\frac{\hb}{I_{1,1}(t)}\frac{d}{dt}\bigg\}
\big(e^{xt/\hb}\cY(\hb,x,e^t)\big).
\end{equation}\\
In particular,
\begin{equation}\label{Tderivdfn_e2}
\cY_0(\hbar,x,e^t)=\cY(\hbar,x,e^t) \qquad\hbox{and}\qquad
\cY_1(\hbar,x,e^t)=\bigg\{x\frac{dt}{dT}+\hb\frac{d}{dT}\bigg\}\cY(\hbar,x,e^t),
\end{equation}
if $T$ and $t$ are related by the mirror transformation~\e_ref{mirmap_e}.
Let
\begin{equation}\label{tildeYdfn_e}
\wt\cY(\hb_1,\hb_2,x,e^t)=
x\!\!\!\!\!\!\sum_{\underset{p,q\ge0}{p+q=n-2}}\!\!\!\!\!\!
\cY_p(\hb_1,x,e^t)\cY_q(\hb_2,x,e^t)
+x^{-(n-1)}\cY_{n-1}(\hb_1,x,e^t)\cY_{n-1}(\hb_2,x,e^t).
\end{equation}

\begin{lmm}
\label{Z1pt_lmm}
The power series $\cZ_i^*(\hb,u)$ is rational in $\hb\!\in\!S^2$ and 
vanishes to second order at $\hb\!=\!\i$.
It has simple poles at $\hb\!=\!(\al_j\!-\!\al_i)/d$, with $j\!\neq\!i$ and $d\!\in\!\Z^+$,
and another pole at $\hb\!=\!0$.
Furthermore,
\begin{equation}\label{Z1pt_lmm_e1}
\fR_{\hb=(\al_j-\al_i)/d}\big\{\cZ_i^*(\hb,u)\big\}
=
\sum_{\Ga}\int_{\cZ_{\Ga}}\frac{\E(\V_0')\ev_1^*\phi_i}{\E\big(\N\cZ_{\Ga}\big)},
\end{equation}
where the sum is taken over the two-pointed trees $\Ga$ as in~\e_ref{decortgraphdfn_e}
such that the marked point $1$ is attached to a vertex $v_0\!=\!\eta(1)$ of valence~$2$,
$\mu(v_0)\!=\!i$, $\mu(v)\!=\!j$ for the unique vertex $v$ adjacent to~$v_0$,
and $\d(\{v_0,v\})\!=\!d$.
Finally, for all $z\!=\!0,\i,-n\al_i$ and $a\!\in\!\Z$,
\begin{equation}\label{Z1pt_lmm_e2}
\fR_{\hb=z}\Big\{\frac{\hb^a}{n\al_i\!+\!\hb}\big(1\!+\!\cZ_i^*(\hb,e^T)
-e^{(t-T)\al_i/\hb}\cY(\hb,\al_i,e^t)\big)\Big\}
\in \big(\cI\cdot\ti\Q_i[\al]^{S_{n-1}}\big)\big[\big[e^t\big]\big].
\end{equation}
\end{lmm}

\begin{lmm}
\label{Z1pt2_lmm}
The power series $\hb\cZ_{ji}^*(\hb,u)$ is rational in $\hb\!\in\!S^2$ and 
vanishes at $\hb\!=\!\i$.
It has simple poles at $\hb\!=\!(\al_l\!-\!\al_j)/d$, with $l\!\neq\!j$ and $d\!\in\!\Z^+$,
and another pole at $\hb\!=\!0$.
Furthermore,
\begin{equation}\label{Z1pt2_lmm_e1}
\fR_{\hb=(\al_l-\al_j)/d}\big\{\hb\cZ_{jl}^*(\hb,u)\big\}
=
\sum_{\Ga}\int_{\cZ_{\Ga}}\frac{\E(\V_0')\ev_1^*\phi_j\ev_2^*\phi_i}{\E\big(\N\cZ_{\Ga}\big)},
\end{equation}
where the sum is taken over the two-pointed trees $\Ga$ as in~\e_ref{decortgraphdfn_e}
such that the marked point $1$ is attached to a vertex $v_0\!=\!\eta(1)$ of valence~$2$,
$\mu(v_0)\!=\!j$, $\mu(v)\!=\!l$ for the unique vertex $v$ adjacent to~$v_0$,
and $\d(\{v_0,v\})\!=\!d$.
Finally, for all $z\!=\!0,\i,-n\al_j$ and $a\!\in\!\Z$,
\begin{gather}\label{Z1pt2_lmm_e2}
\fR_{\hb=z}\Big\{\frac{\hb^a}{n\al_j\!+\!\hb}\big(\al_i^{n-2}\al_j\!+\!\hb\cZ_{ji}^*(\hb,e^T)
-\al_i^{n-2}e^{(t-T)\al_j/\hb}\cY_1(\hb,\al_j,e^t)\big)\Big\}
\in  \K_{ij}\big[\big[e^t\big]\big],\\
\hbox{where}\qquad
\K_{ij}=\big\{1,\al_i,\ldots,\al_i^{n-3},\al_i^{n-1}\}\cdot\ti\Q_j[\al]^{S_{n-1}}
\oplus\al_i^{n-2}\cI\cdot\ti\Q_j[\al]^{S_{n-1}}.\notag
\end{gather}
\end{lmm}

\begin{lmm}
\label{Z2pt_lmm}
The power series $\hb_1\hb_2\wt\cZ_{ii}^*(\hb_1,\hb_2,u)$ is rational in $\hb_1\!\in\!S^2$ and 
vanishes at $\hb_1\!=\!\i$.
It has simple poles at $\hb_1\!=\!(\al_j\!-\!\al_i)/d$, with $j\!\neq\!i$ and $d\!\in\!\Z^+$,
and another pole at $\hb_1\!=\!0$.
Furthermore,
\begin{equation}\label{Z2pt_lmm_e1}
\fR_{\hb_1=(\al_j-\al_i)/d}\big\{2\hb_1\hb_2\wt\cZ_{ii}^*(\hb_1,\hb_2,u)\big\}
=\sum_{\Ga}\int_{\cZ_{\Ga}}
\frac{\E(\V_0')\ev_1^*\phi_i\ev_2^*\phi_i}{(\hb_2\!-\!\psi_2)\E\big(\N\cZ_{\Ga}\big)},
\end{equation}
where the sum is taken over the two-pointed trees $\Ga$ as in~\e_ref{decortgraphdfn_e}
such that the marked point $1$ is attached to a vertex $v_0\!=\!\eta(1)$ of valence~$2$,
$\mu(v_0)\!=\!i$, $\mu(v)\!=\!j$ for the unique vertex $v$ adjacent to~$v_0$,
and $\d(\{v_0,v\})\!=\!d$.
The analogous statements hold for~$\hb_2$.
Finally, for all $a_1,a_2\!\in\!\Z^-$,
\begin{equation}\label{Z2pt_lmm_e2}\begin{split}
&\fR_{\hb_1=0}\fR_{\hb_2=0}\Bigg\{\hb_1^{a_1}\hb_2^{a_2}
\bigg(2\hb_1\hb_2\wt\cZ_{ii}^*(\hb,e^T)\\
&\hspace{1.5in}
-\frac{e^{(t-T)\al_i(h_1^{-1}+\hb_2^{-1})}}{\hb_1\!+\!\hb_2}\wt\cY(\hb_1,\hb_2,\al_i,e^t)
\bigg)\Bigg\} 
\in \big(\cI\cdot\ti\Q_i[\al]^{S_{n-1}}\big)\big[\big[e^t\big]\big].
\end{split}\end{equation}\\
\end{lmm}

\noindent
All statements concerning rationality of $\cZ_i^*$, $\cZ_{ji}^*$, and $\cZ_{ii}^*$
in these three lemmas refer to rationality of the coefficients of the powers of~$e^t$.
Lemma~\ref{Z1pt_lmm} is proved in \cite[Chapter~30]{MirSym};
Lemmas~\ref{Z1pt2_lmm} and~\ref{Z2pt_lmm} are proved 
in~\cite{bcov0}.\footnote{In fact, Lemma~\ref{Z1pt_lmm} is essentially the main result of
\cite[Chapter~30]{MirSym}, while Lemmas~\ref{Z1pt2_lmm} and~\ref{Z2pt_lmm} are essentially
the main results of~\cite{bcov0}.}
For example, the conclusion of \cite[Section~30.4]{MirSym} is that 
$$1\!+\!\cZ_i^*(\hb,e^T)=e^{C(e^t)\si_1/\hb}e^{(t-T)\al_i/\hb}\cY(\hb,\al_i,e^t),
\qquad\hbox{where}\quad
C(u)=-\sum_{d=1}^{\i}u^d\bigg(\frac{(nd)!}{(d!)^n}\sum_{r=1}^{r=d}\frac{1}{r}\bigg).$$
This statement clearly implies \e_ref{Z1pt_lmm_e2}, provided $n\!\ge\!2$ 
so that $\si_1\!\in\!\cI$.\\

\noindent
The differences in~\e_ref{Z1pt_lmm_e2}, \e_ref{Z1pt2_lmm_e2}, and~\e_ref{Z2pt_lmm_e2} 
are of course symmetric in the $\al$'s in every appropriate sense.
For example, any of the differences in~\e_ref{Z1pt_lmm_e2} is the evaluation at $x\!=\!\al_i$ 
of a power series in $e^t$ with coefficients in the rational functions in $x$, 
$\si_1,\ldots,\si_{n-1}$.
This is immediate from the explicit formulas for $\cZ_i^*$, $\cZ_{ji}^*$, and 
$\wt\cZ_{ij}$ in \cite[Chapter~30]{MirSym} and in~\cite{bcov0}.
This symmetry is used in the next subsection in the computation of 
the contributions of $\ti\cA_{ij}$ and~$\ti\cB_{ij}$.

\subsection{Proof of Theorem~\ref{main_thm}}
\label{mainthmpf_subs}

\noindent
We will use Lemma~\ref{algebr_lmm1}, along with Lemmas~\ref{Z1pt_lmm}-\ref{Z2pt_lmm},
to extract the coefficients of $\al_i^{n-2}$  from the expressions of 
Propositions~\ref{contr_prp1} and~\ref{contr_prp2}  modulo~$\K_i[[u]]$.
In the notation of Theorem~\ref{main_thm}, the two coefficients are
$\frac{d}{dT}\ti{A}(e^t)$ and $\frac{d}{dT}\ti{B}(e^t)$.
Let $\cong_i$ be as in~\e_ref{cong_e1}, 
with its meaning extended to power series in $e^t$ in the natural way.\\

\noindent
We begin by defining the analogues of the power series $\cY(\hb,x,u)$ and 
$\cY_p(\hb,x,u)$ without the $\al$'s.
Let
\begin{equation}\label{Ysimpdfn_e}\
Y(\hb,x,e^t) =\frac{1}{I_{0,0}(t)}\sum_{d=0}^{\i}e^{dt}
\frac{\prod_{r=1}^{r=nd}(nx\!+\!r\hb)}
{\prod_{r=1}^{r=d}\big((x\!+\!r\hb)^n-x^n\big)}.
\end{equation}
If $p\!\in\!\bar\Z^+$, let 
\begin{equation}\label{Tsimpderivdfn_e}
Y_p(\hbar,x,e^t)=e^{-xt/\hb}
\bigg\{\frac{\hb}{I_{p,p}(t)}\frac{d}{dt}\bigg\}\ldots
\bigg\{\frac{\hb}{I_{1,1}(t)}\frac{d}{dt}\bigg\}
\big(e^{xt/\hb}Y(\hb,x,e^t)\big).
\end{equation}\\
Similarly to~\e_ref{Tderivdfn_e2}
\begin{equation}\label{Tsimpderivdfn_e2}
Y_0(\hbar,x,e^t)=Y(\hbar,x,e^t) \qquad\hbox{and}\qquad
Y_1(\hbar,x,e^t)=\bigg\{x\frac{dt}{dT}+\hb\frac{d}{dT}\bigg\}Y(\hbar,x,e^t),
\end{equation}
if $T$ and $t$ are related by the mirror transformation~\e_ref{mirmap_e}.
Let
\begin{equation}\label{tildeYsimpdfn_e}
\wt{Y}(\hb_1,\hb_2,x,e^t)=
x\!\!\!\!\!\!\sum_{\underset{p,q\ge0}{p+q=n-2}}\!\!\!\!\!\!
Y_p(\hb_1,x,e^t)Y_q(\hb_2,x,e^t)
+x^{-(n-1)}Y_{n-1}(\hb_1,x,e^t)Y_{n-1}(\hb_2,x,e^t).
\end{equation}\\

\noindent
We note that 
\begin{equation}\label{symmetr_e}
\bigg(\prod_{k=1}^{k=n}(x\!-\!\al_k\!+\!r\hb)
-\prod_{k=1}^{k=n}(x\!-\!\al_k)\bigg)
-\big((x\!+\!r\hb)^n-x^n\big)\in\cI[\hb,x].
\end{equation}
It follows that 
\begin{alignat}{1}\label{resdiff_e}
&\fR_{\hb=z}\bigg\{\frac{\hb^a}{n\al_i\!+\!\hb}\cY(\hb,\al_i,u)\bigg\}-
\Bigg(\fR_{\hb=z'}\bigg\{\frac{\hb^a}{nx\!+\!\hb}Y(\hb,x,u)\bigg\}\Bigg)\bigg|_{x=\al_i}\\
&\quad
=\fR_{\hb=z}\bigg\{\frac{\hb^a}{n\al_i\!+\!\hb}\cY(\hb,\al_i,u)\bigg\}-
\fR_{\hb=z}\bigg\{\frac{\hb^a}{n\al_i\!+\!\hb}Y(\hb,\al_i,u)\bigg\}
\in \big(\cI\cdot\ti\Q_i[\al]^{S_{n-1}}\big)\big[\big[u\big]\big]\notag
\end{alignat}
for all $a\!\in\!\Z$, $z\!=\!0,\i,-n\al_i$ and the corresponding $z'\!=\!0,\i,-nx$.
The equality in~\e_ref{resdiff_e} holds because the evaluation at $x\!=\!\al_i$
of the residue of $\hb^aY(\hb,x,u)/(nx\!+\!\hb)$ at $\hb\!=\!z'$ is defined, 
i.e.~the evaluation of $Y(\hb,x,u)$ at $x\!=\!\al_i$ does not change
the order of the pole of (the coefficient of each $u^d$ in) $Y(\hb,x,u)$.
This is  the reason we modify the denominators of \cite[Chapter~30]{MirSym} and 
of~\cite{bcov0} by subtracting off $\prod_{k=1}^{k=n}\!(x\!-\!\al_i)$.
This modification has no effect on the evaluation maps $x\!\lra\!\al_i$.\\

\noindent
We now use Lemmas~\ref{Z1pt_lmm}-\ref{Z2pt_lmm} and \e_ref{resdiff_e}
to extract the ``relevant'' part from the expressions of Proposition~\ref{contr_prp1} and~\ref{contr_prp2}.
If $\eta_i(u)$ and $\mu(u)$ are as in~\e_ref{etadfn_e} and~\e_ref{mudfn_e0},
respectively, then
\begin{equation}\label{etadiff_e}
\eta_i(e^T)-\big(t-T+\mu(e^t)\big)\al_i\in  
\big(\cI\cdot\ti\Q_i[\al]^{S_{n-1}}\big)\big[\big[e^t\big]\big]
\end{equation}
by~\e_ref{Z1pt_lmm_e2} and~\e_ref{resdiff_e}.
Similarly, if $\Phi_0(\al_i,u)$ and $\Phi_0(u)$ are as in~\e_ref{Phi0dfn_e}
and~\e_ref{hypergeom_prp2e2}, respectively, then
\begin{equation}\label{Phi0diff_e}
\Phi_0(\al_i,e^T)-\frac{\Phi_0(e^t)}{I_{0,0}(e^t)}
\in  \big(\cI\cdot\ti\Q_i[\al]^{S_{n-1}}\big)\big[\big[e^t\big]\big]
\end{equation}
by~\e_ref{Z1pt_lmm_e2}, \e_ref{resdiff_e}, and \e_ref{etadiff_e}.
By~\e_ref{contr_prp1e1}, \e_ref{Z2pt_lmm_e2}, \e_ref{etadiff_e}, 
and~\e_ref{Phi0diff_e},
\begin{equation}\label{Adiff_e}
\cA_i(e^T)-A(e^t)\al_i^{n-2}\in\big(\cI\cdot\ti\Q_i[\al]^{S_{n-1}}\big)\big[\big[e^t\big]\big],
\end{equation}
where
\begin{equation}\label{Adiff_e2}\begin{split}
A(e^t)=\frac{I_{0,0}(e^t)}{2\Phi_0(e^t)}
\fR_{h_1=0}\fR_{h_2=0}\bigg\{\hb_1^{-1}\hb_2^{-1}
\frac{e^{-\mu(e^t)(\hb_1^{-1}+\hb_2^{-1})}}{\hb_1\!+\!\hb_2}
\wt{Y}(\hb_1,\hb_2,1,e^t)\bigg\}. 
\end{split}\end{equation}\\

\noindent
On the other hand, by~\e_ref{Z1pt2_lmm_e2}, \e_ref{resdiff_e}, and~\e_ref{etadiff_e},
\begin{equation}\label{pt12diff_e}
\fR_{\hb=0}\big\{e^{-\eta_j(e^T)/\hb}\cZ_{ji}^*(\hb,e^T)\big\}
\cong_i \big(-1+\De(e^t)\big)\al_i^{n-2}\al_j,
\end{equation}
where 
\begin{equation}\label{pt12diff_e2}\begin{split}
\De(e^t)&\equiv\fR_{\hb=0}\Bigg\{\hb^{-1}e^{-\mu(e^t)/\hb}Y_1(\hb,1,e^t)\Bigg\}\\
&=\fR_{\hb=0}\bigg\{\hb^{-1}
\Big\{\frac{d(t\!+\!\mu(e^t))}{dT}+\hb\frac{d}{dT}\Big\}
\big(e^{-\mu(e^t)/\hb}Y(\hb,1,e^t)\big)\Bigg\}\\
&=\frac{d(t\!+\!\mu(t))}{dT}
\fR_{\hb=0}\Big\{\hb^{-1}e^{-\mu(e^t)/\hb}Y(\hb,1,e^t)\Big\}
=\frac{d(t\!+\!\mu(e^t))}{dT}\frac{\Phi_0(e^t)}{I_{0,0}(t)}.
\end{split}\end{equation}
The first equality in~\e_ref{pt12diff_e2} follows from~\e_ref{Tsimpderivdfn_e2},
the second from the holomorphicity statement of Proposition~\ref{hypergeom_prp2},
and the last from~\e_ref{Ysimpdfn_e} and~\e_ref{hypergeom_prp2e2}.
By~\e_ref{contr_prp1e2}, \e_ref{Adiff_e}, and~\e_ref{pt12diff_e},
\begin{alignat}{1}\label{pt12diff_e2b}
\sum_{j=1}^{j=n}\ti\cA_{ij}(e^T)
&\cong_i \al_i^{n-2}\big(-1+\De(e^t)\big)
\sum_{j=1}^{j=n}\frac{\al_j\cA_j(e^T)}{\prod_{k\neq j}(\al_j\!-\!\al_k)} \notag\\
&\cong_i \al_i^{n-2}\big(-1+\De(e^t)\big)A(e^t)
\sum_{j=1}^{j=n}\fR_{z=\al_j}\Bigg\{\frac{z^{n-1}}{\prod_{k=1}^{k=n}(z\!-\!\al_k)}\Bigg\}\\
&=-\al_i^{n-2}\big(-1+\De(e^t)\big)A(e^t)
\fR_{z=\i}\Bigg\{\frac{z^{n-1}}{\prod_{k=1}^{k=n}(z\!-\!\al_k)}\Bigg\}
=\al_i^{n-2}\big(-1+\De(e^t)\big)A(e^t).\notag
\end{alignat}
The first equality above follows from the Residue Theorem on~$S^2$.
Thus, by~\e_ref{Adiff_e} and~\e_ref{pt12diff_e2},
\begin{equation}\label{Aterm_e1}
\cA_i(e^T)+\sum_{j=1}^{j=n}\ti\cA_{ij}(e^T)
\cong_i \al_i^{n-2}\frac{d(t\!+\!\mu(t))}{dT}\frac{\Phi_0(e^t)}{I_{0,0}(t)}A(e^t),
\end{equation}
with $A(e^t)$ defined by~\e_ref{Adiff_e2}.\\

\noindent
We next reduce the right-hand side of~\e_ref{Adiff_e2} to 
the explicit form of Theorem~\ref{main_thm}.
Let
\begin{equation}\label{notat_e}
L(e^t)=\big(1-n^ne^t\big)^{1/n} \qquad\hbox{and}\qquad
f_p(e^t)=\frac{1}{L(e^t)I_{p,p}(t)} ~~~\forall\,p\!=\!0,1,\ldots.
\end{equation}
By~\e_ref{Qdfn_e}, \e_ref{Tsimpderivdfn_e}, and~\e_ref{hypergeom_prp2e1},
\begin{equation*}\begin{split}
e^{-\mu(e^t)/\hb}Y_p(\hb,1,e^t)&=
\frac{1}{I_{p,p}(t)}\bigg\{\frac{d(t\!+\!\mu(t))}{dt}+\hb\frac{d}{dt}\bigg\}  
\ldots
\frac{1}{I_{1,1}(t)}\bigg\{\frac{d(t\!+\!\mu(t))}{dt}+\hb\frac{d}{dt}\bigg\}  
\Big(\frac{Q(\hb,e^t)}{I_{0,0}(t)}\Big) \\
&=f_p(e^t)\bigg\{1+\hb L(e^t)\frac{d}{dt}\bigg\}   \ldots
f_1(e^t)\bigg\{1+\hb L(e^t)\frac{d}{dt}\bigg\}\Big(\frac{Q(\hb,e^t)}{I_{0,0}(t)}\Big)
\end{split}\end{equation*}
for all $p\!\ge\!0$. 
Thus, by \e_ref{hypergeom_prp2e2}, \e_ref{hypergeom_prp2e3}, 
and the regularity statement of Proposition~\ref{hypergeom_prp2},
\begin{alignat}{1}\label{Aterm_e4a}
\Th_p^{(0)}(e^t)&\equiv
\fR_{h=0}\bigg\{\hb^{-1}e^{-\mu(e^t)/\hb}Y_p(\hb,1,e^t)\bigg\}
=\prod_{r=0}^{r=p}f_r(e^t);\\
\Th_p^{(1)}(e^t)&\equiv
\fR_{h=0}\bigg\{\hb^{-2}e^{-\mu(e^t)/\hb}Y_p(\hb,1,e^t)\bigg\}\notag\\
\label{Aterm_e4b}
&=L(e^t)\bigg(\prod_{r=0}^{r=p}f_r(e^t)\bigg)\bigg(\Phi_1(e^t)+
\sum_{r=0}^{p-1}(p\!-\!r)\frac{f_r'(e^t)}{f_r(e^t)}\bigg),
\end{alignat}
where $'$ denotes the derivative with respect to $t$.
Note that by~\e_ref{hypergeom_prp_e1} and~\e_ref{hypergeom_prp_e2},  
$$\Th_{n-1}^{(0)}(e^t)=1    \qquad\hbox{and}\qquad
\Th_{n-1}^{(1)}(e^t)=L(e^t)\Phi_1(e^t).$$
Thus, by \e_ref{tildeYsimpdfn_e}, \e_ref{Aterm_e4a}, and~\e_ref{Aterm_e4b},
\begin{equation}\label{Aterm_e5}\begin{split}
\fR_{h_1=0}\fR_{h_2=0}\bigg\{
\frac{e^{-\mu(e^t)(\hb_1^{-1}+\hb_2^{-1})}}{\hb_1\hb_2(\hb_1\!+\!\hb_2)}
\wt{Y}(\hb_1,\hb_2,1,e^t)\bigg\}
&=\!\!\!\!\!\!\sum_{\underset{p,q\ge0}{p+q=n-2}}\!\!\!\!\!\!
\Th_p^{(1)}(e^t)\Th_q^{(0)}(e^t)+\Th_{n-1}^{(1)}(e^t)\Th_{n-1}^{(0)}(e^t)\\
&= L(e^t)\bigg(n\Phi_1(e^t)+
\sum_{p=0}^{n-2}\sum_{r=0}^{p-1}(p\!-\!r)\frac{f_r'(e^t)}{f_r(e^t)}\bigg).
\end{split}\end{equation}
The last equality uses \e_ref{hypergeom_prp_e3}, followed by \e_ref{hypergeom_prp_e1}.\\

\noindent
By~\e_ref{Aterm_e1}, \e_ref{Adiff_e2}, \e_ref{Aterm_e5}, \e_ref{mirmap_e}, 
\e_ref{hypergeom_prp2e1}, and \e_ref{hypergeom_prp2e3},
the contribution of the fixed loci of Proposition~\ref{contr_prp1} is
\begin{equation*}\begin{split}
&\frac{1}{2I_{1,1}(t)}\Bigg(\frac{(n\!-\!2)(n\!+\!1)}{24}\Big(L(e^t)^{-1}-L(e^t)^{-n}\Big)
+\sum_{p=0}^{n-2}\binom{n\!-\!1\!-\!p}{2}\frac{f_p'(e^t)}{f_p(e^t)}\Bigg)\\
&\qquad 
=\frac{1}{2I_{1,1}(t)}\frac{d}{dt}\Bigg(
\frac{(n\!-\!2)(n\!+\!1)}{24}\Big(\big(t\!+\!\mu(e^t)\big)-
\big(t-n\ln L(e^t)\big)\Big)+\sum_{p=0}^{n-3}\binom{n\!-\!1\!-\!p}{2}\ln f_p(e^t)\Bigg)\\
&\qquad
=\frac{1}{2}\frac{d}{dT}\Bigg(\frac{(n\!-\!2)(n\!+\!1)}{24}\mu(e^t)
-\frac{(n\!-\!2)(3n\!-\!5)}{24}\ln\big(1\!-\!n^ne^t\big)
-\sum_{p=0}^{n-3}\binom{n\!-\!1\!-\!p}{2}\ln I_{p,p}(t)\Bigg).
\end{split}\end{equation*}
The first equality above uses \e_ref{notat_e} and~\e_ref{hypergeom_prp2e1};
the second one also uses~\e_ref{Tden_e} and~\e_ref{mirmap_e}.
We have now proved the first statement of Theorem~\ref{main_thm};
the second form of the expression in~\e_ref{main_thm_e1} is easily obtainable
from the first using \e_ref{hypergeom_prp_e1}-\e_ref{hypergeom_prp_e3}.\\

\noindent
We compute the contributions of the terms $\cB_i(u)$
and $\ti\cB_{ij}(u)$ of Proposition~\ref{contr_prp2} to~$\F_0(u)$ similarly.
However, before proceeding, we observe that the $\al$-free analogue of 
the term $\prod_{k=1}^{k=n}(\al_i\!-\!\al_k\!+\!\hb)$ in the numerators 
in~\e_ref{contr_prp2e1} and~\e_ref{contr_prp2e2} is $(x\!+\!\hb)^n\!-\!x^n$.
The reason is the $r\!=\!1$ case of~\e_ref{symmetr_e} and
that subtracting off $\prod_{k=1}^{k=n}(x\!-\!\al_k)$ from the numerators
has no effect on the evaluation maps $x\!\lra\!\al_i$.\\

\noindent
By~\e_ref{contr_prp2e1}, \e_ref{Z1pt_lmm_e2}, and~\e_ref{resdiff_e},
\begin{gather*}
\cB_i(e^T)\cong_i \al_i^{n-2}B(e^t), \qquad\hbox{where}\\
B(e^t)=\frac{n}{24} \fR_{\hb=0,\i,-n}\bigg\{
\frac{(1\!+\!\hb)^n-1}{(n\!+\!\hb)\hb^3}\,
\frac{e^{(t-T)/\hb}Y(\hb,1,e^t)-1}{e^{(t-T)/\hb}Y(\hb,1,e^t)}\bigg\}.
\end{gather*}
On the other hand, by~\e_ref{contr_prp2e2}, \e_ref{Z1pt_lmm_e2}, \e_ref{Z1pt2_lmm_e2}, 
\e_ref{resdiff_e}, and~\e_ref{Tsimpderivdfn_e2},
\begin{equation*}\begin{split}
\ti\cB_{ij}(e^T) &\cong_i -\frac{n\al_j\al_i^{n-2}}{24\prod_{k\neq j}(\al_j\!-\!\al_k)}
\fR_{\hb=0,\i,-n\al_j}\bigg\{
\frac{(\al_j+\!\hb)^n\!-\!\al_j^n}{(n\al_j\!+\!\hb)\hb^3}\,
\frac{e^{(t-T)\al_j/\hb}Y_1(\hb,\al_j,e^t)-\al_j}
{e^{(t-T)\al_j/\hb}Y(\hb,\al_j,e^t)}\bigg\}\\
&=-\frac{n\al_i^{n-2}\al_j^{n-1}}{24\prod_{k\neq j}(\al_j\!-\!\al_k)}\fR_{\hb=0,\i,-n}\bigg\{
\frac{(1+\!\hb)^n\!-\!1}{(n\!+\!\hb)\hb^3}\,
\frac{\big\{1+\hb\frac{d}{dT}\big\}\big(e^{(t-T)/\hb}Y(\hb,1,e^t)\big)-1}
{e^{(t-T)/\hb}Y(\hb,1,e^t)}\bigg\}\\
&=\frac{\al_i^{n-2}\al_j^{n-1}}{\prod_{k\neq j}(\al_j\!-\!\al_k)}
\Big(-B(e^t)+\frac{d}{dT}\ti{B}(e^t)\Big),
\end{split}\end{equation*}
where
\begin{equation}\label{Bterm_e2}\begin{split}
\ti{B}(e^t)&\equiv \ti{B}_0(e^t)+\ti{B}_{\i}(e^t)+\ti{B}_{-n}(e^t),\\
\ti{B}_z(e^t)&=-\frac{n}{24}\fR_{\hb=z}\bigg\{
\frac{(1+\!\hb)^n\!-\!1}{(n\!+\!\hb)\hb^2}\,\ln\Big(e^{(t-T)/\hb}Y(\hb,1,e^t)\Big)\bigg\}.
\end{split}\end{equation}
Thus, applying the Residue Theorem on $S^2$ as in~\e_ref{pt12diff_e2b},
we obtain
\begin{equation}\label{Bterm_e3}
\cB_i(e^T)+\sum_{j=1}^{j=n}\ti\cB_{ij}(e^T)\cong_i\al_i^{n-2}\frac{d}{dT}\ti{B}(e^t).
\end{equation}\\

\noindent
It remains to compute the three residues $\ti{B}_z(e^t)$.
For $z\!=\!-n$, the pole is simple.
Since 
$$e^{t/h}Y(\hb,1,e^t)=\bar{R}(\hb^{-1},t)
\qquad\hbox{and}\qquad \bar{R}(-n^{-1},t)=e^{-t/n}\big/I_{0,0}(t)$$
with $\bar{R}$ as in Subsection~\ref{lowdeg_subs}, we obtain
\begin{equation}\label{Bterm_e5}
\ti{B}_{-n}(e^t)=-\frac{n}{24}\cdot\frac{(1\!-\!n)^n-1}{n^2}
\bigg(\frac{t\!-\!T}{-n}-\ln I_{0,0}(t)\bigg).
\end{equation}
The residue at $z\!=\!0$ is computed using Proposition~\ref{hypergeom_prp2}:
\begin{equation}\label{Bterm_e6}
\ti{B}_0(e^t)=-\frac{n}{24}\bigg(
\frac{(n\!-\!2)(n\!+\!1)}{2n}\big(t\!-\!T\!+\!\mu(e^t)\big)
+\ln\Phi_0(e^t)-\ln I_{0,0}(e^t)\bigg).
\end{equation}
Finally, note that
$$e^{wt}Y(w,1,e^t) =\bar{R}(w,t)
=1+Tw+\sum_{q=2}^{\i}\frac{I_{0,q}(t)}{I_{0,0}(t)}w^q.$$
Thus,
\begin{equation}\label{Bterm_e7}\begin{split}
\ti{B}_{\i}(e^t) &=\frac{n}{24}\fR_{w=0}
\bigg\{\frac{(1\!+\!w)^n}{(1\!+\!nw)w^{n-1}}\big(-Tw+\ln \bar{R}(w,t)\big)\bigg\}\\
&=\frac{n}{24}\Bigg\{-T\,\cD_w^{n-3}\bigg(\frac{(1\!+\!w)^n}{(1\!+\!nw)}\bigg)
+\cD_w^{n-2}\bigg(\frac{(1\!+\!w)^n}{(1\!+\!nw)}\ln\bar{R}(w,t)\bigg)\Bigg\}\\
&=\frac{n}{24}
\sum_{p=2}^{n-2}\bigg(\cD_w^{n-2-p}\frac{(1\!+\!w)^n}{(1\!+\!nw)}\bigg)
\big(\cD_w^p\ln\bar{R}(w,t)\big).
\end{split}\end{equation}
The remaining statement of Theorem~\ref{main_thm} is obtained by adding
up~\e_ref{Bterm_e5}-\e_ref{Bterm_e7}.

\appendix

\section{Some Combinatorics}
\label{comb_app}

\noindent
This appendix contains the computational steps omitted in the proof of 
Lemma~\ref{regular_lmm}, as well as the proof of Lemma~\ref{prodres_lmm}.

\begin{lmm}
\label{comb_l0}
For all $b\!\ge\!0$, $N\!\ge\!0$, and $q_1,\ldots,q_N\!\ge\!0$,
\begin{equation}\label{comb_l0e1}
\sum_{\underset{\be_l\ge0}{\sum_{l=1}^{l=N}\be_l=b}}
\prod_{l=1}^{l=N}\binom{q_l}{\be_l}=\binom{q_1\!+\!\ldots\!+q_N}{b}.
\end{equation}
For all $q\!\ge\!0$ and $a\!\ge\!1$, 
\begin{equation}\label{comb_l0e2}
\sum_{b=0}^{\i}(-1)^b\binom{q}{b}\frac{1}{a\!+\!b}
=\frac{(a\!-\!1)!q!}{(a\!+\!q)!}.
\end{equation}
For all $q\!\ge\!0$ and $a,s\!\ge\!0$, 
\begin{equation}\label{comb_l0e3}
\sum_{b=0}^{\i}(-1)^b\binom{q}{b}\prod_{r=a-s+1}^{r=a}\!\!\!\!\!\!(r\!+\!b)
=(-1)^qs!\binom{a}{s\!-\!q}.
\end{equation}\\
\end{lmm}

\noindent
{\it Proof:} (1) Each summand on the left side of \e_ref{comb_l0e1} is the total number
of ways to choose $\be_l$ elements from a $q_l$-element set for $l\!=\!1,\ldots,N$.
Thus, the number on the left side of \e_ref{comb_l0e1} is the number of ways to choose
$b\!=\!\sum_{l=1}^{l=N}\!b_l$ elements from a set with $\sum_{l=1}^{l=N}\!q_l$ elements.\\

\noindent
(2) The identity~\e_ref{comb_l0e2} is satisfied for $q\!=\!0$.
Suppose \e_ref{comb_l0e2} holds for all $a\!\ge\!1$ and some $q\!\ge\!0$.
Then,
\begin{equation*}\begin{split}
\sum_{b=0}^{\i}\frac{(-1)^b}{a\!+\!b}\binom{q\!+\!1}{b}
&=\sum_{b=0}^{\i}\frac{(-1)^b}{a\!+\!b}
\bigg(\binom{q}{b}+\binom{q}{b\!-\!1}\bigg)\\
&=\sum_{b=0}^{\i}\frac{(-1)^b}{a\!+\!b}\binom{q}{b}
-\sum_{b=0}^{\i}\frac{(-1)^b}{a\!+\!1\!+\!b}\binom{q}{b}\\
&=\frac{(a\!-\!1)!\,q!}{(a\!+\!q)!}-\frac{a!\,q!}{(a\!+\!1\!+\!q)!}
=\frac{(a\!-\!1)!\,(q\!+\!1)!}{(a\!+\!q\!+\!1)!},
\end{split}\end{equation*}
as needed.\\

\noindent
(3) With $q$, $a$, and $s$ as in~\e_ref{comb_l0e3},
\begin{equation*}\begin{split}
\sum_{b=0}^{\i}(-1)^b\binom{q}{b}\prod_{r=a-s+1}^{r=a}\!\!\!\!\!\!(r\!+\!b)
&=\bigg\{\frac{d}{dx}\bigg\}^s\bigg(\sum_{b=0}^{\i}(-1)^b\binom{q}{b}x^{a+b}\bigg)
\bigg|_{x=1}\\
&=\bigg\{\frac{d}{dx}\bigg\}^s\big((1\!-\!x)^qx^a\big)\Big|_{x=1}
=\binom{s}{q}(-1)^q q!\cdot\frac{a!}{(a\!-\!s\!+\!q)!},
\end{split}\end{equation*}
as claimed.\\

\noindent
For each $k\!\in\!\Z$, let
$$[k]=\big\{l\!\in\!\Z^+\!:1\!\le\!l\!\le\!k\big\}$$
as before.
If $\al\!=\!(\al_1,\ldots,\al_k)$ is a tuple of nonnegative integers, define
$$|\al|=\sum_{i=1}^{i=k}\al_i, \qquad
\al!=\prod_{i=1}^{i=k}\!\al_i!\,, \qquad
S(\al)=\big\{(i,j)\!:\,j\!=\!1,\ldots,\al_i;~i\!=\!1,\ldots,k\big\}.$$
If $b$ is an integer, possibly negative, let
$$\binom{|\al|\!+\!b}{\al,b}
=\binom{|\al|+b}{\al_1,\ldots,\al_k,b}.$$
Denote by $(\bar\Z^+)^{\eset}$ the $0$-dimensional lattice.\\

\noindent
{\it Proof of~\e_ref{regprp_e1}:}
Suppose $k\!\in\!\bar\Z^+$, $\al\!\in\!(\Z^+)^k$, and $q\!\in\!(\bar\Z^+)^k$
is a tuple of distinct nonnegative integers.
With $C_q$ as in~\e_ref{bZexp_e}, we will compare the coefficients of 
$$C_q^{\al}\equiv \prod_{i=1}^{i=k}C_{q_i}^{\al_i}$$
on the two sides of~\e_ref{regprp_e1}.
By~\e_ref{Rdfn_e}, for each $\be\!\in\!(\bar\Z^+)^{S(\al)}$ and every choice of $k$
distinguished disjoint subsets of $[a\!+\!2\!+\!|\be|]$ of cardinalities
$\al_1,\ldots,\al_k$, the term $C_q^{\al}$ appears in the $m\!=\!a\!+\!2\!+\!|\be|$ 
summand on the left side of~\e_ref{regprp_e1} with the coefficient
\begin{equation}\label{combpf_e1}
\eta^{m-|\al|}\prod_{(i,j)\in S(\al)}\!\!
\bigg(\frac{(-1)^{\be_{i,j}}}{\be_{i,j}!}\cdot
\frac{\eta^{q_i+1-\be_{i,j}}}{(q_i\!+\!1\!-\!\be_{i,j})!}\bigg)
=\frac{\eta^{a+2+\al\cdot q}}{(q\!+\!1)!^{\al}}
\prod_{(i,j)\in S(\al)}\!\!(-1)^{\be_{i,j}}\binom{q_i\!+\!1}{\be_{i,j}},
\footnote{In the first product, the $(i,j)$-factor is defined to be zero 
if $\be_{i,j}\!>\!(q_i\!+\!1)$}
\end{equation}
where
$$\al\cdot q=\sum_{i=1}^{i=k}\al_iq_i \qquad\hbox{and}\qquad
(q\!+\!1)!^{\al}=\prod_{i=1}^{i=k}\big((q_i\!+\!1)!\big)^{\al_i}.$$
Since the number of above choices is
$$\binom{a\!+\!2\!+\!|\be|}{\al,a\!+\!2\!+\!|\be|\!-\!|\al|},$$
it follows that the coefficient of $C_q^{\al}$ on the left-hand side 
of~\e_ref{regprp_e1} is
\begin{equation}\label{combpf_e3}
\frac{\eta^{a+2+\al\cdot q}}{(q\!+\!1)!^{\al}}\sum_{\be\in(\bar\Z^+)^{S(\al)}}\!
\frac{1}{(a\!+\!2\!+\!|\be|)(a\!+\!1\!+\!|\be|)}
\binom{a\!+\!2\!+\!|\be|}{\al,a\!+\!2\!+\!|\be|\!-\!|\al|}
\prod_{(i,j)\in S(\al)}\!\! (-1)^{\be_{i,j}}\binom{q_i\!+\!1}{\be_{i,j}}.
\end{equation}
If $k\!=\!0$ and thus $(\bar\Z^+)^{S(\al)}\!=\!\{0\}$, this expression
reduces to $a!\eta^{a+2}/(a\!+\!2)!$.
By~\e_ref{Rdfn_e}, this is the term on the right-hand side of~\e_ref{regprp_e1}
that does not involve any~$C_q$.
If $k\!\ge\!1$, \e_ref{combpf_e3} becomes
\begin{equation*}\begin{split}
&\frac{\eta^{a+2+\al\cdot q}}{(q\!+\!1)!^{\al}\al!}
\sum_{b=0}^{\i}(-1)^b\binom{\al\cdot q+|\al|}{b}
\frac{1}{(a\!+\!2\!+\!b)(a\!+\!1\!+\!b)}
\frac{(a\!+\!2\!+\!b)!}{(a\!+\!2\!+\!b\!-\!|\al|)!}\\
&\qquad\qquad\qquad\qquad\qquad\qquad
=a!\frac{\eta^{a+2+|q|}}{(a\!+\!|q|\!+\!2)!} ~\times~
\begin{cases} 1,&\hbox{if}~|\al|\!=\!1;\\
0,&\hbox{if}~|\al|\!\ge\!2; \end{cases}
\end{split}\end{equation*}
by Lemma~\ref{comb_l0}.
By~\e_ref{Rdfn_e}, this is also the coefficient of $C_q^{\al}$ on 
right-hand side of~\e_ref{regprp_e1}.\\

\noindent
{\it Proof of~\e_ref{regprp_e2}:}
The coefficient of $C_q^{\al}$ in 
the $m\!=\!a\!+\!|\be|$ summand on the left side of~\e_ref{regprp_e2} 
is again given by the first expression in~\e_ref{combpf_e1}.
Thus, the coefficient of $C_q^{\al}$ on the left-hand side 
of~\e_ref{regprp_e2} is
\begin{equation}\label{combpf_e5}
\frac{\eta^{a+\al\cdot q}}{(q\!+\!1)!^{\al}}\sum_{\be\in(\bar\Z^+)^{S(\al)}}\!
\binom{a\!+\!|\be|}{\al,a\!+\!|\be|\!-\!|\al|}
\prod_{(i,j)\in S(\al)}\!\! (-1)^{\be_{i,j}}\binom{q_i\!+\!1}{\be_{i,j}}.
\end{equation}
If $k\!=\!0$, this expression reduces to $\eta^a$.
If $k\!\ge\!1$, \e_ref{combpf_e5} becomes
$$\frac{\eta^{a+\al\cdot q}}{(q\!+\!1)!^{\al}\al!}
\sum_{b=0}^{\i}(-1)^b\binom{\al\cdot q+|\al|}{b}
\frac{(a\!+\!b)!}{(a\!+\!b\!-\!|\al|)!}
=(-1)^{|\al|}\eta^a ~\times~
\begin{cases} 1,&\hbox{if}~|q|\!=\!0;\\
0,&\hbox{if}~|q|\!\ge\!1; \end{cases}$$
by Lemma~\ref{comb_l0}, as required.\\

\noindent
{\it Proof of Lemma~\ref{prodres_lmm}:} For each $e\!\in\!E$, let
$$r_e=\fR_{\la=0}\big\{f_e(\la)\}, \qquad
g_e(\la)=f_e(\la)-r_e\la^{-1}, \qquad 
h_e(\la)=\la f_e(\la)=r_e+\la g_e(\la).$$
Since $f_e$ has at most a simple pole at $\la\!=\!0$,
$g_e$ and $h_e$ are holomorphic at $\la\!=\!0$ and
$$\cD_\la^{j_e+1}h_e=\cD_\la^{j_e}g_e \qquad\forall~j_e\!\ge\!0.$$
Since $\prod_{e\in E}\!f_e(\la)$ has a pole of order at most $|E|$ at $\la\!=\!0$,
\begin{equation*}\begin{split}
\fR_{\la=0}\bigg\{\prod_{e\in E}\!f_e(\la)\bigg\}
&=\cD_{\la}^{|E|-1}\bigg\{\la^{|E|}\prod_{e\in E}\!f_e\bigg\}
=\cD_{\la}^{|E|-1}\bigg\{\prod_{e\in E}\!h_e\bigg\}
=\sum_{\sum\limits_{e\in E}j_e=|E|-1}\prod_{e\in E}\!\cD_\la^{j_e}h_e\\
&=\sum_{E_+\subset E}\Bigg\{\bigg(\prod_{e\in E_+}\!\!\!r_e\bigg)
\sum_{\sum\limits_{e\not\in E_+}\!\!\!j_e=|E_+|-1} \!\!\!
\bigg(\prod_{e\not\in E_+}\!\!\cD_\la^{j_e+1}h_e \bigg)\Bigg\}\\
&=\sum_{E_+\subset E}\Bigg\{\bigg(\prod_{e\in E_+}\!\!\!r_e\bigg)
\sum_{\sum\limits_{e\not\in E_+}\!\!\!j_e=|E_+|-1} \!\!\!
\bigg(\prod_{e\not\in E_+}\!\!\cD_\la^{j_e}g_e \bigg)\Bigg\}\\
&=\sum_{E_+\subset E}\Bigg\{\bigg(\prod_{e\in E_+}\!\!\!r_e\bigg)
\cD_\la^{|E_+|-1}\bigg(\prod_{e\not\in E_+}\!\!\!g_e\bigg)\Bigg\},
\end{split}\end{equation*}
as claimed.

\section{Comparison of Mirror Symmetry Formulations}
\label{comparison_app}

\noindent
In this appendix we compare a number of mirror symmetry formulations for
genus~$0$ and genus~$1$ curves in a quintic threefold.
In all cases, the predictions are of the form 
$$F_g^{top}(T)=F_g(t),$$
where $g\!=\!0,1$, $F_g^{top}(T)$ is a generating function for the genus~$g$ 
GW-invariants of a quintic threefold 
(related to an $A$-model correlation function), 
$F_g(t)$ is an explicit function of~$t$ 
(related to a $B$-model correlation function), 
and $T\!=\!\J(t)$ for some function~$\J$ (called a mirror transformation).\\

\noindent
In~\cite{CaDGP}, the variables on the $B$-side and $A$-side are $\psi$ and $t$,
respectively.
Let
$$\vp_0(\psi)=1+\sum_{d=1}^{\i}\frac{(5d)!}{(d!)^5}(5\psi)^{-5d};$$
this is equation (3.8) in~\cite{CaDGP} with $n$ replaced by~$d$.
The mirror transformation is defined by equation~(5.9):
$$t=\J(\psi)\equiv
-\frac{5}{2\pi\I}\bigg\{\ln(5\psi)- \frac{1}{\vp_0(\psi)}
\sum_{d=1}^{\i}\frac{(5d)!}{(d!)^5}(5\psi)^{-5d}
\Big(\sum_{r=d+1}^{5d}\frac{1}{r}\Big)\bigg\}.$$
The mirror symmetry prediction for genus~$0$ curves is given in \cite{CaDGP}
in equation~(5.13):
\begin{equation}\label{CaDGP_e}
5+\sum_{d=1}^{\i}n_{0,d}d^3\frac{e^{2\pi\I dt}}{1\!-\!e^{2\pi\I dt}}
=\frac{5\psi^2}{(1\!-\!\psi^5)\vp_0^2(\psi)}
\bigg(\frac{5}{2\pi\I}\Big(\frac{d\psi}{dt}\Big)\bigg)^3,
\end{equation}
where $n_{0,d}$ is the genus~$0$ degree~$d$ instanton number of a quintic.
These numbers are related to the genus~$0$ GW-invariants by the formula
\begin{equation}\label{instaton_e0}
N_{0,d}=\sum_{d_1d_2=d}\frac{n_{0,d_1}}{d_2^3}
\qquad\Llra\qquad
\sum_{d=1}^{\i}N_{0,d}d^3q^d=\sum_{d=1}^{\i}n_{0,d}d^3\frac{q^d}{1\!-\!q^d}.
\end{equation}
The right-hand side of~\e_ref{CaDGP_e} replaces the function $\ka_{ttt}$ appearing 
in \cite[(5.13)]{CaDGP}, using \cite[(5.11)]{CaDGP} and the following two lines.\\

\noindent
In \cite[Chapter~2]{CoKa}, the $B$-side variable is $x$ and 
the variables on $A$-side are $s$ and~$q\!=\!e^s$.
Let 
\begin{equation*}
y_0(x)=1+\sum_{d=1}^{\i}\frac{(5d)!}{(d!)^5}(-x)^d,\qquad
y_1(x)=y_0(x)\ln(-x)+
5\sum_{d=1}^{\i}\frac{(5d)!}{(d!)^5}(-x)^d
\Big(\sum_{r=d+1}^{5d}\frac{1}{r}\Big).
\end{equation*}
These equations are (2.23) and (2.24) in~\cite{CoKa}.
The mirror transformation is given~by
\begin{equation}\label{CoKmap_e}
s=\J(x)\equiv y_1(x)\big/y_0(x) \qquad\hbox{and}\qquad
q=e^s=e^{y_1(x)/y_0(x)}.
\end{equation}
The mirror symmetry prediction for genus~$0$ curves is given in \cite{CoKa}
in equation~(2.26):
\begin{equation}\label{CoK_e}
5+\sum_{d=1}^{\i}n_{0,d}d^3\frac{q^d}{1\!-\!q^d}
=\frac{5}{(1\!+\!5^5x)y_0^2(x)}
\bigg(\frac{q}{x}\frac{dx}{dq}\bigg)^3.
\end{equation}
The relation with variables in~\cite{CaDGP} is
\begin{equation}\label{CaDGPCoK_e}
x=-(5\psi)^{-5}, \qquad s=2\pi\I t, \qquad  q=e^{2\pi\I t}.
\end{equation}
With these identifications, $y_0(x)\!=\!\vp_0(\psi)$, $\J(x)$ of~\cite{CoKa}
is $\J(\psi)$ of \cite{CaDGP} times $2\pi\I$, and the right-hand sides of
\e_ref{CaDGP_e} and~\e_ref{CoK_e} are the same.\\

\noindent
In \cite[Chapters~29,30]{MirSym}, the variables on the $B$ and $A$ sides
are $t$ and~$T$.
The mirror transformation is given~by
$$T=\J(t)\equiv J_1(t)=I_1(t)\big/I_0(t),$$
where $I_q$ and $J_q$ are as in~\e_ref{Ifuncdfn5_e} and~\e_ref{Jfuncdfn5_e}.
The mirror symmetry prediction for genus~$0$ curves is formulated in 
\cite[(29.2)]{MirSym} as
\begin{equation}\label{MirSym_e}
e^{HT}+\frac{H^2}{5}\sum_{d=1}^{\i}n_{0,d}d^3
\sum_{k=1}^{\i}\frac{e^{(H+kd)T}}{(H\!+\!kd)^2}
=\sum_{i=0}^{i=3}J_i(t)H^i \qquad\mod H^4.
\end{equation}
Using~\e_ref{instaton_e0}, \e_ref{MirSym_e} can be re-written as
\begin{equation}\label{MirSym_e2}
\sum_{i=0}^{i=3}\frac{H^iT^i}{i!}+
\frac{H^2}{5}\sum_{d=1}^{\i}N_{0,d}(d\!-\!2H)e^{(H+d)T}
=\sum_{i=0}^{i=3}J_i(t)H^i \qquad\mod H^4.
\end{equation}
The relation~\e_ref{g0ms_e} is obtained by extracting the $H^2$ and $H^3$-terms 
from~\e_ref{MirSym_e2}; it is the statement of \cite[Exercise~29.2.2]{MirSym}, 
minus a typo.
The relation with variables in~\cite{CoKa} is
\begin{equation}\label{CoKMirSym_e}
x=-e^t, \qquad s=T, \quad\hbox{and}\quad q=e^T.
\end{equation}
With these identifications, $I_i(t)\!=\!y_i(x)$ and $\J(t)$ of~\cite{MirSym}
is precisely $\J(x)$ of \cite{CoKa}.
It is shown at the very end of  \cite[Chapter~2]{CoKa} that 
the third derivative of the right-hand side of~\e_ref{g0ms_e} with respect 
to~$s\!=\!T$ is the right-hand side of~\e_ref{CoK_e}.
We note that throughout Section~2.6 of~\cite{CoKa} (in contrast to Section~2.4),
$$Y(q)=\frac{5}{(1\!+\!5^5x)y_0^2(x)}
\bigg(\frac{q}{x}\Big(\frac{dx}{dq}\Big)\bigg)^3.$$
The relation between $x$ and $q$ is described in the previous paragraph.
In Section~2.5, this function~$Y$ is called 
\sf{the normalized Yukawa coupling regarded as a function of~$q$}.\\

\noindent
As in~\cite{CaDGP}, the $B$-side and $A$-side variables in~\cite{BCOV}
are $\psi$ and~$t$.
However, they are now related by the mirror transformation
$$t=\J(\psi)\equiv
5\bigg\{\ln(5\psi)-\frac{1}{\vp_0(\psi)}
\sum_{d=1}^{\i}\frac{(5d)!}{(d!)^5}(5\psi)^{-5d}
\Big(\sum_{r=d+1}^{5d}\frac{1}{r}\Big)\bigg\}=-T,$$
with $T$ as in \cite[Chapters 29,30]{MirSym}.
This is contrary to the suggestion in the paper that $\psi$ and $t$
are related in the same way as in~\cite{CaDGP}.
Let
\begin{equation}\label{BCOV_eB}
F_1(\psi)=\ln\bigg(\Big(\frac{\psi}{\vp_0(\psi)}\Big)^{62/3}
(1\!-\!\psi^5)^{-1/6}\Big(\frac{d\psi}{dt}\Big)\bigg);
\end{equation}
see equation (23) in~\cite{BCOV}.
The mirror symmetry prediction for genus~$1$ curves is given in \cite{BCOV}
in equation~(24):
\begin{equation}\label{BCOV_e}
\frac{25}{6}-
2\!\!\sum_{d_1,d_2=1}^{\i}\!\!n_{1,d_1}d_1d_2\frac{q^{d_1d_2}}{1\!-\!q^{d_1d_2}}
-\frac{1}{6}\sum_{d=1}^{\i}n_{0,d}d\frac{q^d}{1\!-\!q^d}
=\partial_t F_1(\psi),
\end{equation}
where $q\!=\!e^{-t}$ and $n_{1,d}$ is the genus~$1$ degree-$d$ instanton number 
of the quintic.
These numbers are related to the genus~$1$ GW-invariants by the formula
\begin{equation}\label{instaton_e1}\begin{split}
N_{1,d}&=\sum_{d_1d_2=d}n_{1,d_1}\frac{\si_{d_2}}{d_2}
+\frac{1}{12}\sum_{d_1d_2=d}n_{0,d_1}\frac{1}{d_2}\\
&\qquad\Llra\qquad
\sum_{d=1}^{\i}N_{1,d}dq^d
=\!\!\sum_{d_1,d_2=1}^{\i}\!\!n_{1,d_1}d_1d_2\frac{q^{d_1d_2}}{1\!-\!q^{d_1d_2}}
+\frac{1}{12}\sum_{d=1}^{\i}n_{0,d}d\frac{q^d}{1\!-\!q^d},
\end{split}\end{equation}
where $\si_r$ is the number of degree~$r$ unramified connected
covers of a smooth genus~$1$ surface or equivalently of subgroups of $\Z^2$
of index~$r$. Since this number is the same as the sum of positive integer
divisors of~$r$,
\begin{equation}\label{toruscovers_e}
\sum_{r=1}^{\i}\si_rq^r=\sum_{r=1}^{\i}r\frac{q^r}{1\!-\!q^r}.
\end{equation}
This identity implies equivalence of the two equalities in~\e_ref{instaton_e1}.\\

\noindent
Integrating both sides of~\e_ref{BCOV_e} with respect to $t$ and using~\e_ref{instaton_e1},
we find that~\e_ref{BCOV_e} is equivalent~to
\begin{equation}\label{BCOV_e2}
C+\frac{25}{6}t+2\sum_{d=1}^{\i}N_{1,d}q^d
=\ln\bigg(\Big(\frac{\psi}{\vp_0(\psi)}\Big)^{62/3}
(1\!-\!\psi^5)^{-1/6}\Big(\frac{d\psi}{dt}\Big)\bigg)
\end{equation}
for some constant $C$.
This equality should be interpreted by moving $25t/6$ to 
the right-hand side and expanding as a power series in~$q$ at~$q\!=\!0$.
The relation between the variables in~\cite{BCOV} and~\cite{CoKa} is 
$$x=-(5\psi)^{-5}, \qquad s=-t, \quad\hbox{and}\quad q=q.$$
Thus, \e_ref{BCOV_e2} is equivalent~to
\begin{equation}\label{BCOV_e3}
C'-\frac{25}{6}s+2\sum_{d=1}^{\i}N_{1,d}q^d
=\ln\bigg(x^{-25/6}y_0(x)^{-62/3}\big(1\!+\!5^5x\big)^{-1/6}
\Big(\frac{q}{x}\frac{dx}{dq}\Big)\bigg),
\end{equation}
with $q\!=\!e^s$ and $s$ and $x$ related by the mirror transformation~\e_ref{CoKmap_e}.
In the notation of~\cite{MirSym}, i.e.~with identifications~\e_ref{CoKMirSym_e},
\e_ref{BCOV_e3} becomes
\begin{equation}\label{BCOV_e4}
C''-\frac{25}{6}T+2\sum_{d=1}^{\i}N_{1,d}e^{dT}
=-\frac{25}{6}t+
\ln\bigg(I_0(t)^{-62/3}\big(1\!-\!5^5e^t\big)^{-1/6}J_1'(t)^{-1}\bigg).
\end{equation}
It is straightforward to see that $C''\!=\!0$. 
Thus, \e_ref{BCOV_e} and \e_ref{BCOV_e4} are equivalent to~\e_ref{g1ms_e}.\\

\noindent
{\it Remark:} The conventions used in~\cite{KlPa} to formulate a mirror symmetry prediction
for the genus~$1$ GW-invariants of a sextic fourfold  are the same as in~\cite{BCOV},
except $5$ above is replaced by~$6$.\\

\vspace{.2in}

\noindent
{\it Department of Mathematics, SUNY Stony Brook, NY 11794-3651\\
azinger@math.sunysb.edu}

\end{document}